\documentclass[a4paper,11pt,oneside]{book}
\usepackage[english]{babel}
\usepackage{latexsym}
\usepackage{amssymb,amsmath,amstext,amsthm,amsfonts}
\usepackage[latin1]{inputenc}
\usepackage{amsthm}
\usepackage{graphicx}
\usepackage{graphics}

\input pictex.tex

\newcommand{\esp}{\hspace{0.05cm}}
\newcommand{\Q}{\mathbb Q}
\newcommand{\N}{\mathbb N}
\newcommand{\R}{\mathbb R}
\newcommand{\efe}{\mathrm{F}}
\newcommand{\Z}{\mathbb Z}

\theoremstyle{definition}

\newtheorem{thm}{Theorem}[section]

\newtheorem{prop}[thm]{Proposition}

\newtheorem{lem}[thm]{Lemma}

\newtheorem{defn}[thm]{Definition}

\newtheorem{cor}[thm]{Corollary}

\newtheorem{ex}[thm]{Example}

\newtheorem{rem}[thm]{Remark}

\newtheorem{qs}[thm]{Question}

\newcommand{\ce}{\mathcal{C}}
\newcommand{\vsp}{\vspace{0.1cm}}
\newcommand{\vs}{\vspace{0.35cm}}

\begin{document}


\begin{center}
\null \vspace{3cm}


{\Huge {\sc Orderable Groups}}

 \vs\vs\vs

 Tesis

 entregada a la

 Facultad de Ciencias

 de la

 Universidad de Chile

 en cumplimiento parcial de los requisitos

 para optar al grado de

 \vs\vs

 Doctor en Ciencias con mención en Matemáticas

\vs

 Diciembre 2010

\vs\vs\vs

por

\vs\vs

Cristóbal Luciano Rivas Espinosa

\vs\vs

Director de Tesis: Dr. Andrés Ignacio Navas Flores
\end{center}

\thispagestyle{empty}

\newpage
\null \vspace{2cm} $\hfill$ {\Large\em Para Antonia, por todo
\quad\quad}

\newpage

\begin{center}
{\Large{ \textbf{Agradecimientos:}}}
\end{center}

\vs

Quiero aprovechar esta oportunidad para agradecer a todos los que
han contribuído y me han acompañado, en lo que considero han sido
los mejores años de mi (corta) vida. Agradecer al departamento,
profesores, alumnos y funcionarios, por generar un muy grato
ambiente de trabajo. Por ser una constante fuente de inspiración,
energía y alegrías. Agradezco especialmente al profesor Manuel
Pinto, por su suculento apoyo en momentos difíciles. Por su entrega
para con las matemáticas y con las personas. Quiero agradecer
también a Adam Clay y Dale Rolfsen, por su gran hospitalidad en
nuestra estadía en Vancouver, y por numerosas e interesantes
discusiones sobre temas relacionados con esta tesis.

\vsp

Quiero agradecer a mis padres. A él, por su tenaz influencia... por
enseñarme a trabajar duro y con pasión. A ella, por constantemente
traerme los pies sobre la tierra con un simple por qué.

\vsp

Finalmente agradezco fuertemente a Andrés Navas. Por su enorme
generosidad, entusiasmo y dedicación. Para mí, ha sido una
experiencia fascinante poder trabajar con él.

\vsp

Este trabajo a sido parcialmente financiado por la Comición Nacional
de Investigación Científica y Tecnológica (CONICYT) a través de la
beca Apoyo a la Realización de Tesis Doctoral.

\newpage

\begin{center}
{\Large{ \textbf{Resumen:}}}
\end{center}

\vs

En este trabajo estudiamos grupos ordenables. Ponemos especial
énfasis en órdenes de tipo Conrad.

\vsp

En el Capítulo 1 recordamos algunos resultados y definiciones
básicos. También damos una nueva caracterización de la propiedad de
Conrad.

\vsp

En el Capítulo 2 usamos dicha nueva caracterización para obtener una
clasificación de los grupos que admiten solo una cantidad finita de
órdenes Conrad \S 2.1. Con esta clasificación en la mano, somos
capaces de mostrar que el espacio de órdenes Conrad es finito, o
bien no contiene órdenes aislados \S 2.2. Finalmente, la nueva
caracterización de órdenes Conrad nos permite dar un teorema de
estructura para el espacio de órdenes a izquierda, esto tras
analizar la posibilidad de aproximar un orden por sus conjugados \S
2.3.

\vsp

En el Capítulo 3, mostramos que, para grupos que admiten solo una
cantidad finita de órdenes Conrad, es equivalente tener un orden a
izquierda aislado que tener finitos órdenes a izquierda.

\vsp

En el Capítulo 4 probamos que el espacio de órdenes a izquierda del
grupo libre a dos o mas generadores, tiene una órbita densa bajo la
acción natural de éste grupo en dicho espacio. Esto resulta en una
nueva demostración del hecho que el espacio de órdenes a izquierda
del grupo libre en dos o mas generadores no tiene órdenes aislados.

\vsp

En el Capítulo 5 describimos el espacio de bi-órdenes del grupo de
Thompson $\efe$. Mostramos que este espacio está compuesto de 8
puntos aislados junto con 4 copias canónicas del conjunto de Cantor.

\newpage
\begin{center}
{\Large{ \textbf{Abstract:}}}
\end{center}

\vs

In this work we study orderable groups. We put special attention to
Conradian orderings.

\vsp

In Chapter 1 we give the basic background and notations. We also
give a new characterization of the Conrad property for orderings.

\vsp

In Chapter 2, we use the new characterization of the Conradian
property to give a classification of groups admitting finitely many
Conradian orderings \S 2.1. Using this classification we deduce a
structure theorem for the space of Conradian orderings \S 2.2. In
addition, we are able to give a structure theorem for the space of
left-orderings on a group by studying the possibility of
approximating a given ordering by its conjugates \S 2.3.

\vsp

In Chapter 3 we show that, for groups having finitely many Conradian
orderings, having an isolated left-ordering is equivalent to having
only finitely many left-orderings.

\vsp

In Chapter 4, we prove that the space of left-orderings of the free
group on $n\geq2$ generators have a dense orbit under the natural
action of the free group on it. This gives a new proof of the fact
that the space of left-orderings of the free group in at least two
generators have no isolated point.


\vsp

In Chapter 5, we describe the space of bi-orderings of the
Thompson's group $\efe$. We show that this space contains eight
isolated points together with four canonical copies of the Cantor
set.


 \tableofcontents


\chapter[{\small Introduction}]{Introduction}
\pagenumbering{arabic}


A \textit{left-ordered group} $G$ is a group $G$ with a (total)
order relation $\preceq$ that is invariant under left
multiplication. That is, $f\prec g$ implies $hf\prec hg$ for any
$f,g,h$ in $G$. If, in addition, we have that $f\prec g$ implies
$hfh^{-1}\prec hgh^{-1}$ for all $f,g,h$ in $G$, then we say that
$G$ is \textit{bi-ordered} or that $G$ has a bi-invariant ordering.
We will use the term \textit{ordered} when there is no harm of
ambiguity (\textit{e.g.} when $G$ is Abelian).

\vsp

The theory of orderable groups is a venerable subject of mathematics
whose starting point are the works of R. Dedekind and O. H\"older at
the end of XIX century and at the beginning of XX century,
respectively. Dedekind characterizes the real numbers  as a complete
ordered Abelian group, while H\"older proves that any
\textit{Archimedean}\footnote{An ordering is Archimedean if for any
$a\prec b, a\not= id$, there exists $n\in \Z$ such that $b\prec
a^n$.} Abelian ordered group is order isomorphic to a subgroup of
the additive real numbers with the standard ordering; see
\cite{holder} or \cite{ghys} for a modern version of this.

\vsp

Besides the two different kinds of orderings described above, there
is a third type which will be shown to be of great importance in
this work. These are left-orderings satisfying $$f \succ id \text{
and } g \succ id \; \Rightarrow f g^n \succ g \text{ for some } n
\in \mathbb{N}= \{1,2, \ldots \}.$$ These so-called {\em Conradian}
orderings (or $\ce$-orderings) were introduced, in the late fifties,
by P. Conrad in his seminal work \cite{conrad}. There, Conrad shows
that the above condition on a left-ordered group is equivalent to
the fact that the conclusion of H\"older's theorem holds ``locally"
(see $(4)$ below). Since their introduction, Conradian orderings
have played a fundamental role in the theory of left-orderable
groups; see, for instance, \cite{botto, linnell,witte,navas,
rr,zenkov}. Actually, for some time, it was an open question whether
any left-orderable group admits a Conradian ordering. To the best of
our knowledge, the first example of a left-orderable group admitting
no $\ce$-ordering appears in \cite{thurston}, but, apparently, this
was not widely known (among people mostly interested in ordered
groups) until \cite{bergman} appeared.

\vsp

For the statement of Conrad's theorem recall that, in a left-ordered
group $(G,\preceq)$, a subset $S$ is {\em convex} if whenever \esp
$f_1 \prec h \prec f_2$ \esp for some $f_1, f_2$ in $S$, we have $h
\in S$. As it is easy to check, the family of convex subgroups is
linearly ordered under inclusion \cite{glass,kopytov,rr}. In
particular, (arbitrary) unions and intersections of convex subgroups
is also a convex subgroup. Therefore, for every $g\in G$, there
exists $G_g$ (resp. $G^g$), the largest (resp. smaller) convex
subgroup that does not contain $g$ (resp. does contain $g$). The
inclusion $G_g\subset G^g$ is typically referred to as the {\em
$\preceq$-convex jump} associated to $g$.

\begin{thm}[{\bf Conrad}] \label{teo C}{\em A left-ordering $\, \preceq \,$
on a group $G$ is said to be \textit{Conradian} if the following
four equivalent properties hold:

\vsp\vsp

\noindent (1) For all $f \succ id$ and $g \succ id \; $ (for all
{\em positive} $f,g$, for short), we have $f g^n \succ g$ for some
$n \in \mathbb{N}$.

\vsp\vsp

\noindent (2) If $\;1 \prec g  \prec f$, then $g^{-1} f^n g \succ f$
for some $n \in \mathbb{N} $.

\vsp\vsp

\noindent (3) For all positive $g \in G$, the set $S_g = \{f \in G
\mid \esp f^n \prec g, \mbox{ for all } n \in \mathbb{Z} \}$ is a
\textit{convex} subgroup.

\vsp\vsp

\noindent (4) For every $g$, we have that $G_g$ is normal in $G^g$,
and there exists a non-decreasing group homomorphism (to be referred
to as the {\em Conrad homomorphism}) \esp $\tau_{\preceq}^{g} \!:
G^g \rightarrow \mathbb{R}$ \esp whose kernel coincides with $G_g$.
Moreover, this homomorphism is unique up to multiplication by a
positive real number.}

\end{thm}

\vsp

About a decade ago, a new tool for studying left-orderable groups,
the so-called \textit{space of left-orderings} of a left-orderable
group, was introduced by Ghys and, independently, by Sikora
\cite{sikora}. Roughly, the space of left-orderings of a group $G$
is the set of all left-orderings of $G$, where we declare two
left-orderings to be ``close" if they coincide on a large finite
subset of $G$. This object turns out to be a Hausdorff, totally
disconnected and compact topological space on which $G$ acts by
conjugacy: given $\preceq$, a left-ordering on $G$, and $f,g$ in
$G$, we define $\preceq_f$ by $id \preceq_f g$ if and only if
$id\preceq fgf^{-1}$; see \S \ref{pre-1} for details. Although this
object appears for the first time in the literature in
\cite{sikora},  it was in \cite{linnell} and specially \cite{witte},
that the full strength of this object was stressed. In
\cite{linnell}, Linnell put to great use the compactness of the
space of left-orderings to show that if a group admits infinitely
many left-orderings, then it admits uncountably infinitely many. On
the other hand, in \cite{witte}, Morris-Witte squeezes the dynamics
of a group acting on its space of left-orderings, to show that an
amenable, left-orderable group must admit a Conradian ordering.

\vsp

As it was noticed in \cite{leslie,navas}, in (1) and (2) above one
may actually take $n \!=\! 2$. The topological counterpart of this
is the fact that the set of $\ce$-orderings (of a given group) is
compact when it is endowed with the natural topology; see \S
\ref{pre-1}. This leads, for instance, to a new and short proof of
the fact, first proved by Brodskii in \cite{brodskii}, that {\em
locally indicable}\footnote{A group $\Gamma$ is {\em locally
indicable} if for any nontrivial finitely generated subgroup $H$,
there exists a nontrivial group homomorphism from $H$ to the group
of real numbers under addition.} groups are $\ce$-orderable
\cite[Proposition 3.11]{navas}; see also \cite[Corollary
3.2.2]{kopytov}. In particular, the class of $\ce$-orderable groups
contains the class of torsion-free one-relator groups; see
\cite{brodskii}. Note that, from (4) above, the converse to this
result also holds, that is, $\ce$-orderable groups are locally
indicable. Indeed, in a finitely generated, $\ce$-ordered group
$(G,\preceq)$, the homomorphism $\tau^g_\preceq:G\to \R$, where
$g\in G$ is the ``largest" element in the generating set of $G$, is
nontrivial.

\vsp

One last important ingredient, fundamental for our work, is the
dynamical content of the Conradian property for left-orderings
revealed by Navas in \cite{navas}. There, Navas shows that a
left-ordering on a countable group is Conradian if and only if some
natural action on the real line, the so-called {\em dynamical
realization of a left-ordering} (which we trace back to
\cite{ghys}), has no {\em crossings}; see Proposition \ref{real din}
for the definition of the action and \S \ref{sec crossings} for the
definition of crossings.

\vsp

What the concept of crossings is encoding, is the fact that the
action of left translation of $G$ on itself has some sort of
well-behaved ``levels" structure, which puts great constraint to the
dynamics on the cosets of convex subgroups. Nevertheless, as
illustrated in \S \ref{basic contr}, Conradian orderings, unlike
bi-invariant orderings, shares many nice properties with
left-orderings, especially those related to possible modifications.
It is this mixture between rigidity and flexibility what makes
Conradian orderings a good stand point in the study of the more
general left-orderings, and also, what makes them a very nice object
of study.

\newpage


\section{Description of the main results}


Our first main result is the generalization for uncountable groups
of Navas dynamical characterization of Conradian orderings. To do
this, we had to understand the concept of crossings in an intrinsic
way. For us, a crossing for an action by order preserving bijections
of a group $G$ on a totally ordered space $(\Omega,\leq)$, is a
5-uple $(f,g,u,v,w)$, where $f,g$ (resp. $u,v,w$) belong to $G$
(resp. $\Omega$), that satisfies:

\vsp

\noindent $(i)$ \esp $u < w < v$.

\vsp

\noindent $(ii)$ \esp For every $\,n \in \mathbb{N} , \,$ we have
$\,g^n u < v $ \esp and \esp $f^n v > u \,.$

\vsp

\noindent $(iii)$ \esp There exist $M,N$ in $\mathbb{N}$ such that
\esp $f^N v < w < g^M u$.

\vsp

In \S \ref{sec crossings} we show

\vsp\vsp

\noindent {\bf Theorem A.} {\em A left-ordering $\preceq$ on a group
$G$ is Conradian if and only if the action by left translation on
itself contains no crossings.}

\vsp\vsp

We point out that, besides the four equivalences of the Conrad
property given in Theorem \ref{teo C}, many more can be found in
\cite[\S 7.4]{botto}. Unlike ours, all of them are algebraic
descriptions. The investigation of the consequences of this new
characterization of Conradian orderings concerns almost two third of
our work.

\vsp

A major consequence of Theorem A is the classification of groups
admitting only finitely many Conradian orderings. This can be though
of as an analogue of Tararin's classification of groups admitting
only finitely many left-orderings \cite[Theorem 5.2.1]{kopytov}. For
the statement of both results, recall that a series
$$\{ id \} = G_0 \lhd G_{1} \lhd \ldots \lhd G_{n-1} \lhd G_n = G$$
is said to be {\em rational} if it is subnormal ({\em i.e.,} each
$G_{i}$ is normal in $G_{i+1}$) and each quotient $G_{i+1} / G_{i}$
is torsion-free rank-one Abelian. The series is called {\em normal}
if, in addition, each $G_i$ is normal in $G$. In \S \ref{finitos
C-ord} we show

\vsp\vsp

\noindent {\bf Theorem B.}  {\em Let $G$ be a $\ce$-orderable group.
If $G$ admits only finitely many $\ce$-orderings, then $G$ admits a
unique (hence normal) rational series. In this series, no quotient
$G_{i+2}/G_{i}$ is Abelian. Conversely, if $G$ is a group admitting
a  normal rational series
$$\{ id \} = G_0 \lhd G_{1} \lhd \ldots \lhd G_{n-1} \lhd G_n = G$$
so that no quotient $G_{i+2} / G_{i}$ is Abelian, then the number of
$\ce$-orderings on $G$ equals $2^n$.}

\vsp\vsp

We state Tararin's classification as

\begin{thm}[{\bf Tararin}] \label{teo T} {\em Let $G$ be a left-orderable
group. If $G$ admits only finitely many left-orderings, then $G$
admits a unique (hence normal) rational series. In this series, no
quotient $G_{i+2}/G_{i}$ is bi-orderable. Conversely, if $G$ is a
group admitting a normal rational series
$$\{ id \} = G_0 \lhd G_{1} \lhd \ldots \lhd G_{n-1} \lhd G_n = G$$
so that no quotient $G_{i+2} / G_{i}$ is bi-orderable, then the
number of left-orderings on $G$ equals $2^n$.}
\end{thm}

Note that the statement of Tararin's theorem is the same as that of
Theorem B though changing ``$\ce$-orderings" by ``left-orderings",
and the condition ``$G_{i+2} / G_{i}$ non Abelian" by ``$G_{i+2} /
G_{i}$ non bi-orderable".

\vsp

In the late nineties, using Tararin's classification, Zenkov was
able to deduce that if a locally indicable group admits infinitely
many left-orderings, then it admits uncountably many of them; see
\cite[Theorem 5.2.5]{kopytov} or \cite{zenkov}. Here, in \S \ref{sec
struct thm for C}, we use our classification of groups admitting
only finitely many $\ce$-orderings to show

\vsp\vsp

\noindent {\bf Theorem C.} {\em Let $G$ be a $\ce$-orderable group.
If $G$ admits infinitely many $\ce$-orderings, then it admits
uncountably many $\ce$-orderings. Moreover, none of these
$\ce$-orderings is isolated in the space of $\ce$-orderings.}

\vsp\vsp

We remark that the second statement of Theorem C is much stronger
than the first one. For instance, if $G$ is countable, then its
space of $\ce$-orderings is either finite or a Cantor set. Moreover,
as it will be exemplified below, the absence of isolated
$\ce$-orderings when there are infinitely many of them, is a
behavior not shared with left-orderings nor with bi-orderings.
Actually, knowing when a given left-orderable group admits an
isolated left-ordering is one of the main open problems in this
theory.

\vsp

Theorem C corroborates a general principle concerning
$\ce$-orderings. On the one hand, these are sufficiently rigid in
that they allow deducing structure theorems for the underlying group
({\em e.g.}, local indicability). However, they are still
sufficiently malleable in that, starting with a $\ce$-ordering on a
group, one may create very many $\ce$-orderings, which turn out to
be different from the original one with the only exception of the
pathological cases described in Theorem B.

\vsp

Motivated by the high regularity of Conradian orderings, in Chapter
3, we study the space of left-orderings of groups admitting only
finitely many $\ce$-orderings. This chapter is motivated by
\cite{rivas}, where an explicit description of the space of
left-orderings of the Baumslag-Solitar group $B(1,2)$ -a group with
only $2^2$ $\ce$-orderings, but infinitely many left-orderings- is
made. In \S \ref{sec metabelian case} we use the machinery
developed/exposed in \cite[\S 2]{book}, to extend the argument of
\cite{rivas}, and give an explicit description of the space of
left-orderings of any group admitting only four $\ce$-orderings. We
show that any group $G$ admitting only four $\ce$-orderings but
infinitely left-orderings can be embedded in the (real) affine
group, and that any left-ordering  of $G$ is an induced orderings
(in the sense of \S \ref{sec dynanical versions}) of this affine
action, or one of the four possible Conradian orderings; see Theorem
\ref{laprop}. Once the case ``$n=2$" is solved, a simple induction
argument shows

\vsp\vsp

\noindent {\bf Theorem D.} {\em If $G$ is a $\ce$-orderable group
admitting only finitely many $\ce$-orderings, then its space of
left-orderings is either finite or homeomorphic to the Cantor set.}

\vsp\vsp

\vsp

As it was already mentioned, in the bi-ordered case the picture is
totally different. At the time of this writing, there is no
classification of groups admitting only finitely many bi-orderings.
Actually the range of groups admitting only finitely many
bi-orderings should be very large, and just a few results give
partial descriptions of this situation; see for instance
\cite[Chapter VI]{botto} and \cite[\S 5.3]{kopytov}. Indeed, this
class contain all the groups fitting in Theorem B that are not in
Tararin's classification, but also a lot of groups algebraically
very different from those. For example, the commutator subgroup of
the group of piecewise affine homeomorphisms of the unit interval,
and many other similar groups (such as the commutator subgroup of
Thompson's group $\efe$), have finitely many bi-orderings; see
\cite{dlab,ZM} and the remark at the end of Chapter 5. In addition,
there are examples of bi-orderable groups admitting infinitely but
(only) countably many bi-orderings; see \cite[Chapter VI]{botto} or
\cite{butts} for an example of a group admitting only countably
infinitely many bi-orderings. For an example of a family of solvable
groups admitting only finitely many bi-orderings see \cite[\S
5.3]{kopytov}. This two ``strange" behaviors are mainly caused by
the strong rigidity of bi-invariant orderings, and our methods seem
not well adapted to investigate this situation.

\vsp

As mentioned earlier, for the case of left-orderings we have
Tararin's classification of groups admitting only finitely many
left-orderings, and also Linnell's result \cite{linnell} saying
that, if infinite, then the number of left-orderings admitted by a
group must be uncountable. In \S \ref{sec a struc thm for left ord},
we give a new proof of Linnell's result, but using a quite different
approach. Our method relies strongly on the nature of Conradian
orderings described in Theorem A. This new characterization of
Conradian orderings is used to detect the so-called {\em Conradian
soul of a left-ordering $\preceq$}, which was introduced by Navas in
\cite{navas} as the maximal convex subgroup for which $\preceq$ is
Conradian; see \S \ref{detecting the soul}. The Conradian soul plays
a fundamental role when dealing with the problem of approximating a
given left-ordering by its conjugates (in the sense of \S
\ref{acting on the space}). We show that, in most cases, this can be
done; see for instance Theorem \ref{trivial soul}. In the few cases
this can not be done, we show that we still have enough information
to conclude

\vsp\vsp

\noindent {\bf Theorem E (Linnell).}  {\em If a left-orderable group
admits infinitely many left-orderings, then it admits uncountably
many left-orderings.}

\vsp\vsp

To finish the discussion, we have to point out that the space of
left-orderings of a group may be infinite and still have isolated
left-orderings. Actually, this is the case of braid groups,
\cite{braids,dd}, and a particular central extension of Hecke groups
\cite{navas-hecke}.

\vs

In the last two chapters of our dissertation, we analyze the spaces
of orderings of two remarkable groups.

\vsp

Chapter 4 is devoted to the study of the space of left-orderings of
the free group on two or more generators, $F_n$, $n\geq 2$. This has
a long history. In \cite{mccleary}, McCleary studies an object
introduced by Conrad in \cite{conrad2}, called the free
lattice-ordered group (in this case) over the free group, which is
an universal object in the class of lattice-ordered groups\footnote{
A lattice-ordered group $(G,\preceq)$, is a group together with a
partial bi-invariant ordering of $G$ satisfying that for all $f,g$
in $G$ there exist $f\vee g \in G$ (resp. $f\wedge g \in G$), the
least upper (resp. greatest lower) bound of $f$ and $g$; see for
instance \cite{glass}.}. He is able to prove that no left-ordering
on the free group on two or more generators is {\em finitely
determined}. In our language, this is equivalent to saying that the
space of left-orderings of the free group on two or more generators
has no isolated points. In \cite{navas}, Navas gives a different and
easier proof of this fact. He shows that small perturbations of the
dynamical realization of a left-ordering of $F_n$, made outside
large compact intervals in $\R$, can be used to approximate the
given left-ordering.

\vsp

In \cite{clay 2}, Clay establishes a strong connection between some
representations of the free lattice-ordered group over a group $G$,
and the dynamics of the action of $G$ on its space of
left-orderings. Using this connection, together with the previous
work of Kopitov \cite{kopytov2}, he showed

\vsp\vsp

\noindent {\bf Theorem F (Clay).} {\em The space of left-orderings
of the free group on two or more generators has a dense orbit under
the natural conjugacy action of $F_n$.}

\vsp\vsp

However, Clay's proof is highly non constructive. Moreover,
Kopytov's work \cite{kopytov2} also involves the free
lattice-ordered group over the free group. In Chapter 4 of this
work, we give an explicit construction of a left-ordering of $F_n$,
whose set of conjugates is dense in the space of left-orderings of
$F_n$. Our proof uses a very simple idea which resembles a lot
McCleary's and Kopitov's originals constructions from
\cite{mccleary} and \cite{kopytov2}, respectively. Nevertheless, as
we avoid the use of of any lattice structure, we don't have to take
care of certain unpleasant technical details which make
\cite{mccleary} and \cite{kopytov2} hard to read. Thus, our
construction is easier to follow.

\vsp

The rough idea for proving Theorem F is the following. Since the
space of left-orderings of $F_n$ is compact, it contains a dense
countable subset. Now, we can consider the dynamical realization of
each of these left-orderings, and cut large pieces of each one of
them. Since we are working with a free group, we can glue these
pieces of dynamical realizations together in a sole action of $F_n$
on the real line. Moreover, if the gluing is made with a little bit
of care, then we can ensure very nice conjugacy properties, from
which we can deduce Theorem F.

\vsp

Theorem F gives a new proof of the fact that no left-ordering of
$F_n$ is isolated. Indeed, as also shown by Clay, any group $G$
admitting a dense orbit in the natural action of $G$ on its space of
left-orderings has no isolated left-orderings. This is shown here in
Proposition \ref{clay -> mccleary}.

\vs

Finally, in Chapter 5, we study the space of bi-orderings of the
(remarkable!) Thompson's group $\efe$. This group is usually
represented as the group of piecewise affine homeomorphisms of the
closed interval, such that the break points of any element in $\efe$
are diadic rationals, there are only finitely many of them, and the
slopes of the elements are integer powers of 2. The standard
reference on $\efe$ is \cite{CFP}.

\vsp

Using the theory of Conradian orderings together with the internal
structure of Thompson's group $\efe$, we are able to show

\vsp\vsp

\noindent {\bf Theorem G.} {\em The space of bi-orderings of the
Thompson group $\efe$ is made up of eight isolated points, together
with four canonical copies of the Cantor set.}

\vsp\vsp

An important intermediate step in proving Theorem G, is the
description of all bi-orderings on $[\efe,\efe]$, the commutator
subgroup of $\efe$. We show that $[\efe,\efe]$ admits only four
bi-orderings. All of them can be easily described. This particular
result is strongly related to \cite{dlab}, where Dlab shows that a
large family of piecewise affine groups, where the slopes of the
elements are contained in a rank-one Abelian group of the
multiplicative groups of positive real numbers, have only four
bi-orderings. However, the commutator subgroup of the Thompson group
is not included in this family, since in \cite{dlab}, the ``break
points" of the elements considered can accumulate on the right. In
particular, elements from \cite{dlab} can have infinitely many break
points. By counterpart, in Thompson's group $\efe$, the break points
are contained in the diadic numbers and, for each element, there are
only finitely many of them. This additional assumption on the brake
points implies for instance that conjugacy classes on $\efe$ are way
smaller than in the groups considered by Dlab. It also implies that
$\efe$ is not a lattice-ordered group under the natural pointwise
partial ordering defined by $f\succeq id$ if and only if $f(x)\geq
x$ for all $x\in [0,1]$. In turns, the groups considered by Dlab are
lattice-ordered with the point wise ordering; see \cite{ZM} for a
discussion on that subject (see also \cite{Me}). As a consequence,
our method of proof is different form the one used in \cite{dlab},
and actually, the method for describing the bi-orderings on
$[\efe,\efe]$ is essentially the same as the one used to describe
the bi-orderings on $\efe$.

\vsp

We point out that in Chapter 5, besides the description of the
topology of the space of bi-orderings of $\efe$, we give an explicit
description of all the bi-orderings on $\efe$. The four canonical
copies of the Cantor set arise form extending the bi-orderings on
$\efe/[\efe,\efe]\simeq\Z^2$ with each of the four bi-orderings on
$[F,F]$. The eight isolated bi-orderings are explicitly described
too. For instance, the left-ordering on $\efe$, whose set of
positive elements are the elements whose first non trivial slope is
greater than 1, determines a bi-invariant ordering, which is
isolated in the space of bi-orderings. The other seven isolated
bi-orderings are similar to this one, but we do not have space to
describe them here.

\vs

We have compiled the main results of our work. The rest of Chapter 1
is devoted to give the basic background and notation. In \S 1.2 we
illustrate some basic constructions for producing new orderings
staring with a given one. In \S \ref{pre-1} we recall the concept of
spaces of orderings. The main dynamical tools for carrying out our
study are recalled in \S 1.4. More importantly, also in \S 1.4, we
develop the concept of \textit{crossings} for an action of a group
on an ordered set, and we prove Theorem A.


\section{Some basic constructions for producing new orderings}

\label{basic contr}


In this section we describe some basic constructions for creating
new (left, C, or bi) orderings starting with a given one. The main
idea is to exploit the flexibility given by the convex subgroups.

\vsp

If $C$ is a proper convex subgroup of a left-ordered group
$(G,\preceq)$, then $\preceq$ induces a total order on the set of
left-cosets of $C$ by
\begin{equation} \label{def convex
ord}g_1C\prec g_2C \Leftrightarrow g_1c_1\prec g_2c_2\;\; \text{for
all $c_1, c_2$ in C}.\end{equation} More importantly, this order is
preserved by the left action of $G$; see for instance \cite[\S
2]{kopytov}. This easily implies

\begin{prop} \label{nice prop} {\em Let $(G,\preceq)$ be a left-ordered (resp. $\ce$-ordered) group and let
$C$ be a convex subgroup. Then any left-ordering (resp.
$\ce$-ordering) on $C$ may be extended, via $\preceq$, to a total
left-ordering (resp. $\ce$-ordering) on $G$. In addition, in this
new left-ordering, $C$ is still a convex subgroup.}
\end{prop}

\noindent{\em Proof:} We denote $\preceq_1$ the induced ordering on
the cosets of $C$ (from equation (\ref{def convex ord})). Let
$\preceq_2$ be any left-ordering on $C$. For $g\in G$, we define $id
\preceq^\prime g$ if and only if $g\in C$ and $id \preceq_2 g$ or
$g\notin C$ and $C\prec_1 gC$. We claim that $\preceq^\prime$ is a
left-ordering.

\vsp

Indeed, let $f,h$ in $G$ such that $f\succ^\prime id$ and
$h\succ^\prime id$. If both $f,h$ belong to $C$, then clearly $fh\in
C$ and $fh \succ^\prime id$. If neither $f$ nor $h$ belong to $C$,
then, since the $G$ action on the cosets of $C$ preserves
$\preceq_1$, we have that $fhC\succ_1 C$, thus $fh\succ^\prime id$.
Finally, if $h\in C$ and $f\notin C$, we have that $fC=fhC\succ_1
C$, so $fh\succ^\prime id$. To check that $hf\succ^\prime id$ we
note that $fC\succ_1 C$ implies $hfC\succ_1 hC=C$. This shows the
left-invariance of $\preceq^\prime$. To see that $C$ is convex in
$\preceq^\prime$, we note that $id\prec^\prime h\prec^\prime c$, for
$c\in C$, is equivalent to $id \prec^\prime h^{-1}c$. We claim that
in this case, $h$ belongs to $C$. We have two possibilities. Either
$h^{-1}c\in C$, in which case we conclude $h\in C$, or $h^{-1}c
\notin C$, in which case we have that $h^{-1}cC=h^{-1}C \succ_1 C$,
therefore, $h^{-1}\succ^\prime id$, which contradicts the fact that
$id\prec^\prime h$. This shows the convexity of $C$.

\vsp

We now show that $\preceq^\prime$ is Conradian when $\preceq$ and
$\preceq_2$ are Conradian. Let $f$, $g$ in $G$, be such that
$id\prec^\prime f\prec^\prime g$. We have to show that
$fg^2\succ^\prime g$ (note that $id \prec f\prec g$ easily implies
$gf^2\succ f$ in any left-ordering!), or, equivalently, that
$g^{-1}fg^2\succ^\prime id$. If both $f$ and $g$ belongs to $C$,
then the conclusion follows by the assumption on $\preceq_2$. If it
is the case that $g$ does not belong to $C$, then we claim that
$g^{-1}fg^2\notin C$. Indeed, since $\preceq$ is Conradian, the
Conrad homomorphism $\tau_\preceq^g$, defined in Theorem \ref{teo
C}, is an order preserving homomorphism whose kernel is $G_g$, and
in this case $C\subseteq G_g$. Therefore
$\tau_\preceq^g(g^{-1}fg^2)=\tau_\preceq^g(fg)>0$. In particular
$g^{-1}fg^2\notin C$. This latter statement, together with the fact
that $\preceq$ is Conradian, implies $g^{-1}fg^2 C\succ_1 C$ which,
in turns, implies $g^{-1}fg^2\succ^\prime id$. $\hfill\square$

\vsp\vsp

\begin{ex} \label{first ex} Let $\, \preceq \,$ be a left-ordering on $G$.
Recall that the {\em reverse} (or ``flipped") ordering, denoted
$\overline{\preceq}$, is the ordering that satisfies
$f\overline{\prec}g \Leftrightarrow f\succ g$. Showing that if
$\preceq$ is a left, $\ce$, or bi- ordering then
$\overline{\preceq}$ is of the same kind is routine. Moreover, the
convex series in $\preceq$ coincides with the convex series in
$\overline{\preceq}$. Now, suppose there is a nontrivial convex
subgroup $C$ of $G$. Then, by Proposition \ref{nice prop}, there is
a (left- or $\ce$-) ordering $\preceq_{C}$ of $G$ defined by $id
\prec_{C} f$, where $f\in G$, if and only if either

\vsp

\noindent -- \esp  $f\succ id$ and $f  \not \in C$, or

\vsp

\noindent -- \esp $f\overline{\succ} id$ and $f\in C$.

\vsp

In the case $\preceq$ is a bi-ordering, the preceding construction
does not imply that $\preceq_{C}$ is also a bi-ordering.
Nevertheless, if $C$ is a convex and normal subgroup, then the
conclusion follows. Indeed, for $f,g$ in $G$ with  $id \prec_{C} f$,
we have to prove that $id\prec_{C} gfg^{-1}$. If $f\notin C$, then
by the normality of $C$ we have that $gfg^{-1}\notin C$ and the
conclusion follows from the bi-invariance of $\preceq$. If $f\in C$,
then we have that $id \succ f$ and $id\succ gfg^{-1} \in C$, which
is the same to say that $id\prec_{C} gfg^{-1}$, so $\preceq_{C}$ is
a bi-ordering.

\end{ex}

The following example will serve us to approximate (in the sense of
\S\ref{pre-1}) a given ordering when the series of convex subgroups
is long enough.

\begin{ex} \label{exflip} Let $g\in G\setminus\{id\}$ and let $\preceq$ be a left-ordering on $G$.
Consider the (perhaps infinite) series of $\, \preceq$-convex
subgroups.
$$\{id\} = G^{id} \subset \ldots \subset G_g \subset G^g \subset\ldots \subset G.$$
We will use Example \ref{first ex} to produce the a new ordering
$\preceq^g$ by ``flipping" the ordering on $G^g\setminus G_g$. More
precisely, we let $\preceq^g=(\preceq_{G^g})_{G_g}$, that is, we
flip the ordering on $G^g$ and then we flip again on $G_g$. It is
easy to see that, for $f\in G$, $id\prec^g f$ if and only if

\vsp

\noindent -- \esp  $f\succ id$ and $f  \not \in G^g$,

\vsp

\noindent -- \esp $f\overline {\succ} id$ and $f\in G^g\setminus
G_g$,

\vsp

\noindent -- \esp  $f\succ id$ and $f  \not \in G_g$.

\vsp

\noindent Clearly, $\preceq^g$ is Conradian when $\preceq$ is
Conradian.

\end{ex}

\vsp

Notice that if $C$ is normal in $G$, equation (\ref{def convex ord})
defines a left-ordering on the group $G/C$. In this case, we have
even more flexibility for producing new (left-, $\ce$- or bi-)
orderings, since we can change our ordering not only on the subgroup
$C$, but also on the quotient group $G/C$. We state this as

\begin{lem} \label{lema 1.2.4} {\em Let $(G,\preceq)$ be a
left-ordered group. Let $C$ be a normal and convex subgroup of $G$.
We denote by $\preceq_1$ the induced ordering on $G/C$ (from
equation (\ref{def convex ord})) and by $\preceq_2$ the restriction
to $C$ of $\preceq$. We have,
 \vsp

\noindent $(i)$ \esp For $f\in G$, $id\prec f$ if and only if $f
\not \in C$ and $f\, C \succ_1 C$, or \esp $f\in C$ and $f\succ_2
id$.

\vsp

\noindent $(ii)$ Let $f\in G$. For any left-ordering
$\preceq^{\prime\prime}$ on $C$, and any left-ordering
$\preceq^{\prime}$ on $G/C$, there is a left-ordering
$\tilde{\preceq}$ on $G$ defined by $f\, \tilde{\succ}\, id$ if and
only if $f  \not \in C$ and $f\, C \succ^\prime C$ , or \esp $f\in
C$ and $f\succ^{\prime\prime} id$.

\vsp

\noindent $(iii)$ $\preceq$ is Conradian if and only if $\preceq_1$
and $\preceq_2$ are Conradian.

\vsp

\noindent $(iv)$ $\preceq$ is a bi-ordering if and only if
$\preceq_1$ and $\preceq_2$ are bi-orderings and  for $c\in C$,
$id\prec_2 c \Rightarrow id\prec_2 fcf^{-1}$ for all $f\in G$.}
\end{lem}

\noindent {\em Proof:} Items $(i)$, $(ii)$ and $(iii)$ follow
arguing as in the proof of Proposition \ref{nice prop}.

\vsp

We show item $(iv)$. Clearly, if $\preceq$ is a bi-ordering, then
$\preceq_2$ is invariant under the whole group $G$. That is, $c\in
C$ and $id \prec_2 c$ implies $fcf^{-1}\in C$ and $id\prec_2
fcf^{-1}$ for all $f\in G$. To see that $\preceq_1$ is bi-invariant
we note that, if not, then there are $f,g$ in $G$ such that $C
\preceq_1 fC$ and $gC\, fC\, g^{-1}C= gfg^{-1}C\prec_1 C$. In
particular, $gfg^{-1}\notin C$, and item $(i)$ implies
$gfg^{-1}\prec id$. This contradiction implies that $\preceq_1$ is a
bi-ordering.

\vsp

For the converse, suppose that $\preceq$ is not a bi-ordering. Then
there are $f$ and $g$ in $G$ such that $id \prec f $ and
$gfg^{-1}\prec id$. If $f\in C$ we have that $\preceq_2$ is not
invariant under the action of $g$, which contradicts the assumption
on $\preceq_2$. If $f\notin C$ then, since $C$ is a normal subgroup,
$gfg^{-1}\notin C$. But, by item $(i)$, in this case we have that
$C\prec_1 fC$ and $gfg^{-1} C\prec_1 C$, so $\preceq_1$ is not
bi-invariant, contrary to our assumption on $\preceq_2$.
$\hfill\square$

\vsp

What follows is a rewording of the previous lemma.

\begin{cor} \label{nice coro} {\em Suppose that $G$ is a left-orderable
(resp. $\ce$-orderable) group and $C$ is a normal subgroup of $G$.
Then for any left-ordering (resp. $\ce$-ordering) $\preceq_1$ on
$G/C$ and any left-ordering (resp. $\ce$-ordering) $\preceq_2$ on
$C$, there is a left-ordering (resp. $\ce$-ordering) $\preceq$ on
$G$ such that $\preceq$ coincides with $\preceq_2$ on $C$ and the
induced ordering of $\preceq$ on $G/C$ coincides with $\preceq_1$.

\vsp

If, in addition, $G$ is bi-orderable, then for any bi-ordering
$\preceq_1$ of $G/C$ and any bi-ordering $\preceq_2$ of $C$ with the
additional property that, for $c\in C$ and $f\in G$ we have
$id\prec_2 c\Rightarrow id \prec_2 fcf^{-1}$, there is a bi-ordering
$\preceq$ on $G$ such that $\preceq$ coincides with $\preceq_2$ on
$C$ and the induced ordering of $\preceq$ on $G/C$ coincides with
$\preceq_1$.}
\end{cor}


\section{The space of orderings of a group}


\label{pre-1}

Recall that, given a left-ordering $\preceq$ on a group $G$, we say
that $f\in G$ is \textit{positive} or {\em $\preceq$-positive}
(resp. \textit{negative} or {\em $\preceq$-negative}) if $f\succ id$
(resp. $f\prec id$). We denote $P_\preceq$ the set of
$\preceq$-positive elements in $G$, and we sometimes call it {\em
the positive cone} of $\preceq$. Clearly, $P_\preceq$ satisfies the
following properties:

\vsp

\noindent $(i)$ $P_\preceq  P_\preceq \subseteq P_\preceq \;$, that
is, $P_\preceq$ is a semi-group;

\vsp

\noindent $(ii)$ $G=P_\preceq \sqcup P^{-1}_\preceq \sqcup \{id\}$,
where the union is disjoint, and $P^{-1}_\preceq=\{g^{-1}\in G\mid
g\in P_\preceq\}=\{g\in G\mid g\prec id\}$.

\vsp

Moreover, given any subset $P\subseteq G$ satisfying the conditions
$(i)$ and $(ii)$ above, we can define a left-ordering $\preceq_P$ by
$f\prec_P g$ if and only if $f^{-1}g\in P$. We will usually identify
$\preceq$ with its positive cone $P_\preceq$.

\vsp

Given a left-orderable group $G$ (of arbitrary cardinality), we
denote the set of all left-orderings on $G$ by $\mathcal{LO} (G)$.
This set has a natural topology first introduced by Sikora in
\cite{sikora}. This topology can be defined by identifying $P\in
\mathcal{LO}(G)$ with its characteristic function $\chi_P \in
\{0,1\}^G $. In this way, we can view $\mathcal{LO}(G)$ embedded in
$\{0,1\}^G$. This latter space, with the product topology, is a
Hausdorff, totally disconnected, and compact space. It is not hard
to see that (the image of) $\mathcal{LO}(G)$ is closed inside, and
hence compact as well (see \cite{navas,sikora} for details).

\vsp

A basis of neighborhoods of $\,\preceq \,$ in $\, \mathcal{LO} (G)$
is the family of the sets $\, U_{g_1,\ldots,g_k}$ of all
left-orderings $\,\preceq'\,$ on $G$ that coincide with $\, \preceq
\,$ on $\{g_1,\ldots,g_k\}$, where $\{g_1,\ldots,g_k\}$ runs over
all finite subsets of $G$. Another basis of neighborhoods is given
by the sets $\, V_{f_1,\ldots,f_k} \,$ of all left-orderings $\,
\preceq' \,$ on $G$ such that all the $\,f_i\,$ are
$\,\preceq'$-positive, where $\{f_1,\ldots,f_k\}$ runs over all
finite subsets of $\,\preceq$-positive elements of $G$. The (perhaps
empty) subspaces $\, \mathcal{BO}(G) \,$ and $\, \mathcal{CO}(G) \,$
of bi-orderings and $\ce$-orderings on $G$ respectively, are closed
inside $\mathcal{LO}(G)$, hence compact; see \cite{navas}.

\vsp

If $G$ is countable, then this topology is metrizable: given an
exhaustion $G_0 \subset G_1 \subset \ldots$ of $G$ by finite sets,
for different $\, \preceq \,$ and $\,\preceq'\,, \,$ we may define
$\,dist(\preceq,\preceq') = 1 / 2^n$, where $n$ is the first integer
such that $\, \preceq \,$ and $\, \preceq'\,$ do not coincide on
$G_n$. If $G$ is finitely generated, we may take $G_n$ as the ball
of radius $n$ with respect to a fixed finite system of generators.

\begin{ex} \label{LOZ^2} It was shown in \cite{sikora} that the space of (bi-) orderings
on a torsion free Abelian group of rank greater than one is
homeomorphic to the Cantor set. We now describe the space of
orderings of $\Z^2$.

\vsp

Let $e_1$, $e_2$ be the standard basis on $\Z^2$. For every $x\in
\overline{\R}=\R\cup \{\infty\}$ define $\psi_x: \Z^2\to \R$ to be
the homomorphism defined by $\psi_\infty(e_1)=0$,
$\psi_\infty(e_2)=1$ and $\psi_x(e_1)=1$, $\psi_x(e_2)=x$ if $x\in
\R$.

\vsp

Clearly, if $x$ is irrational, then $\psi_x$ is injective and
$P_x=\{g\in \Z^2\mid \psi_x(g)> 0\}$ defines a positive cone in
$\Z^2$. The associated ordering is said to be of {\em irrational
type}. Note that this ordering has no proper convex subgroup. These
irrational orderings are dense in $\mathcal{LO}(\Z^2)$.

\vsp

If $x$ is rational or $x=\infty$, then $\psi_x$ is not injective and
$ker(\psi_x)\simeq\Z$. Thus, the set $P_x=\{g\in \Z^2\mid \psi_x(g)>
0\}$ defines only a partial ordering which can be completed into two
total orderings $P_x^+$ and $P_x^-$, where $P_x^+$ (resp. $P_x^-$)
corresponds to the limit of $P_{x_n}$ where $(x_n)$ is a sequence of
irrational numbers converging to $x$ from the right (resp. left).
These orderings are called of {\em rational type} ({\em e.g,} a
lexicographic ordering). In these orderings, $ker(\psi_x)$ is the
unique proper convex subgroup.

\vsp

Finally, we show that any ordering on $\Z^2$ is one of the orderings
just described. Let $\preceq$ be any ordering on $\Z^2$. Since
$\Z^2$ is finitely generated and Abelian, $\preceq$ is Conradian,
and, for $g=\max_\preceq\{{\pm }e_1, \pm e_2\}$, there exists
$\tau=\tau_\preceq^g:\Z^2\to \R$ (as defined in Theorem \ref{teo
C}). Now, let
$$y=  \left\{ \begin{array}{c c} \tau(e_2)/\tau(e_1) & \text{ if }
\tau(e_1)\not=0,
 \\  \infty   & \text{ if } \tau(e_1)=0.
\end{array} \right.$$
Then there is a positive real number $\alpha$ such that
$\alpha\tau=\psi_y$. This shows that $\preceq$ must coincide with
$P_y$ or $P_y^\pm$.

\vsp

There is another, more geometric, way to see the orderings on
$\Z^2$. Since the orderings of rational type are limits of orderings
of irrational type, we just describe the latter type of orderings.
Let $x\in \R\setminus \Q$. Consider $\Z^2$ embedded in $\R^2$ in the
usual way. The unique $\R$-linear extension of $\psi_x$ from $\Z^2$
to $\R^2$ will be denoted $\hat{\psi}_x$. Let $L_x=\{ w\in \R^2\mid
\hat{\psi}_x(w)=0\}$. Take  $w_0\in \hat{\psi}_x^{-1}(1)$ and let
$H_x$ be the open half plane with boundary $L_x$ that contains
$w_0$. Then we have that $P_x=\Z^2\cap H_x$. Moreover, if $\preceq$
is the ordering on $\Z^2$ corresponding to $P_x$, then for $g_1$ and
$g_2$ in $P_x$, we have that $g_1\prec g_2 \Leftrightarrow
dist(g_1,L_x)<dist(g_2,L_x)$, where $dist$ is the Euclidean distance
in $\R^2$.

\end{ex}


\subsection{An action on the space of orderings}
\label{acting on the space}


One of the most interesting properties of $\mathcal{LO}(G)$ is that
$Aut(G)$, the group of automorphism of $G$, naturally acts on
$\mathcal{LO}(G)$ by continuous transformations. More precisely,
given any $\varphi \in Aut(G)$, we define
$\varphi(\preceq)=\preceq_\varphi \in \mathcal{LO}(G)$ by letting
$h\preceq_\varphi f$ if and only if $\varphi^{-1}(h) \prec
\varphi^{-1}(f)$, where $h,f$ belong to $G$. One easily checks that
$\varphi(U_{g_1,\ldots,g_n})=U_{\varphi(g_1),\ldots, \varphi(g_n)}$.

\vsp

In particular, we obtain an action of $G$ on $\mathcal{LO}(G)$ which
factors throughout the group of {\em inner
automorphisms}\footnote{Inner automorphisms are automorphisms
induced by conjugation of $G$.}. The above condition reads
$\;g(\preceq)=\preceq_g$, where by definition, $h\succ_g f$ if and
only if $ghg^{-1}\succ gfg^{-1}$. We say that $\preceq_g$ is the
conjugate of $\preceq$ under $g$.

\vsp

It immediately follows that the global fixed points of $G$ for this
action are precisely the bi-orderings of $G$. In particular, the
action of $Aut(G)$ on $\mathcal{BO}(G)$ factors throughout
$Out(G)=Aut(G)/Inn(G)$, the group of outer automorphism of $G$.

\vsp

Another typical object in dynamics, namely periodic orbits, also
plays an important role in the theory of orderable groups.

\begin{prop} \label{periodic orbits} {\em In the action of $G$ on
$\mathcal{LO}(G)$, every periodic point (that is, a point whose
orbit is finite) corresponds to a $\ce$-ordering.}
\end{prop}

\noindent {\em Proof:} Suppose that $\preceq$ has a periodic orbit.
Then $Stab_G(\preceq)=\{g\in G\mid g(\preceq)=\preceq\}$ has finite
index in $G$. Moreover, the restriction of $\preceq$ to
$Stab_G(\preceq)$ is a bi-invariant ordering. In particular it is
Conradian. Now, by a result of Rhemtulla and Rolfsen \cite[Theorem
2.4]{rr}, here Corollary \ref{rem-rolf}, we have that $\preceq$ is
Conradian. $\hfill\square$

\vsp

One may wish that any Conradian ordering is periodic, but, as shown
in \cite[Example 4.6]{witte}, there is a group such that no
Conradian ordering is periodic. Our next example shows that on the
Heisenberg group -a group in which every left-ordering is Conradian-
both phenomena appears, {\em i.e.} it admits periodic and non
periodic orbits.

\begin{ex} Let $H=\langle a,b,c\mid [a,b]=c, ac=ca, bc=cb\rangle$
be the discrete Heisenberg group.

\vsp

The group $H$ is left-orderable, since $H/[H,H]\simeq \Z^2$ and
$[H,H]\simeq \Z$. Therefore, we can produce an ordering $\preceq$ on
$H$ by defining an ordering $\preceq_1$ on $H/[H,H]$ and an ordering
$\preceq_2$ on $[H,H]$, and then declaring $id\prec g$ if and only
if $[H,H] \prec_1 g\,[H,H]$ or $g\in [H,H]$ and $id\prec_2 g$. The
ordering $\preceq$ is easily shown to be bi-invariant, so it is a
fixed point for the action of $G$ in $\mathcal{LO}(G)$. Moreover,
due to Proposition \ref{nilp_prop}, we have that every ordering on
$H$ is Conradian.

\vsp

Note that we have the freedom to choose any ordering on
$H/[H,H]\simeq \Z^2$, so we can assume that $\langle b[H,H] \rangle
\lhd H/[H,H]$ is convex in $\preceq_1$. This implies that in the
ordering $\preceq$, the subgroup $\langle b,c\rangle\simeq \Z^2$ is
normal and convex in $H$ (and $[H,H]=\langle c \rangle$ is convex in
$\langle b,c \rangle$). Then, using Corollary \ref{nice coro}, we
can define an ordering $\preceq^\prime$ on $H$ by choosing any
ordering $\preceq^\prime_2$ on $\langle b,c \rangle\simeq\Z^2$ and
an ordering $\preceq_1^\prime$ on $H/\langle b,c\rangle\simeq \Z$.
For concreteness, we let $\preceq_2^\prime$ be an ordering of
irrational type of $\Z^2$, say $P_x$, $x\in \R\setminus \Q$; see
Example \ref{LOZ^2}.

\vsp

We now let $X\subset \mathcal{LO}(H)$ be the orbit of
$\preceq^\prime$. Since $\langle b,c \rangle$ is normal, convex and
Abelian, it acts trivially on $X$. Therefore,
$X=\{\preceq^\prime_{a^n}\mid n\in \Z\}$. To see that
$\preceq^\prime$ is not periodic, it is enough to see that, if
$n\not=0$, then the restrictions of $\preceq^\prime_{a^n}$ and
$\preceq^\prime$ to $\langle b,c \rangle $ do not coincide. To see
this, note that, since $aba^{-1}= bc$, and $aca^{-1}=c$, making the
identifications $e_1=c$ and $e_2=b$, we have that the action by
conjugation of $a$ on $\langle b,c\rangle$  corresponds to the
action on $\Z^2\subset \R^2$ given by the matrix
$$ M_a= \left( \begin{array}{c c} 1 & 1 \\ 0 &1 \end{array}
\right).$$ Now, recall from Example \ref{LOZ^2} that to any ordering
on $\Z^2$ we associate a one dimensional vector space
$L_x\subset\R^2$, namely the kernel of the $\R$-linear map
$e_1\mapsto 1$, $e_2\mapsto x$. Moreover, if $L_x\not= L_y$ then
$P_x\not= P_y$. Checking that $M_a^n(L_x)\not= L_x$ for $x\in
\R\setminus\Q$ and every $n\not=0$ is an easy exercise. In
particular, $\preceq^\prime$ is not periodic.
\end{ex}


\section{Dynamical versions of group orderability}
\label{sec dynanical versions}

Though orderability may look as a very algebraic concept, it has a
deep (one-dimensional) dynamical content. The following theorem, due
to P. Cohn, M. Zaitseva, and P. Conrad, goes in this direction (see
\cite[Theorem 3.4.1]{kopytov}):

\vsp

\begin{thm}\label{Cohn} {\em A group $G$ is left-orderable if and only if it
embeds in the group of (order-preserving) automorphisms of a totally
ordered set.}
\end{thm}

\vsp

Both implications of this theorem are easy. In one direction, note
that a left-ordered group acts on itself by order preserving
automorphisms, namely, left translations. Conversely, to create a
left-ordering on a group $G$ of automorphisms of a totally ordered
set $(\Omega,\leq)$, we construct the what is called \textit{induced
ordering} from the action as follows. Fix a well-order $\, \leq^*
\,$ on the elements of $\,\Omega\,, \,$ and, for every $f\in G\,,
\,$ let $\, w_f = \min_{\leq^*} \{w \in \Omega \mid f(w)\not= w\}$.
Then we define an ordering $\, \preceq \,$ on   $G$ by letting
$f\succ id\,$ if and only if $\;f(w_f)>w_f$. It is not hard to check
that this order relation is a (total) left-ordering on $G$.

\vs

For the case of countable groups, we can give more dynamical
information since we can take  $\Omega$ as being the real line (see
\cite[Theorem 6.8]{ghys}, \cite[\S 2.2.3]{book}, or \cite{navas} for
further details). This characterization will be used in \S 3 and \S
4.1.

\vsp

\begin{prop} \label{real din}{\em For a countable infinite group $G$, the following two properties are
equivalent:\\

\noindent -- $G$ is left-orderable,\\

\noindent -- $G$ acts faithfully on the real line by
orientation-preserving homeomorphisms. That is, there is an
homomorphic embedding $G\to Homeo_+(\R)$.}
\end{prop}

\noindent\textit{Sketch of proof: } The fact that a group of
orientation-preserving homeomorphisms of the real line is
left-orderable is a direct consequence of Theorem \ref{Cohn}.

\vsp

For the converse, we construct what is called \textit{a dynamical
realization of a left-ordering $\preceq$.} Fix an enumeration
$(g_i)_{i \geq 0}$ of $G$, and let $t_\preceq(g_0)=0$. We shall
define an order-preserving map $t_\preceq: G \to \R$ by induction.
Suppose that $t_\preceq(g_0), t_\preceq(g_1),\ldots,t_\preceq(g_i)$
have been already defined. Then if $g_{i+1}$ is greater (resp.
smaller) than all $g_0,\ldots, g_i$, we define $t_\preceq(g_{i+1})=
max\{t_\preceq(g_0),\ldots, t_\preceq(g_i)\}+1$ (resp.
$min\{t_\preceq(g_0),\ldots, t_\preceq(g_i)\}-1$). If $g_{i+1}$ is
neither greater nor smaller than all $g_0,\ldots,g_i$, then there
are $g_n,g_m\in\{g_0,\ldots , g_i \}$ such that $g_n\prec
g_{i+1}\prec g_m$ and no $g_j$ is between $g_n,g_m$ for $0\leq j\leq
i$. Then we set
$t_\preceq(g_{i+1})=(t_\preceq(g_n)+t_\preceq(g_m))/2$.

\vsp

Note that $G$ acts naturally on $t_\preceq(G)$ by $g(t(g_i)) =
t_\preceq(gg_i)$. It is not difficult to see that this action
extends continuously to the closure of $t_\preceq(G)$. Finally, one
can extend the action to the whole real line by declaring the map
$g$ to be affine on each interval in the complement of
$\overline{t_\preceq(G)}$. $\hfill\square$

\vs

We have constructed an embedding $ G \to Homeo_+(\R)$. We call this
embedding a dynamical realization of the left-ordered group
$(G,\preceq)$. The order preserving map $t_\preceq$ is called the
reference map. When the context is clear, we will drop the subscript
$\preceq$ of the map $t_\preceq$.

\vsp

\begin{rem} \label{rem conj1} As constructed above, the dynamical realization depends
not only on the left-ordering $\preceq$, but also on the enumeration
$(g_i)_{i\geq 0}$. Nevertheless, it is not hard to check that
dynamical realizations associated to different enumerations (but the
same ordering) are \textit{topologically conjugate}.\footnote{Two
actions $\phi_1\!: G \to \mathrm{Homeo}_+(\R)$ and $\phi_2\!:G \to
\mathrm{Homeo}_+(\R)$ are topologically conjugate if there exists
$\varphi \in \mathrm{Homeo}_+(\R)$ such that $\varphi\circ \phi_1(g)
= \phi_2(g) \circ \varphi$ for all $g \in G$.} Thus, up to
topological conjugacy, the dynamical realization depends only on the
ordering $\preceq$ of $G$.

\vsp

An important property of dynamical realizations is that they do not
admit global fixed points (\textit{i.e.,} no point is stabilized by
the whole group). Another important property is that $t_\preceq(id)$
has a {\em free orbit} (\textit{i.e} $\{g\in G\mid
g(t_\preceq(id))=t_\preceq(id)\}=\{id\}$ ). Hence $g\succ id$ if and
only if $g(t_\preceq(id))> t_\preceq(id)$, which allows us to
recover the left-ordering from its dynamical realization.
\end{rem}


\subsection{Crossings: a new characterization of Conrad's property}
\label{sec crossings}

The Conrad property has many characterizations; see Theorem \ref{teo
C} and \cite[\S 7.4]{botto}. All of them are algebraic descriptions.
We finish this introductory chapter giving a new characterization of
the Conrad property for left-orderings which has a strong dynamical
flavor. The dynamical object to look at are the so-called
\textit{crossings}. We will make a strong use of this
characterization in Chapter 2.

\vsp

The following definition, first introduced in \cite{navas} for the
case of countable groups and latter in \cite{crossings} for the
general case, will be of great importance in this work. Let $G$ be a
group acting by order preserving bijections on a totally ordered
space $(\Omega,\leq)$. A {{\em crossing}} for the action of $G$ on
$\Omega$ is a 5-uple $\, (f,g,u,v,w) \,$ where $f,g$ (resp. $u,v,w$)
belong to $G$ (resp. $\Omega$) and satisfy:

\vsp

\noindent $(i)$ \esp $u < w < v$.

\vsp

\noindent $(ii)$ \esp For every $\,n \in \mathbb{N} , \,$ we have
$\,g^n u < v $ \esp and \esp $f^n v > u \,.$

\vsp

\noindent $(iii)$ \esp There exist $M,N$ in $\mathbb{N}$ such that
\esp $f^N v < w < g^M u$.

\vspace{0.2cm}


\beginpicture

\setcoordinatesystem units <1cm,1cm>


\putrule from 1.5 -2.5 to 6.5 -2.5 \putrule from 1.5 2.5 to 6.5 2.5
\putrule from 1.5 -2.5 to 1.5 2.5 \putrule from 6.5 -2.5 to 6.5 2.5

\plot 1.5 0 1.625 0.05296 1.75 0.109 1.875 0.16056 2 0.216 2.125
0.26056 2.25 0.3009 2.375 0.35296 2.5 0.4 /

\plot 2.5 0.4 2.625 0.45296 2.75 0.509 2.875 0.56056 3 0.616 3.125
0.66056 3.25 0.7009 3.375 0.75296 3.5 0.8 /

\plot 3.5 0.8 3.625 0.85296 3.75 0.909 3.875 0.96056 4 1.016 4.125
1.06056 4.25 1.1009 4.375 1.15296 4.5 1.2 /

\plot 4.5 1.2 4.625 1.25296 4.75 1.309 4.875 1.36056 5 1.416 5.125
1.46056 5.25 1.5009 5.375 1.55296 5.5 1.6 5.625 1.65296 /


\plot 6.5 -0.82519 6.375 -0.85296 6.125 -0.909 5.875 -0.96056 5.625
-1.016 5.375 -1.06056 5.125 -1.1009 4.875 -1.15296 4.625 -1.2 /

\plot 4.625 -1.2 4.375 -1.25296 4.125 -1.309 3.875 -1.36056 3.625
-1.416 3.375 -1.46056 3.125 -1.5009 2.875 -1.55296 2.625 -1.6 2.375
-1.65296 /


\setdots

\plot 1.5 -2.5 6.5 2.5 /

\putrule from 2.32 -1.68 to 5.68 -1.68 \putrule from 2.32 -1.68 to
2.32 1.68 \putrule from 2.32 1.68 to 5.68 1.68 \putrule from 5.68
-1.68 to 5.68 1.68

\put{Figure 1.1: A crossing} at 4 -3.5 \put{} at -4.2 0

\small


\put{$u$} at 1.5 -2.8 \put{$v$} at 6.5 -2.8 \put{$w$} at 4   -2.8
\put{$f^N v$} at 2.8 -2.8 \put{$g^M u$} at 5.2 -2.8 \put{$\bullet$}
at 1.5 -2.5 \put{$\bullet$} at 6.5 -2.5 \put{$\bullet$} at 4   -2.5
\put{$\bullet$} at 2.8 -2.5 \put{$\bullet$} at 5.2 -2.5

\put{$f$} at 4.7 -0.9 \put{$g$} at 3.4 1

\endpicture


\vspace{0.5cm}

\vsp \vsp

The reason why this definition is so important is because it
actually characterizes the $\ce$-orderings, as is shown in
\cite[Theorem 1.4]{crossings}. We quote the theorem and the proof
below.

\vsp

\vsp\vsp

\noindent {\bf Theorem A.} {\em A left-ordering $\, \preceq \,$ on
$G$ is Conradian if and only if the action of $G$ by left
translations on itself admits no crossing.}

\vsp\vsp

\noindent {\em Proof:} Suppose that $\preceq$ is not Conradian, and
let $f,g$ be positive elements so that $fg^n \prec g$ for every $n
\in \mathbb{N}$. We claim that $(f,g,u,v,w)$ is a crossing for
$(G,\preceq)$ for the choice $u = 1$, and $v = f^{-1}g$, and $w =
g^2$. Indeed:

\vsp

\noindent -- From $fg^2 \prec g$ one obtains $g^2 \prec f^{-1} g$,
and since $g \succ 1$, this gives \esp $1 \prec g^2 \prec f^{-1}g,$
\esp that is, \esp $u \prec w \prec v$.

\vsp

\noindent -- From $f g^n \prec g$ one gets \esp $g^n \prec f^{-1}
g$, \esp that is, \esp $g^n u \prec v$ \esp (for every $n \in
\mathbb{N}$); on the other hand, since both $f,g$ are positive, we
have $f^{n-1} g \succ 1$, and thus \esp $f^n (f^{-1}g) \succ 1$,
\esp that is, \esp $f^n v \succ u$ \esp (for every $n \in
\mathbb{N}$).

\vsp

\noindent -- The relation \esp $f (f^{-1}g) = g \prec g^2$ \esp may
be read as \esp $f^N v \prec w$ \esp for $N=1$; on the other hand,
the relation \esp $g^2 \prec g^3$ \esp is \esp \esp $w \prec g^M u$
\esp for $M = 3$.

\vsp

Conversely, assume that $(f,g,u,v,w)$ is a crossing for
$(G,\preceq)$ so that \esp $f^N v \prec w \prec g^M u$ \esp (with
$M,N$ in $\mathbb{N}$). We will prove that $\preceq$ is not
Conradian by showing that, for $h = g^Mf^N$ and $\bar{h} = g^M$,
both elements $w^{-1} h w$ and $w^{-1} \bar{h} w$ are positive, but
$$(w^{-1} h w) (w^{-1} \bar{h} w)^n \prec w^{-1} \bar{h} w
\quad \mbox{ for all } \esp n \in \mathbb{N}.$$ To show this, first
note that \esp $gw \succ w$. \esp Indeed, if not, then we would have
$$w \prec g^N u \prec g^N w \prec g^{N-1} w \prec \ldots \prec gw \prec w,$$
which is absurd. Clearly, the inequality \esp $gw \succ w$ \esp
implies
$$g^M w \succ g^{M-1} w \succ \ldots \succ gw \succ w,$$
and hence
\begin{equation}
w^{-1} \bar{h} w = w^{-1} g^M w \succ 1. \label{zero}
\end{equation}
Moreover,
$$h w = g^M f^N w \succ g^M f^N f^N v
= g^M f^{2N} v \succ g^{M} u \succ w.$$ and hence
\begin{equation}
w^{-1} h w \succ 1. \label{one}
\end{equation}
Now note that for every $n \in \mathbb{N}$,
$$h \bar{h}^n w = h g^{Mn} w \prec h g^{Mn} g^M u
= h g^{Mn+M} u \prec h v = g^M f^N v \prec g^M w = \bar{h} w.$$
After multiplying by the left by $w^{-1}$, the last inequality
becomes
$$(w^{-1} h w) (w^{-1} \bar{h} w)^n =
w^{-1} h \bar{h}^n w \prec w^{-1} \bar{h} w,$$ as we wanted to
check. Together with (\ref{zero}) and (\ref{one}), this shows that
$\preceq$ is not Conradian. \hfill$\square$


\chapter[{\small Applications of the new characterization}]{Applications of the new characterization of Conrad's property}

With the aid of Theorem A, we are able to prove two structure
theorems for $\ce$-orderings. One is the algebraic description of
groups with finitely many $\ce$-orderings; Theorem B in \S 2.1. The
second tell us that the spaces of Conradian groups orderings are
either finite or without isolated points; Theorem C in \S 2.2. This
is essentially taken from \cite{rivas}. As an application, we give a
new proof of a theorem of Navas in \cite{navas} asserting that
torsion free nilpotent groups have no isolated left-orderings unless
they are rank-one Abelian; see \S 2.2.2.

\vsp

Theorem A also serves us to detect the so-called {\em Conradian soul
of a left-ordering}; see \S 2.3.1. This allows us to ensure many
good conjugacy properties of the left action of the group on itself
and to investigate the possibility of approximating a given
left-ordering by its conjugates, which in many cases can be done; \S
2.3.2. This is used to give a new proof of the fact, first proved by
Linnell, that a left-orderable group admits either finitely or
uncountably many left-orderings; Theorem E in \S 2.3.3. This is
essentially taken from \cite{crossings}, which, in turn, is inspired
from \cite{navas}.


\section{On groups with finitely many Conradian orderings}
\label{finitos C-ord}

In this section we give an algebraic description of groups admitting
only finitely many Conradian orderings \cite{rivas}. This may be
considered as an analogue of Tararin's classification of groups
admitting only finitely many left-orderings; see Theorem \ref{teo T}
or \cite[Theorem 5.2.1]{kopytov}. For the statement, recall that a
series
$$\{ id \} = G_0 \lhd G_{1} \lhd \ldots \lhd G_{n-1} \lhd G_n = G$$
is said to be {\em rational} if it is subnormal ({\em i.e.,} each
$G_{i}$ is normal in $G_{i+1}$) and each quotient $G_{i+1} / G_{i}$
is torsion-free rank-one Abelian. The series is called {\em normal}
if, in addition, each $G_i$ is normal in $G$.

\vsp\vsp

\noindent {\bf Theorem B.} {\em Let $G$ be a $\ce$-orderable group.
If $G$ admits finitely many $\ce$-orderings, then $G$ admits a
unique (hence normal) rational series. In this series, no quotient
$G_{i+2}/G_{i}$ is Abelian. Conversely, if $G$ is a group admitting
a  normal rational series
$$\{ id \} = G_0 \lhd G_{1} \lhd \ldots \lhd G_{n-1} \lhd G_n = G$$
so that no quotient $G_{i+2} / G_{i}$ is Abelian, then the number of
$\ce$-orderings on $G$ equals $2^n$.}

\vsp\vsp

\subsection{Some lemmata}

\vsp

The following crucial lemma is essentially proved by Navas in
\cite{navas} in the case of countable groups, but the proof therein
rests upon very specific issues about the dynamical realization of
an ordered group. Here we give a general algebraic proof.

\vsp

\begin{lem} \label{pullback} {\em Suppose $G$ is faithfully acting by
order preserving bijections on a totally ordered set
$(\Omega,\leq)$. Then, the action has no crossings if and only if
any induced ordering is Conradian.} \label{generalization}
\end{lem}


\noindent{\em Proof:} Suppose that the ordering $\, \preceq \,$ on
$G$ {\em induced} from some well-order $\, \leq^* \,$ on $\Omega$ is
not Conradian (comments after Theorem \ref{Cohn} explain how to
induce an ordering from an action, and shows the relation between
$\preceq$ and $\leq^*$). Then there are $\,\preceq$-positive
elements $f,g$ in $G$ such that $fg^{n} \prec g$, for every $n \in
\mathbb{N}$. This easily implies $f \prec g$, since, in the case
$id\prec g \prec f$, we get $id \prec g\prec g^n$, so $id \prec
g\prec f\prec fg \prec fg^n$.

\vsp

Let $\bar{w} = \min_{\leq^*} \{ w_f,\, w_g \}$, where, for $h\in G$,
$w_h=\min_{\leq^*}\{w\in \Omega\mid h(w)\not=w\}$.  We claim that
$(fg,fg^2, \bar{w}, g(\bar{w}), fg^2(\bar{w}))$ is a crossing (see
Figure 2.1). Indeed, the inequalities $id \prec f \prec g$ imply
that $\bar{w}=w_g\leq^* w_f$ and $g(\bar{w})>\bar{w}$. Moreover,
$f(\bar{w}) \geq \bar{w} \,, \,$ which together with $fg^{n} \prec g
\,$ yield $\bar{w}< fg^2(\bar{w})< g(\bar{w})$, hence condition
$(i)$ of the definition of crossing is satisfied. Note that the
preceding argument actually shows that $fg^n(\bar{w})<g(\bar{w})$,
for all $n \in \mathbb{N}\,. \,$ Thus $fg^2fg^2(\bar{w})<
fg^3(\bar{w})< g(\bar{w})$. A straightforward induction argument
shows that $(fg^2)^n(\bar{w})< g(\bar{w})$, for all $n \in
\mathbb{N}$, which proves the first part of condition $(ii)$. For
the second part, from $g(\bar{w})>\bar{w}\,$ and
$\,f(\bar{w})\geq\bar{w} \,$ we conclude that $\, \bar{w}<
(fg)^n(g(\bar{w}))\,. \,$ Condition $(iii)$ follows because $\bar{w}
< fg^2(\bar{w})$ implies $fg^2(\bar{w}) <
fg^2(fg^2(\bar{w}))=(fg^2)^2(\bar{w})$, and $fg^2(\bar{w}) <
g(\bar{w})$ implies $(fg)^2(g(\bar{w}))= fg (fg^2(\bar{w})) < fg
(g(\bar{w}))=fg^2(\bar{w})$.

\vsp

For the converse, suppose that $(f,g,u,v,w)$ is a crossing for the
action. Let $N,M$ in $\N$ be such that $f^N(v)<w$ and $g^M(u)> w$;
see Figure 1.1. We let $\preceq_u$ be any induced ordering on $G$
with $u$ as first reference point. In particular, $g^M\succ_u id$
and $f^Ng^M \succ_u id$. We claim that $\preceq_u$ is not Conradian.
Indeed, we have that $(f^Ng^M) g^{2M} \prec_u g^M$, since
$w<g^M(u)<v$ and $f^Mg^{3N}(u)<f^M(v)<w$. $\hfill\square$

\vspace{0.5cm}


\beginpicture

\setcoordinatesystem units <1cm,1cm>


\putrule from 1.5 -2.5 to 6.5 -2.5 \putrule from 1.5 2.5 to 6.5 2.5
\putrule from 1.5 -2.5 to 1.5 2.5 \putrule from 6.5 -2.5 to 6.5 2.5

\plot 1.5 0 1.625 0.05296 1.75 0.109 1.875 0.16056 2 0.216 2.125
0.26056 2.25 0.3009 2.375 0.35296 2.5 0.4 /

\plot 2.5 0.4 2.625 0.45296 2.75 0.509 2.875 0.56056 3 0.616 3.125
0.66056 3.25 0.7009 3.375 0.75296 3.5 0.8 /

\plot 3.5 0.8 3.625 0.85296 3.75 0.909 3.875 0.96056 4 1.016 4.125
1.06056 4.25 1.1009 4.375 1.15296 4.5 1.2 /

\plot 4.5 1.2 4.625 1.25296 4.75 1.309 4.875 1.36056 5 1.416 5.125
1.46056 5.25 1.5009 5.375 1.55296 5.5 1.6 5.625 1.65296 /


\plot 6.5 -0.82519 6.375 -0.85296 6.125 -0.909 5.875 -0.96056 5.625
-1.016 5.375 -1.06056 5.125 -1.1009 4.875 -1.15296 4.625 -1.2 /

\plot 4.625 -1.2 4.375 -1.25296 4.125 -1.309 3.875 -1.36056 3.625
-1.416 3.375 -1.46056 3.125 -1.5009 2.875 -1.55296 2.625 -1.6 2.375
-1.65296 /


\setdots

\plot 1.5 -2.5 6.5 2.5 /

\putrule from 2.32 -1.68 to 5.68 -1.68 \putrule from 2.32 -1.68 to
2.32 1.68 \putrule from 2.32 1.68 to 5.68 1.68 \putrule from 5.68
-1.68 to 5.68 1.68

\put{Figure 2.1: The crossing} at 4 -3.5 \put{} at -4.2 0

\small


\put{$\bar{w}$} at 1.5 -2.8 \put{$g(\bar{w})$} at 6.5 -2.8
\put{$fg^2(\bar{w})$} at 4 -2.8  \put{$\bullet$} at 1.5 -2.5
\put{$\bullet$} at 6.5 -2.5 \put{$\bullet$} at 4   -2.5

\put{$fg$} at 4.7 -0.9 \put{$fg^2$} at 3.4 1.1

\endpicture


\vspace{0.5cm}

As an application of Theorem A and/or Lemma \ref{pullback}, we give
a new proof of a theorem proved in \cite{rr}.

\begin{cor}\label{rem-rolf} {\em Let $\preceq$ be an ordering on a group $G$.
Let $H$ be a subgroup of finite index such that $\preceq$ restricted
to $H$ is Conradian. Then $\preceq$ is a $\ce$-ordering.}
\end{cor}

\noindent \textit{Proof:} Suppose, by way of a contradiction, that
$\preceq$ is a non Conradian ordering of $G$. Then, there are $f,g$
in $G$, both $\preceq$-positive, such that $fg^n\prec g$ for all
$n\in\N$. By (the proof of) Theorem A, $(f,g,\,id,f^{-1}g,g^2)$ is a
crossing for the left translation action of $G$ on itself. But then,
for any $n\in \N$, we have that $(f^n,g^n,\, id,f^{-n}g^{n},g^{2n})$
is still a crossing for the action (see for instance Figure 2.2
below). But this is a contradiction, since, for certain $n$ big
enough, $f^n$ and $g^n$ belongs to $H$, thus $(f^n,g^n,\,
id,f^{-n}g^{n},g^{2n})$ is a crossing for the left translation
action of $H$ on itself, which, by Theorem A, implies that the
restriction of $\preceq$ to $H$ is non Conradian, contrary to our
assumption. $\hfill\square$

\vspace{0.5cm}


\beginpicture

\setcoordinatesystem units <1cm,1cm>


\putrule from 1.5 -2.5 to 6.5 -2.5 \putrule from 1.5 2.5 to 6.5 2.5
\putrule from 1.5 -2.5 to 1.5 2.5 \putrule from 6.5 -2.5 to 6.5 2.5

\plot 1.5 0 1.625 0.05296 1.75 0.109 1.875 0.16056 2 0.216 2.125
0.26056 2.25 0.3009 2.375 0.35296 2.5 0.4 /

\plot 2.5 0.4 2.625 0.45296 2.75 0.509 2.875 0.56056 3 0.616 3.125
0.66056 3.25 0.7009 3.375 0.75296 3.5 0.8 /

\plot 3.5 0.8 3.625 0.85296 3.75 0.909 3.875 0.96056 4 1.016 4.125
1.06056 4.25 1.1009 4.375 1.15296 4.5 1.2 /

\plot 4.5 1.2 4.625 1.25296 4.75 1.309 4.875 1.36056 5 1.416 5.125
1.46056 5.25 1.5009 5.375 1.55296 5.5 1.6 5.625 1.65296 /


\plot 6.5 -0.82519 6.375 -0.85296 6.125 -0.909 5.875 -0.96056 5.625
-1.016 5.375 -1.06056 5.125 -1.1009 4.875 -1.15296 4.625 -1.2 /

\plot 4.625 -1.2 4.375 -1.25296 4.125 -1.309 3.875 -1.36056 3.625
-1.416 3.375 -1.46056 3.125 -1.5009 2.875 -1.55296 2.625 -1.6 2.375
-1.65296 /


\setdots

\plot 1.5 -2.5 6.5 2.5 /

\putrule from 2.32 -1.68 to 5.68 -1.68 \putrule from 2.32 -1.68 to
2.32 1.68 \putrule from 2.32 1.68 to 5.68 1.68 \putrule from 5.68
-1.68 to 5.68 1.68

\put{Figure 2.2} at 4 -3.5  \put{} at -4.2 0

\small


\put{$id$} at 1.5 -2.8 \put{$f^{-1}g$} at 6.5 -2.8 \put{$g^2$} at
4.9 -2.8  \put{$\bullet$} at 1.5 -2.5 \put{$\bullet$} at 6.5 -2.5
\put{$\bullet$} at 4.8   -2.5

\put{$f$} at 4.7 -0.9 \put{$g$} at 3.4 1.1

\endpicture


\vspace{0.5cm}

Note that if in (the proof of) Lemma \ref{pullback} we let $w_0$ be
the smallest element (with respect to $\leq^*$) of $\Omega \,, \,$
then the stabilizer of $\, w_0 \,$ is $\,\preceq$-convex. Indeed, if
$id\prec g \prec f$, with $f(w_0)=w_0 , \,$ then $\, w_0 <^* w_f
\leq^* w_g$, thus $g(w_0)=w_0\,. \,$ Actually, it is not hard to see
that the same argument shows the following.

\vsp \vsp

\begin{prop} \label{laclave} {\em Let $\Omega$ be a set endowed with a well-order
$\,\leq^* $. If a group $G$ acts faithfully on $\Omega$ preserving a
total order on it, then there exists a left-ordering on $G$ for
which the stabilizer $G_{\Omega_0}$ of any initial segment $\Omega_0
\,$ of $\, \Omega \,$ (w.r.t. $\leq^*$) is convex. Moreover, if the
action has no crossing, then this ordering is Conradian.}
\end{prop}

\vsp

\begin{ex} A very useful example of an action without crossings is
the action by left translations on the set of left-cosets of any
subgroup $H$ which is convex with respect to a $\ce$-ordering
$\,\preceq \,$ on $G$. Indeed, it is not hard to see that, due to
the convexity of $H$, the order $\, \preceq \,$ induces a total
order $\, \preceq_H \,$ on the set of left-cosets $G / H$. Moreover,
$\, \preceq_H \,$ is $G-$invariant. Now suppose that $\, (f, g, uH,
vH, wH) \,$ is a crossing for the action. Since $w_1H\prec_{H} w_2
H$ implies $\,w_1\prec w_2 \,, \,$ for all $w_1,\, w_2$ in $G$, we
have that $(f,g,u,v,w)$ is actually a crossing for the action by
left translations of $G$ on itself. Nevertheless, this contradicts
Theorem A. \label{convex-cosets}
\end{ex}

\vsp

The following is an application of the preceding example. For the
statement, we will say that a subgroup $H$ of a group $G$ is
\textit{$\ce$-relatively convex} if there exists a $\ce$-ordering on
$G$ for which $H$ is convex.

\vsp

\begin{lem} {\em For every $\ce$-orderable group, the intersection of any
family of $\ce$-relatively convex subgroups is $\ce$-relatively
convex.} \label{int-conv}
\end{lem}

\noindent{\em Proof:} We consider the action of $G$ by left
multiplications on each coordinate of the set $\Omega=\prod_\alpha G
/ H_\alpha \,, \,$ where $(G / H_\alpha,\;\preceq_{H_\alpha})$ is
the ($G-$invariant ordered) set of left-cosets of the
$\ce$-relatively convex subgroup $H_\alpha$. Putting the (left)
lexicographic order on $\Omega$ and using Example
\ref{convex-cosets}, it is not hard to see that this action has no
crossing. Moreover, since $\{id\}$ is $\ce$-convex, the action is
faithful.

\vsp

Now consider an arbitrary family
$\Omega_0\subset\{H_\alpha\}_\alpha$ of $\ce$-relatively convex
subgroups of $G$, and let $\leq^*$ be a well-order on $\Omega$ for
which $\Omega_0$ is an initial segment. For the induced ordering
$\preceq$ on $G$, it follows from Proposition \ref{laclave} that the
stabilizer $G_{\Omega_0} = \bigcap_{H\in\, \Omega_0} H$ is
$\preceq$-convex. Moreover, Lemma \ref{generalization} implies that
$\preceq$ is a $\ce$-ordering, thus concluding the proof. $\hfill
\square$
\newline

We close this section with a simple lemma that we will need later
and which may be left as an exercise to the reader (see also
\cite[Lemma 5.2.3]{kopytov}).

\vsp

\begin{lem} {\em Let $G$ be a torsion-free Abelian group. Then $G$
admits only finitely many $\ce$-orderings if and only if $G$ has
rank-one}. \label{abel-1}
\end{lem}


\subsection{Proof of Theorem B}
\label{sec}


Let $G$ be a $\ce$-orderable group admitting only finitely many
$\ce$-orderings. Obviously, each of these orderings must be isolated
in $\mathcal{CO}(G)$. We claim that, in general, if $\, \preceq \,$
is an isolated $\ce$-ordering, then the series of $\,
\preceq$-convex subgroups
$$\{id\} =  G^{id} \subset \ldots \subset G_g\lhd G^g \subset \ldots \subset G$$
is finite. Indeed, let $\{f_1, \ldots, f_n\}\subset G$ be a set of
$\, \preceq$-positive elements such that $V_{f_1,\ldots,f_n}$
consists only of $\, \preceq . \,$ If the series above is infinite,
then there exists $g\in G$ so that no $f_i$ belongs to $G^g\setminus
G_g$. This implies that the flipped ordering $\,\preceq^g ,$ defined
in Example \ref{exflip},  is Conradian and different from $\,\preceq
. \,$ However, every $f_i$ is still $\, \preceq^g$-positive, which
is impossible because $V_{f_1,\ldots,f_n} = \{\preceq\}$.

\vsp

Next let
$$\{id\} = G_0 \lhd G_1 \lhd \ldots \lhd G_n = G$$
be the series of $\, \preceq$-convex subgroups of $G.\,$ According
to Theorem \ref{teo C}, every quotient $G_i/G_{i-1}$ embeds into
$\mathbb{R}$, and thus it is Abelian. Since every ordering on such a
quotient can be extended to an ordering on $G$ (see for instance
Corollary \ref{nice coro}), the Abelian quotient $G_i/G_{i-1}$ has
only a finite number of orderings. It now follows from Lemma
\ref{abel-1} that it must have rank one. Therefore, the series above
is rational.

We now show that this series is unique. Suppose
$$\{id\}=H_0\lhd H_1\lhd \ldots \lhd H_k =G$$
is another rational series. Since $H_{k-1}$ is $\ce$-relatively
convex, we conclude that
$\, N=G_{n-1}\cap H_{k-1} \,$ is $\,C$-relatively convex by Lemma
\ref{int-conv}. Now $\,G/N \,$ is torsion-free Abelian and has only
a finite number of orderings, thus it has rank one. But $H_{k-1}$
and $G_{n-1}$ have the property that $w^n \in G_{n-1}$ (resp. $w^n
\in H_{k-1}$) implies $w \in G_{n-1}$ (resp. $w \in H_{k-1}$). This
clearly yields $H_{k-1} = G_{n-1}$. Repeating this argument several
times, we conclude that rational series is unique, hence normal.

\vsp

Now we claim that no quotient $\, G_{i+2}/G_i \,$ is Abelian. If
not, $G_{i+2}/G_i$ would be a rank-two Abelian group, and so an
infinite number of orderings could be defined on it. But since every
ordering on this quotient can be extended to a $\ce$-ordering on
$G$, this would lead to a contradiction.

\vs

We now prove the converse of Theorem B.
\medskip

Suppose that $\,G \,$ has a normal rational series
$$\{id\}=G_0 \lhd G_1\lhd ... \lhd G_n=G$$
such that no quotient $\,G_{i+2}/G_i \,$ is Abelian. Clearly,
flipping the orderings on the quotients $\,G_{i+1}/G_i \,$ we obtain
at least $2^n$ many $\ce$-orderings on $G$. We claim that these are
the only possible $\ce$-orderings on $G$. Indeed, let $\, \preceq
\,$ be a $\ce$-ordering on $G$, and let $a \in G_1$ and $b \in
G_2\setminus G_1$ be two non-commuting elements. Denoting the Conrad
homomorphism of the group $\, \langle a,b\rangle \,$ endowed with
the restriction of $\, \preceq \,$ by $\, \tau , \,$ we have
$\tau(a) = \tau(bab^{-1})$. Since $G_1$ is rank-one Abelian, we have
$bab^{-1}=a^r$ for some rational number $r\not=1$. Thus
$\tau(a)=r\tau(a)$, which implies that $\tau(a)=0$. Since
$\tau(|b|)>0$, where $|b| = max\{b^{-1},b\} $, we have that $\,
a^n\prec |b| \,$ for every $n \in \mathbb{Z}\,$. Since $G_2/G_1$ is
rank-one, this actually holds for every $b \neq id$ in $G_2
\setminus G_1$. Thus $G_1$ is convex in $G_2$.

\vsp

Repeating the argument above, though now with $G_{i+1}/G_i$ and
$G_{i+2}/G_i$ instead of $G_1$ and $G_2$, respectively, we see that
the rational series we began with is none other than the series
given by the convex subgroups of $\,\preceq . \,$ Since each
$G_{i+1}/G_i$ is rank-one Abelian, if we choose $b_i \in G_{i+1}
\setminus G_i$ for each $i=0,\ldots,n-1$, then any $\ce$-ordering on
$G$ is completely determined by the signs of these elements. This
shows that $G$ admits precisely $2^n$ different $\ce$-orderings.


\subsection {An example of a group with $2^n$ Conradian orderings but
infinitely many left-orderings}


The classification of groups having finitely many left-orderings,
here stated as Theorem \ref{teo T}, was obtained by Tararin and
appears in \cite[\S 5.2]{kopytov}. In any of those groups, the
number of left-orderings is $2^n$ for some $n\in \N$. An example of
a group having precisely $2^n$ orders is \esp $T_n = \Z^n$ endowed
with the product rule
$$(\alpha_n, \ldots , \alpha_1) \cdot (\alpha'_n, \ldots, \alpha'_1)=(\alpha_n+\alpha'_n, \, (-1)^{\alpha'_n}
\alpha_{n-1} + \alpha'_{n-1}\, , \ldots , \,(-1)^{\alpha'_2}
\alpha_1+\alpha'_1).$$

A presentation for $T_n$ is $ \langle a_n,\ldots,a_1 \mid
R_n\rangle,$ where the set of relations $R_n$ is
$$\;a_{i+1}a_i a_{i+1}^{-1} = a_i^{-1} \quad \mbox{if}
\quad i < n, \qquad \mbox{and } \qquad a_ia_j=a_ja_i \quad \mbox{if}
\quad |i - j| \geq 2.$$

\vsp

Since the action of a group with only finitely many left-orderings
on its corresponding space of left-orderings has only periodic
orbits, Proposition \ref{periodic orbits} implies

\begin{cor} \label{finitely leftorders}{\em Any left-ordering on a group with only finitely many left-orderings
is Conradian.}
\end{cor}

\vsp

Therefore, it is natural to ask whether for each $n \geq 2$ there
are groups having precisely $2^n$ Conradian orderings but infinitely
many left-orderings. As it was shown in \cite{rivas}, for $n=2$ this
is the case of the Baumslag-Solitar group $B(1,\ell)$, $\ell\geq 2$;
see also \S 3.1. It turns out that, in order to construct examples
for larger $n$ and having $B(1,\ell)$ as a quotient by a normal
convex subgroup, we need to choose an odd integer $\ell$. As a
concrete example, for $n \geq 1$, we endow \esp $C_n = \Z \times
\Z[\frac{1}{3}] \times \Z^{n}$ (where $\Z[\frac{1}{3}]$ denotes the
group of triadic rational numbers) \esp with the multiplication
\begin{multline*}
\Big( \gamma, \,\frac{\eta}{3^{\kappa}}\,, \alpha_n, \ldots,
\alpha_1\Big) \cdot
\Big(\gamma^\prime, \,\frac{\eta^\prime}{3^{\kappa^\prime}}\,, \alpha_n^\prime, \ldots, \alpha_1^\prime \Big) =\\
= \Big( \gamma + \gamma^\prime,
\;3^{-\gamma^\prime}\frac{\eta}{3^{\kappa}}+\frac{\eta^\prime}{3^{\kappa^\prime}}\;,
(-1)^{\eta^\prime} \alpha_n+\alpha^\prime_n \,,
(-1)^{\alpha^\prime_n} \alpha_{n-1} + \alpha^\prime_{n-1}, \ldots,\,
(-1)^{\alpha^\prime_2} \alpha_1 + \alpha^\prime_1 \Big),
\end{multline*}
Note that the product rule is well defined because if \esp
$\eta/3^\kappa = \bar{\eta}/3^{\bar{\kappa}}$, \esp then \esp
$(-1)^{\eta}=(-1)^{\bar{\eta}}$ \esp (it is here where we use the
fact that \esp $\ell = 3$ \esp is odd).

\begin{lem} {\em The group $C_n$ admits the presentation
$$C_n\cong \langle c,b, a_n,\ldots ,a_1 \mid  cbc^{-1}=b^3\,,\;
ca_i=a_ic\,,\; ba_nb^{-1}= a_n^{-1}\,,\; ba_i=a_ib \esp \text{ if }
\esp i\not=n \,  ,\;R_n \rangle,$$ where $R_n$ is the set of
relations of $T_n$ above.}
\end{lem}

\noindent{\em Proof:} Let $\tilde{C_n}$ be the group with
presentation $\langle c,b, a_n,\ldots ,a_1 \mid  cbc^{-1}=b^3\,,\;
ca_i=a_ic\,,\; ba_nb^{-1}= a_n^{-1}\,,\; ba_i=a_ib \esp \text{ if }
\esp i\not=n \,  ,\;R_n \rangle$. Let $e_i=(0, \ldots,0, 1,0,
\ldots,0) \in C_n$, where $1$ is in the $i$-th vector entrance from
the right, that is, $e_1=(0,\ldots,0,1)$. Clearly,
$\{e_i\}_{i=1}^{n+2}$ is a generating set for $C_n$. A direct
calculation shows that
$$\;e_{i+1}e_i e_{i+1}^{-1} = e_i^{-1} \quad \mbox{if}
\quad i \leq n, \quad \mbox{and } \quad e_ie_j=e_je_i \quad
\mbox{if} \quad |i - j| \geq 2, \;\,1\leq i\leq j\leq n+2.$$
Moreover,
$$ e_{n+2}e_{n+1} e^{-1}_{n+2}= e_{n+1}^3\;,\; e_{n+2}e_i=e_ie_{n+2} \;
\text{ if }\; i\leq n \quad \text{and} \quad e_{n+1}e_i=e_ie_{n+1}
\;\text{ if }\; i<n.$$ This means that $\psi:\tilde{C_n}\to C_n$
defined by $\psi(c)=e_{n+2}$, $\psi(b)=e_{n+1}$, and $\psi(a_i)=e_i$
for $1\leq i\leq n$, is a surjective homomorphism.

\vsp

We now let $\varphi:C_n\to \tilde{C_n}$ defined by $\varphi((\gamma,
\frac{\eta}{3^\kappa}, \alpha_n,\ldots,\alpha_1))=c^\gamma
(c^{-\kappa} b^\eta c^\kappa) a_n^{\alpha_n}\ldots a_1^{\alpha_1}$,
where we are assuming that $\eta\in \Z$. Note that, if
$\eta/3^{\kappa}= \bar{\eta}/3^{\bar{\kappa}}$, where both $\eta$
and $\bar{\eta}$ are integers, then $c^{-\kappa}b^\eta
c^\kappa=c^{-\bar{\kappa}}b^{\bar{\eta}}c^{\bar{\kappa}}$. For
instance, if $\kappa\leq \bar{\kappa}$, then
$c^{\bar{\kappa}-\kappa}b^{\eta} c^{-\bar{\kappa}+\kappa}=
b^{3^{\bar{\kappa}-\kappa}\eta}=b^{\bar{\eta}} $. Therefore,
$\varphi$ is a well-defined function. To check that $\varphi$ is a
homomorphism, let $\omega=(\gamma, \frac{\eta}{3^\kappa},
\alpha_n,\ldots,\alpha_1)$ and $ \omega^\prime=(\gamma^\prime,
\frac{\eta^\prime}{3^{\kappa^\prime}},
\alpha^\prime_n,\ldots,\alpha^\prime_1)$. If $\kappa+\gamma^\prime$
and $ {\kappa^\prime}$ are positive, then we have
\begin{eqnarray*}
\varphi(\omega)\varphi(\omega^\prime) &=&  c^\gamma (c^{-\kappa}
b^\eta c^\kappa) a_n^{\alpha_n}\ldots a_1^{\alpha_1}
c^{\gamma^\prime} (c^{-\kappa^\prime} b^{\eta^\prime}
c^{\kappa^\prime}) a_n^{\alpha^\prime_n}\ldots
a_1^{\alpha^\prime_1}\\
 &=& c^{\gamma+\gamma^\prime} (c^{-\kappa-\gamma^\prime} b^\eta
c^{\kappa+\gamma^\prime}) (c^{-\kappa^\prime}b^{\eta^\prime}
c^{\kappa^\prime}  ) \;\;\; a_n^{(-1)^{\eta^\prime}\alpha_n}
a_{n-1}^{\alpha_{n-1}} \ldots a_1^{\alpha_1} a_n^{\alpha^\prime_n}\ldots a_1^{\alpha^\prime_1} \\
&=&  c^{\gamma+\gamma^\prime}
(c^{-\kappa-\gamma^\prime-\kappa^\prime} b^{3^{\kappa^\prime}\eta +
3^{\kappa+\gamma^\prime}\eta^\prime}
c^{\kappa+\gamma^\prime+\kappa^\prime})
a_n^{(-1)^{\eta^\prime}\alpha_n + \alpha_n^\prime}
a_{n-1}^{(-1)^{\alpha_n^\prime}\alpha_{n-1}+\alpha_n^\prime}
 \ldots a_1^{(-1)^{\alpha_2^\prime}\alpha_1 +\alpha_1^\prime} \\
 &=&\varphi(   \gamma + \gamma^\prime,
\;\frac{3^{\kappa^\prime}\eta + 3^{\kappa+\gamma^\prime}
\eta^\prime}{3^{\kappa+\gamma^\prime+\kappa^\prime}},
(-1)^{\eta^\prime} \alpha_n+\alpha^\prime_n \,,
(-1)^{\alpha^\prime_n} \alpha_{n-1} + \alpha^\prime_{n-1}, \ldots,\,
(-1)^{\alpha^\prime_2} \alpha_1 + \alpha^\prime_1  ))\\
 &=& \varphi(\omega\omega^\prime).
\end{eqnarray*}
The equality
$\varphi(\omega)\varphi(\omega^\prime)=\varphi(\omega\omega^\prime)$
in the other cases can be checked similarly. Therefore, $\varphi$ is
a surjective homomorphism. Moreover, $\varphi \circ \psi$ is the
identity of $\tilde{C}_n$ and $\psi\circ \varphi$ is the identity of
$C_n$, thus $\varphi$ and $\psi$ are isomorphisms. This finishes the
proof of the lemma. \hfill $\square$

\vsp\vsp

The group $C_n$ satisfies the hypotheses of Theorem B and has
exactly $2^{n+2}$ Conradian orderings. Indeed, for $1\leq i \leq n$
we let $G^{(i)}=\langle e_1,\ldots ,e_i\rangle\lhd C_n$. Note that
$C_n/G^{(n)}\simeq B(1,3)=\langle f,g\mid fgf^{-1}=g^3\rangle$. We
let $G^{(n+1)}$ be the inverse image of the derived subgroup
$\big(C_n/G^{(n)}\big)^\prime$ under the projection $C_n\to
C_n/G^{(n)}$. Clearly, $C_n/G^{(n+1)}\simeq \Z$, and
$G^{(n+1)}/G^{(n)}\simeq \Z[\frac{1}{3}]$. Moreover, if we let
$G^{(n+2)}=C_n$ then, for $1\leq i \leq n+1$, each quotient
$G^{(i+1)}/G^{(i)}$ is rank-one Abelian. Therefore, the series
$$\{id\} \lhd G^{(1)} \lhd \ldots \lhd\ G^{(n+1)}\lhd
G^{(n+2)}=C_n,$$ is rational. Finally, we have that, for $1\leq i
\leq n$, $G^{(i+2)}/G^{(i)}$ is non Abelian. Thus, the group $C_n$
fits in the classification of groups with only finitely many
$\ce$-orderings. Nevertheless, $C_n$ has $B(1,3)$ as a quotient by a
normal convex subgroup. Since $B(1,3)$ admits uncountably many
left-orderings, the same is true for $C_n$. In fact, it will follow
form Theorem D that no left-ordering on $C_n$ is isolated.







\section{A structure theorem for the space of Conradian orderings}
\label{sec struct thm for C}

As shown by Linnell in \cite{linnell}, the space of left-orderings
of a group is either finite or uncountable; see also \S 2.3 or
\cite{clay 1,crossings}. Although this is no longer true for
bi-orderings, as there are examples of groups having (only)
infinitely countably many bi-orderings; see \cite[\S 6.2]{botto},
\cite{butts}. As announced in the Introduction, in this section we
show

\vsp\vsp

\noindent {\bf Theorem C.} {\em Let $G$ be a $\ce$-orderable group.
If $G$ admits infinitely many $\ce$-orderings, then it has
uncountably many $\ce$-orderings. Moreover, none of these is
isolated in the space of $\ce$-orderings.}

\vsp\vsp

Note that the second statement implies that, if $G$ is countable and
admits infinitely many $\ce$-orderings, then its space of Conradian
orderings is a Cantor set. Note also that for the case of
left-orderings there are group admitting infinitely many
left-orderings together with isolated left-orderings, as it is the
case of braid groups \cite{braids,dd} and the central extensions of
Hecke groups appearing in \cite{navas-hecke}.


\subsection{Finitely many or a Cantor set of Conradian orders}


Let $G$ be a group admitting a $\ce$-ordering $\,\preceq \,$ that is
isolated in the space of $\ce$-orderings. As we have seen at the
beginning of \S \ref{sec}, the series of $\,\preceq$-convex
subgroups must be finite, say
$$\{id\} = G_0 \lhd G_1 \lhd \ldots \lhd G_n = G.$$
Proceeding as in Example \ref{exflip}, any ordering on
$\,G_{i+1}/G_i \,$ may be extended (preserving the set of positive
elements outside of it) to a $\ce$-ordering on $G$. Hence, each
quotient must be rank-one Abelian, so the series above is rational.
We claim that this series of $\,\preceq$-convex subgroups is unique
(hence normal) and that no quotient $\,G_{i+2}/G_i \,$ is Abelian.
In fact, if the series has length 2, then it is normal. Moreover,
since no $\ce$-ordering on a rank-two Abelian group is isolated, we
have that $G_2$ is non Abelian. Then, by Theorem B, the series is
unique. In the general case, we will use induction on the length of
the series. Suppose that every group having an isolated
$\ce$-ordering whose rational series of convex subgroups
$$\{id\} = H_0 \lhd H_1 \lhd \ldots \lhd H_{k}$$
has length $k< n$ admits a unique (hence normal) rational series and
that no quotient $H_{i+2}/H_i$ is Abelian. Let
$$\{id\}=G_0\lhd  \ldots \lhd G_{n-2} \lhd G_{n-1} \lhd G_n=G$$
be a rational series of length $n$ associated to some isolated
$\ce$-ordering $\,\preceq \,$ on $G$. Since $G_{n-1}$ is normal in
$G$, for every $g\in G$, the conjugate series
$$\{id\}=G^g_0\lhd  \ldots \lhd G^g_{n-2} \lhd G^g_{n-1} = G_{n-1}$$
is also a rational series for $G_{n-1}$. Since this series is
associated to a certain isolated $\ce$-ordering, namely the
restriction of $\, \preceq \,$ to $G_{n-1}$, we conclude that it is
unique by the induction hypothesis. Hence the series must coincide
with the original one, or in other words $G_{i}^g = G_{i}$.
Therefore, the series for $G$ is normal. Moreover, every quotient
$G_{i+2}/G_i$ is non Abelian, because if not then $\preceq$ could be
approximated by other $\ce$-orderings on $G$. Thus, by Theorem B,
the rational series for $G$ is unique, and $G$ admits only finitely
many $\ce$-orderings. This completes the proof of Theorem C.


\subsection{An application to orderable nilpotent groups}


In this section we use Theorem C to give a new proof of the
following result which was first proved in \cite{navas}.

\begin{thm}[{\bf Navas}] \label{nilp} \textit{Let $G$ be a torsion-free nilpotent group which
is not rank-one Abelian. Then the space of left-orderings of $G$ is
homeomorphic to a Cantor set.}
\end{thm}

The proof of Theorem \ref{nilp} is a consequence of three facts.

\vsp

The first one is that torsion-free nilpotent groups are
left-orderable. Indeed, as shown in \cite[\S 2.6]{khukhro}, they
admit a filtration
$$ \{id\}\lhd G_1\lhd \ldots \lhd G_{n-1}\lhd G,$$ such that each
quotient $G_i/G_{i-1}$ is torsion-free Abelian (actually they are
bi-orderable). In particular they are $\ce$-orderable. Furthermore,
torsion-free nilpotent group have the (much!) stronger property that
every partial left-ordering can be extended to a (total)
left-ordering; see \cite[Theorem 7.6.4]{botto}.

\vsp

The second fact is a result shown independently by Ault and
Rhemtulla, and appears for instance in \cite[\S 7.5]{botto}. For the
convenience of the reader we give a short proof of this fact:

\begin{prop}\label{nilp_prop} \textit{Every left-ordering on a nilpotent group is Conradian.}

\end{prop}

\noindent \textit{Proof:} Let $(G,\preceq)$ be a left-ordered
nilpotent group. We claim that the action of $G$ on itself has no
crossings.

\vsp

Suppose, by way of a contradiction, that $(f,g,u,v,w)$ is a
crossing. Then, by definition, we have that $f^N(v)\prec w\prec
g^M(u)$. Then, a classical ping-pong argument shows that
$\{f^N,g^M\}$ generates a free semigroup. But this is impossible
since, as it is well known, nilpotent groups can not have free
semigroups (one may think, for instance, in the growth rate of the
subgroup $\langle f^N,g^M\rangle$, see also \cite[exercice
4.47]{book}). $\hfill \square$

\vsp

The third fact is that the only nilpotent, $\ce$-orderable group
with finitely many (Conradian) orderings are the rank-one Abelian
groups. To see this is enough to note that the groups described in
Theorem B, whose rational series has length $2$ or more, have
trivial center.

\vsp

This finishes the proof of Theorem \ref{nilp}.

\begin{rem} A direct consequence of the proof of Proposition \ref{nilp_prop},
is that any left-ordering on a group without free semigroups on two
generators is Conradian. In particular, $H$, the Grigorchuk-Machi
group of intermediate growth from \cite{G-M} (see also
\cite{navas-growth}), is a group admitting only Conradian orderings.
Moreover, from the natural action of $H$ on $\Z^{\N}$, one can be
induced infinitely many left-orderings on it. Therefore, the space
of left-orderings of $H$ is a Cantor set.

\end{rem}


\section{A structure theorem for left-orderings}
\label{sec a struc thm for left ord}

As announced in the Introduction, in this section we explore the
possibility of approximating a given left-ordering by its conjugates
(in the sense of \S \ref{acting on the space}). We will show that in
most cases this can be done. This will give a new proof of Linnell's
result from \cite{linnell}, here stated as Theorem E.

\subsection{Describing the Conradian soul via crossings}
\label{detecting the soul}

The {{\em Conradian soul}} \esp $C_{\preceq}(G)$ of an ordered group
$(G,\preceq)$ corresponds to the maximal (with respect to the
inclusion) subgroup which is convex in $\preceq$, and such that
$\preceq$ restricted to the subgroup is Conradian. This notion was
introduced in \cite{navas}, where a dynamical counterpart in the
case of countable groups was given. To give an analogous
characterization in the general case, we consider the set $S^+$
formed by the elements $h \!\succ\! id$ such that $h \preceq w$ for
every crossing $(f,g,u,v,w)$ satisfying $id \preceq u$. Analogously,
we let $S^{-}$ be the set formed by the elements $h \prec id$ such
that $w \preceq h$ for every crossing $(f,g,u,v,w)$ satisfying $v
\preceq id$. Finally, we let
$$S = \{ id \} \cup S^{+} \cup S^{-}.$$
{\em A priori}, it is not clear that the set $S$ has a nice
structure (for instance, it is not at all evident that it is
actually a subgroup). However, this is largely shown by the next
theorem.

\vspace{0.1cm}

\begin{thm} {\em The Conradian soul of $(G,\preceq)$ coincides
with the set $S$ above.} \label{S=CS}
\end{thm}

\vspace{0.1cm}

Before passing to the proof, we give four general lemmata describing
the flexibility of the concept of crossings for group orderings
(note that the first three lemmata still apply to crossings for
actions on totally ordered spaces). The first one allows us
replacing the ``comparison" element $w$ by its ``images" under
positive iterates of either $f$ or $g$.

\vspace{0.1cm}

\begin{lem} {\em If $(f,g,u,v,w)$ is a crossing,
then $(f,g,u,v,g^n w)$ and $(f,g,u,v,f^n w)$ are also crossings for
every $n \!\in\! \mathbb{N}$.} \label{lema1}
\end{lem}

\noindent {\em Proof:} We will only consider the first 5-uple (the
case of the second one is analogous). Recalling that $g w \succ w$,
for every $n \! \in \! \mathbb{N}$ we have $u \prec w \prec g^n w$;
moreover, $v \succ g^{M+n} u = g^n g^M u \succ g^n w$. Hence, \esp
$u \prec g^n w \prec v$. \esp On the other hand, $f^N v \prec w
\prec g^n w$, while from $g^M u \succ w$ we get $g^{M+n} u \succ g^n
w$. \hfill$\square$

\vspace{0.3cm}

Our second lemma allows replacing the ``limiting" elements $u$ and
$v$ by more appropriate ones.

\vspace{0.15cm}

\begin{lem} {\em Let $(f,g,u,v,w)$ be a crossing. If $f u \succ u$ (resp.
$f u \prec u$) then $(f,g,f^n u,v, w)$ (resp. $(f,g,f^{-n} u, v,
w)$) is also a crossing for every $n > 0$. Analogously, if $g v
\prec v$ (resp. $g v \succ v$), then $(f, g, u, g^n v, w)$ (resp.
$(f,g,u,g^{-n}v,w)$) is also crossing for every $n > 0$.}
\label{lema3}
\end{lem}

\noindent {\em Proof:} Let us only consider the first 5-uple (the
case of the second one is analogous). Suppose that $fu \succ u$ (the
case $fu \prec u$ may be treated similarly). Then $f^n u \succ u$,
which gives $g^M f^n u \succ g^M u \succ w$. To show that $f^n u
\prec w$, assume by contradiction that $f^n u \succeq w$. Then $f^n
u \succ f^N v$, which yields $u \succ f^{N-n} v$, which is absurd.
\hfill$\square$

\vspace{0.3cm}

The third lemma relies on the dynamical insight of the crossing
condition.

\vspace{0.15cm}

\begin{lem}{\em If $(f,g,u,v,w)$ is a crossing, then $(hfh^{-1},
hgh^{-1},hu,hv,hw)$ is also a crossing for every $h \in G$.}
\label{lema2}
\end{lem}

\noindent {\em Proof:} The three conditions to be checked are
nothing but the three conditions in the definition of crossing
multiplied by $h$ by the left. \hfill$\square$

\vspace{0.3cm}

A direct application of the lemma above shows that, if $(f,g,u,v,w)$
is a crossing, then the 5-uples $(f,f^ngf^{-n}, f^n u, f^n v , f^n
w)$ and $(g^nfg^{-n}, g ,g^n u, g^n v, g^n w)$ are also crossings
for every $n \in \mathbb{N}$. This combined with Lemma \ref{lema3}
may be used to show the following.

\vspace{0.15cm}

\begin{lem} {\em If $(f,g,u,v,w)$ is a crossing and $id \preceq h_1 \prec h_2$ are
elements in $G$ such that $h_1 \in S$ and $h_2 \notin S$, then there
exists a crossing
$(\tilde{f},\tilde{g},\tilde{u},\tilde{v},\tilde{w})$ such that $h_1
\prec \tilde{u} \prec \tilde{v} \prec h_2$.} \label{lemapro}
\end{lem}

\noindent {\em  Proof:} Since $id\prec h_2  \notin S$, there must be
a crossing $(f,g,u,v,w)$ such that $id\preceq u \prec w \prec h_2$.
Let $N \in \mathbb{N}$ be such that $f^N v\prec w $. Denote by
$(f,\bar{g},\bar{u},\bar{v},\bar{w})$ the crossing $(f, f^N g
f^{-N}, f^N u, f^N v, f^N w)$. Note that $\bar{v} = f^N v \prec w
\prec h_2$. We claim that $h_1 \preceq \bar{w} = f^N w$. Indeed, if
$f^N u \succ u$, then $f^n u \succ id$, and by the definition of $S$
we must have $h_1 \preceq \bar{w}$. If $f^N u \prec u$, then we must
have $f u \prec u$, so by Lemma \ref{lema3} we know that
$(f,\bar{g},u,\bar{v},\bar{w})$ is also a crossing, which allows
still concluding that $h_1 \preceq \bar{w}$.

Now for the crossing $(f,\bar{g},\bar{u}, \bar{v}, \bar{w})$ there
exists $M \in \mathbb{N}$ such that $\bar{w} \prec \bar{g}^M
\bar{u}$. Let us consider the crossing $(\bar{g}^M f \bar{g}^{-M},
\bar{g}, \bar{g}^M \bar{u}, \bar{g}^M \bar{v}, \bar{g}^M \bar{w})$.
If $\bar{g}^M \bar{v} \prec \bar{v}$, then $\bar{g}^M \bar{v} \prec
h_2$, and we are done. If not, then we must have $\bar{g} \bar{v}
\succ \bar{v}$. By Lemma \ref{lema3}, $(\bar{g}^M f \bar{g}^{-M},
\bar{g}, \bar{g}^M \bar{u}, \bar{g}^M \bar{v}, \bar{w})$ is still a
crossing, and since $\bar{v} \prec h_2$, this concludes the
proof.\hfill$\square$

\vspace{0.38cm}

\noindent{\em Proof of Theorem} \ref{S=CS}. The proof is divided
into several steps.

\vsp\vsp

\noindent {\underbar{Claim 0.}} The set $S$ is convex.

\vsp

This follows directly from the definition of $S$.

\vsp\vsp

\noindent {\underbar{Claim 1.}} If $h$ belongs to $S$, then $h^{-1}$
also belongs to $S$.

\vsp

Assume that $h \in S$ is positive and $h^{-1}$ does not belong to
$S$. Then there exists a crossing $(f,g,u,v,w)$ so that $h^{-1}
\prec w \prec v \preceq id$.

We first note that, if $h^{-1} \preceq u$, then after conjugating by
$h$ as in Lemma \ref{lema2}, we get a contradiction because
$(hgh^{-1}, hfh^{-1}, hu, hv, hw)$ is a crossing with \esp $id
\preceq hu $ \esp and \esp $hw \prec hv \preceq h$. To reduce the
case $h^{-1} \succ u$ to this one, we first use Lemma \ref{lema2}
and we consider the crossing $(g^Mfg^{-M}, g, g^M u, g^M v, g^M w)$.
Since \esp $h^{-1} \prec w \prec g^M u \prec g^M w \prec g^M v$,
\esp if $g^M v \prec v$ then we are done. If not, Lemma \ref{lema3}
shows that $(g^Mfg^{-M}, g, g^M u , g^M v, w)$ is also a crossing,
which still allows concluding.

In the case where $h \in S$ is negative we proceed similarly but we
conjugate by $f^N$ instead of $g^M$. Alternatively, since $id \in S$
and $id\prec h^{-1}$, if we suppose that $h^{-1}\notin S$ then Lemma
\ref{lemapro} provides us with a crossing $(f,g,u,v,w)$ such that
$id\prec u\prec w \prec v \prec h^{-1}$, which gives a contradiction
after conjugating by $h$.

\vsp\vsp

\noindent {\underbar{Claim 2.}} If $h$ and $\bar{h}$ belong to $S$,
then $h\bar{h}$ also belongs to $S$.

\vsp

First we show that for every positive elements in $S$, their product
still belongs to $S$. (Note that, by Claim 1, the same will be true
for products of negative elements in $S$.) Indeed, suppose that
$h,\bar{h}$ are positive elements, with $h \in S$ but $ h \bar{h}
\notin S$. Then, by Lemma \ref{lemapro} we may produce a crossing
$(f,g,u,v,w)$ such that $h \prec u \prec v \prec h \bar{h}$. After
conjugating by $h^{-1}$ we obtain the crossing
$(h^{-1}fh,h^{-1}gh,h^{-1}u,h^{-1}v,h^{-1}w)$ satisfying $id \prec
h^{-1}u \prec h^{-1} w \prec \bar{h}$, which shows that $\bar{h}
\notin S$.

Now, if $h \prec id \prec \bar{h}$, then $h \prec h\bar{h}$. Hence,
if $h\bar{h}$ is negative, then the convexity of $S$ gives $h\bar{h}
\in S$. If $h\bar{h}$ is positive, then $\bar{h}^{-1}h^{-1}$ is
negative, and since $\bar{h}^{-1} \prec \bar{h}^{-1} h^{-1}$, the
convexity gives again that $\bar{h}^{-1}h^{-1}$, and hence
$h\bar{h}$, belongs to $S$. The remaining case $\bar{h} \prec
id\prec h$ may be treated similarly.

\vsp \vsp

\noindent {\underbar{Claim 3.}} The subgroup $S$ is Conradian.

\vsp

In order to apply Theorem A, we need to show that there are no
crossings in $S$. Suppose by contradiction that $(f,g,u,v,w)$ is a
crossing such that $f,g,u,v,w$ all belong to $S$. If $id \preceq w$,
then, by Lemma \ref{lema2}, we have that $(g^n f g^{-n}, g, g^n u,
g^n v, g^n w)$ is a crossing. Taking $n = M$ so that $g^M u \succ
w$, this gives a contradiction to the definition of $S$ because $id
\preceq w \prec g^M u \prec g^M w \prec g^M v \in S$. The case $w
\preceq id$ may be treated in an analogous way by conjugating by
powers of $f$ instead of $g$.

\vsp\vsp

\noindent {\underbar{Claim 4.}} The subgroup $S$ is maximal among
$\preceq$-convex, $\preceq$-Conradian subgroups.

\vsp

Indeed, if $C$ is a subgroup strictly containing $S$, then there is
a positive $h$ in $C \setminus S$. By Lemma \ref{lemapro}, there
exists a crossing $(f,g,u,v,w)$ such that $id\prec u \prec w \prec v
\prec h$. If $C$ is convex, then $u,v,w$ belong to $C$. To conclude
that $C$ is not Conradian, it suffices to show that $f$ and $g$
belong to $C$.

Since $id\prec u$, we have either $id \prec g \prec g u \prec v$ or
$id \prec g^{-1} \prec g^{-1} u \prec v$. In both cases, the
convexity of $C$ implies that $g$ belongs to $C$. On the other hand,
if $f$ is positive, then from $f^N \prec f^N v \prec w$ we get $f
\in C$, whereas in the case of a negative $f$ the inequality
$id\prec u$ gives $id\prec f^{-1} \prec f^{-1} u \prec v$, which
still shows that $f \in C$. $\hfill \square$


\subsection{Approximating a left-orderings by its conjugates}


Recall from \S \ref{pre-1} that the { {\em positive cone}} of a
left-ordering $\preceq$ in $\mathcal{LO}(G)$ is the set $P$ of its
positive elements. Because of the left invariance, $P$ completely
determines $\preceq$. The {{\em conjugate}} of $\preceq$ by $h \in
G$ is the left-ordering $\preceq_h$ having positive cone $hPh^{-1}$.
In other words, $g \succ_h id$ holds if and only if $h g h^{-1}
\succ id$. We will say that $\preceq$ may be approximated by its
conjugates if it is an accumulation point of its set of conjugates.

\vspace{0.2cm}

\begin{thm}\label{trivial soul} {\em Suppose  $(G,\preceq)$  is a non trivial
left-ordered group such that it has trivial Conradian soul. Then
$\preceq$ may be approximated by its conjugates.} \label{primero}
\end{thm}

\noindent {\em Proof:} Let $f_1 \prec f_2 \prec \ldots \prec f_k$ be
finitely many positive elements in $G$. We need to show that there
exists a conjugate of $\preceq$ that is different from $\preceq$ but
for which all the $f_i$'s are still positive. Since $id \!\in\!
C_{\preceq}(G)$ and $f_1 \notin C_{\preceq}(G)$, Theorem \ref{S=CS}
and Lemma \ref{lemapro} imply that there is a crossing $(f,g,u,v,w)$
such that $id\prec u \prec v \prec f_1$. Let $M,N$ in $\mathbb{N}$
be such that $ f^N v \prec w \prec g^M u$. We claim that
$id\prec_{v^{-1}} f_i$ and $id\prec_{w^{-1}} f_i$ for $1\leq i \leq
k$, but $g^M f^N \prec_{v^{-1}} id$ and $g^M f^N \succ_{w^{-1}} id$.
Indeed, since $id \prec v \prec f_i$, we have $v\prec f_i\prec f_i
v$, thus $id \prec v^{-1} f_i v$. By definition, this means that
$f_i\succ_{v^{-1}} id$. The inequality $f_i \succ_{w^{-1}} id$ is
proved similarly. Now note that $g^M f^N v \prec g^M w \prec v$, and
so $g^M f^N \prec_{v^{-1}} id$. Finally, from $g^Mf^N w \succ g^M u
\succ w $ we get $g^Mf^N \succ_{w^{-1}} id$.

Now the preceding relations imply that the $f_i$'s are still
positive for both $\preceq_{v^{-1}}$ and $\preceq_{w^{-1}}$, but at
least one of these left-orderings is different from $\preceq$. This
concludes the proof. \hfill$\square$

\vspace{0.3cm}

Based on the work of Linnell \cite{linnell}, it is shown in
\cite[Proposition 4.1]{navas} that no Conradian ordering is an
isolated point of the space of left-orderings of a group having
infinitely many left-orderings. Note that this result also follows
as a combination of Theorem C and Theorem D. Together with Theorem
\ref{primero}, this shows the next proposition by means of the
convex extension procedure ({\em c.f.}, Corollary \ref{nice coro}).

\vsp

\begin{prop} {\em Let $G$ be a left-orderable group. If $\preceq$
is an isolated point of $\mathcal{LO} (G)$, then its Conradian soul
is nontrivial and has only finitely many left-orderings.}
\label{linelito}
\end{prop}

\vsp

As a consequence of Tararin's theorem, here Theorem \ref{teo T}, the
number of left-orderings on a left-orderable group having only
finitely many left-orderings is a power of 2. Moreover, all of these
left-orderings are necessarily Conradian; see Corollary
\ref{finitely leftorders}. By the preceding theorem, if $\preceq$ is
an isolated point of the space of left-orderings $\mathcal{LO} (G)$,
then its Conradian soul admits $2^n$ different left-orderings for
some $n \geq 1$, all of them Conradian. Let
$\{\preceq_1,\preceq_2,\ldots,\preceq_{2^n}\}$ be these
left-orderings, where $\preceq_1$ is the restriction of $\preceq$ to
its Conradian soul. Since $C_{\preceq}(G)$ is $\preceq$-convex, each
$\preceq_j$ induces a left-ordering $\preceq^j$ on $G$, namely the
convex extension of $\preceq_j$ by $\preceq$. (Note that $\preceq^1$
coincides with $\preceq$.) The Lemma below appears in \cite{navas}.
For the reader's convenience we provide a proof of this fact.

\begin{lem}{\em With the notations above, all the left-orderings
$\preceq^j$ share the same Conradian soul.}
\end{lem}

\noindent {\em Proof:} Consider the left-ordering $\preceq^j$. Since
$\preceq^j$ restricted to $C_\preceq(G)$ is Conradian, and
$C_\preceq(G)$ is convex in $\preceq^j$, we only need to check that
$C_{\preceq^j}(G)\subseteq C_\preceq(G)$. Let $G^*$ be any
$\preceq^j$-convex subgroup strictly containing $C_\preceq(G)$. We
claim that $G^*$ is also $\preceq$-convex. Indeed, since $\preceq^j$
coincides with $\preceq$ outside $C_\preceq(G)$, we have that for
any $f\notin C_\preceq(G)$, $id\prec f$ if and only if $id \prec^j
f$; see for instance \S \ref{basic contr}. In particular,  $id\prec
h \prec g$ for $g \in G^*$ and $g\notin C_\preceq(G)$, implies
$h\prec^j g$, hence $h\in G^*$, and the claim follows.

\vsp

Since $G^*$ is $\preceq$-convex and strictly contains
$C_\preceq(G)$, we have that there are $f,g$ in $G^*$ such that
$id\prec f\prec g$ and $fg^n\prec g$ for all $n\in \N$. Clearly
$g\notin C_\preceq(G)$. We claim that for all $n\in \N$, the element
$g^{-1}fg^n$ does not belong to $C_\preceq(G)$. Indeed, if
$g^{-1}fg^n\in C_\preceq(G)$, then $(g^{-1}fg^{n})^{-1}\prec g$,
which implies that $id\prec gfg^{n+1}$ contrary to our choice of $f$
and $g$.

\vsp

If it is the case that $f\notin C_\preceq(G)$, then we are done.
Indeed, since $id\prec^j f$, $id \prec^j g$ and $g^{-1}fg^n\prec^j
id$ for all $n\in \N$, we have that $\preceq^j$ restricted to $G^*$
is not Conradian. In the case that $f\in C_\preceq(G)$, we let
$h=fg$. Note that $h\notin C_\preceq(G)$ and that $h\succ^j id$.
Moreover, as before, we have that $g^{-1}hg^{n}\prec id$ for all
$n\in \N$ and $g^{-1}hg^n\notin C_\prec(G)$. This shows that
$\preceq^j$ restricted to $G^*$ is not Conradian. \hfill$\square$

\vsp

Below, assume that $\preceq$ is not Conradian.

\begin{thm} {\em  With the notation above, at least one of the left-orderings
$\preceq^j$ is an accumulation point of the set of conjugates of
$\preceq$.} \label{final}
\end{thm}

Before proving this theorem, we immediately state

\begin{cor} {\em At least one of the left-orderings
$\preceq^j$ is approximated by its conjugates.} \label{finalito}
\end{cor}

\noindent {\em Proof:} Assuming Theorem \ref{final}, we have
$\preceq^k \in\! acc (orb(\preceq^1))$ for some $k \in
\{1,\ldots,2^n\}$. Theorem \ref{final} applied to this $\preceq^k$
instead of $\preceq$ shows the existence of $k' \in
\{1,\ldots,2^n\}$ so that $\preceq^{k'} \in acc (orb(\preceq^k))$,
and hence $\preceq^{k'} \in acc(orb(\preceq^1))$. If $k'$ equals
either $1$ or $k$ then we are done; if not, we continue arguing in
this way... In at most $2^n$ steps we will find an index $j$ such
that $\preceq^j \in acc(orb(\preceq^j))$. \hfill$\square$

\vspace{0.2cm}

Theorem \ref{final} will follow from the next

\vsp

\begin{prop} {\em Given an arbitrary finite family
$\mathcal{G}$ of $\preceq$-positive elements in $G$, there exists $h
\in G$ and $id \prec \bar{h} \notin C_{\preceq} (G)$ such that $id
\prec h^{-1} f h \notin C_{\preceq} (G)$ for all $f \in \mathcal{G}
\setminus C_{\preceq} (G)$, but $id \succ h^{-1}\bar{h}h \notin
C_{\preceq} (G)$.} \label{a-probar}
\end{prop}

\vsp

\noindent{\em Proof of Theorem \ref{final}  from Proposition
\ref{a-probar}:} Let us consider the directed set formed by the
finite sets $\mathcal{G}$ of $\preceq$-positive elements. For each
such a $\mathcal{G}$, let $h_{\mathcal{G}}$ and
$\bar{h}_{\mathcal{G}}$ be the elements in $G$ provided by
Proposition \ref{a-probar}. After passing to subnets of
$(h_{\mathcal{G}})$ and $(\bar{h}_{\mathcal{G}})$ if necessary, we
may assume that the restrictions of $\preceq_{h_{\mathcal{G}}^{-1}}$
to $C_{\preceq}(G)$ all coincide with a single $\preceq_j$. Now the
properties of $h_{\mathcal{G}}$ and $\bar{h}_{\mathcal{G}}$ imply:

\vsp

\noindent -- $f \succ^j id$ \esp \esp and \esp \esp $f \esp
(\succ^j)_{h_{\mathcal{G}}^{-1}} \esp id$ \esp \esp for all \esp $f
\in \mathcal{G} \setminus C_{\preceq} (G)$,

\vsp

\noindent -- $\bar{h}_{\mathcal{G}} \succ^j id$, \esp \esp  but \esp
\esp \esp $\bar{h}_{\mathcal{G}} \esp \esp
(\prec^j)_{h_{\mathcal{G}}^{-1}} \prec id$.

\vsp

\noindent This shows Theorem \ref{final}. $\hfill\square$

\vspace{0.5cm}

For the proof of Proposition \ref{a-probar} we will use three
general lemmata.

\begin{lem} {\em For every $id\prec c \notin C_{\preceq} (G)$
there is a crossing $(f,g,u,v,w)$ such that $u,v,w$ do not belong to
$C_{\preceq}(G)$ and $id\prec u\prec w \prec v \prec c$.}
\label{primer-lema}
\end{lem}

\noindent {\em Proof:} By Theorem \ref{S=CS} and Lemma
\ref{lemapro}, for every $id\preceq s\in C_{\preceq} (G)$ there
exists a crossing $(f,g,u,v,w)$ such that $s\prec u\prec w\prec v
\prec c$. Clearly, $v$ does not belong to $ C_{\preceq} (G)$. The
element $w$ is also outside $ C_{\preceq} (G)$, since in the other
case the element $a = w^2$ would satisfy $w \prec a \in C_{\preceq}
(G)$, which is absurd. Taking $M > 0$ so that $g^M u \succ w$, this
gives \esp $g^M u \notin C_{\preceq} (G)$, \esp $g^M w \notin
C_{\preceq} (G)$, \esp and $g^M v \notin C_{\preceq} (G)$. \esp
Consider the crossing $(g^Mfg^{-M}, g, g^M u,g^M v, g^M w)$. If $g^M
v\prec v$, then we are done. If not, then $gv \succ v$, and Lemma
\ref{lema3} ensures that $(g^Mfg^{-M}, g, g^M u,  v, g^Mw)$ is also
a crossing, which still allows concluding. \hfill$\square$


\begin{lem} {\em Given $id \prec c \notin C_{\preceq} (G)$ there exists
$id \prec a \notin C_{\preceq}(G)$ (with $a \prec c$) such that, for
all $id \preceq b \preceq a$ and all \esp $\bar{c} \succeq c$, one
has $id \prec b^{-1} \bar{c} b \notin C_{\preceq} (G)$.}
\label{segundo-lema}
\end{lem}

\noindent {\em Proof:} Let us consider the crossing $(f,g,u,v,w)$
such that $id\prec u\prec w \prec v \prec c$ and such that $u,v,w$
do not belong to $ C_{\preceq} (G)$. We affirm that the lemma holds
for $a = u$ (actually, it holds for $a = w$, but the proof is
slightly more complicated). Indeed, if $id \preceq b \preceq u$,
then from $b \preceq u \prec v \prec \bar{c}$ we get $id \preceq
b^{-1}u \prec b^{-1}v \prec b^{-1} \bar{c}$, and thus the crossing
$(b^{-1}fb,b^{-1}gb,b^{-1}u,b^{-1}v,b^{-1}w)$ shows that $b^{-1}
\bar{c} \notin C_{\preceq} (G)$. Since $id \preceq b$, we conclude
that $id \prec b^{-1} \bar{c} \preceq b^{-1} \bar{c} b$, and the
convexity of $ S$ implies that $b^{-1} \bar{c} b \notin C_{\preceq}
(G)$. \hfill$\square$


\begin{lem} {\em For every $g\in G$ the set $g \esp C_{\preceq} (G)$
is convex. Moreover, for every crossing $(f,g,u,v,w)$ one has $\,u
C_{\preceq} (G) < w C_{\preceq} (G) < v C_{\preceq} (G)$, in the
sense that \esp $uh_1 \prec wh_2 \prec vh_3$ \esp for all
$h_1,h_2,h_3$ in $ C_{\preceq} (G)$.} \label{lema4}\end{lem}

\noindent {\em Proof:} The verification of the convexity of $g
C_{\preceq} (G)$ is straightforward. Now suppose that $uh_1 \succ
wh_2$ for some $h_1,h_2$ in $ C_{\preceq} (G)$. Then since $u \prec
w$, the convexity of both left classes $u C_{\preceq} (G)$ and $w
C_{\preceq} (G)$ gives the equality between them. In particular,
there exists $h \in C_{\preceq} (G)$ such that $uh = w$. Note that
such an $h$ must be positive, so that $id \prec h = u^{-1} w$. But
since $(u^{-1}fu, u^{-1}gu, id,u^{-1}v, u^{-1}w)$ is a crossing,
this contradicts the definition of $ C_{\preceq} (G)$. Showing that
$w C_{\preceq} (G) \prec v C_{\preceq} (G)$ is
similar.$\hfill\square$

\vs

\noindent{\em Proof of Proposition \ref{a-probar}:} Let us label the
elements of $\mathcal{G} \!=\! \{f_1,\ldots,f_r\}$ so that $f_1
\prec \ldots \prec f_r$, and let $k$ be such that $f_{k-1} \in
C_{\preceq} (G)$ but $f_{k} \notin C_{\preceq} (G)$. Recall that, by
Lemma \ref{segundo-lema}, there exists $id \prec a \notin
C_{\preceq} (G)$ such that, for every $id \preceq b \preceq a$, one
has $id \prec b^{-1} f_{k+j} b \notin C_{\preceq} (G)$ for all $j
\geq 0$. We fix a crossing $(f,g,u,v,w)$ such that $id\prec u\prec
v\prec a$ and $u\notin C_{\preceq} (G)$. Note that the conjugacy by
$w^{-1}$ gives the crossing $(w^{-1}fw,
w^{-1}gw,w^{-1}u,w^{-1}v,id)$.

\vsp\vsp

\noindent{\underbar{Case 1.}} One has $w^{-1}v \preceq a$.

In this case, the proposition holds for \esp $h=w^{-1}v$ \esp and
\esp $\bar{h}=w^{-1}g^{M+1}f^N w$. \esp To show this, first note
than neither $w^{-1}gw$ nor $w^{-1}fw$ belong to $C_{\preceq} (G)$.
Indeed, this follows from the convexity of $C_{\preceq} (G)$ and the
inequalities \esp $w^{-1}g^{-M}w \prec w^{-1}u \notin C_{\preceq}
(G)$ \esp and $w^{-1} f^{-N} w \succ w^{-1} v \notin C_{\preceq}
(G)$. \esp We also have $id\prec w^{-1}g^{M}f^{N}w$, and hence
$id\prec w^{-1}gw \prec w^{-1} g^{M+1}f^{N}w$, which shows that
$\bar{h} \notin C_{\preceq} (G)$. On the other hand, the inequality
$w^{-1} g^{M+1}f^{N}w (w^{-1} v) \prec w^{-1}v$ reads as
$h^{-1}\bar{h}h\prec id$. Finally, Lemma \ref{lema1} applied to the
crossing $(w^{-1}fw, w^{-1}gw,w^{-1}u,w^{-1}v,id)$ shows that
$(w^{-1} fw, w^{-1}gw, w^{-1} u, w^{-1} v, w^{-1} g^{M+n}f^{N}w)$ is
a crossing for every $n > 0$. For $n \geq M$ we have $w^{-1}
g^{M+1}f^{N}w (w^{-1}v)\prec w^{-1} g^{M+n}f^{N}w $. Since $w^{-1}
g^{M+n}f^{N}w \prec w^{-1} v$, Lemma \ref{lema4} easily implies that
$w^{-1} g^{M+1}f^{N}w (w^{-1}v) C_{\preceq} (G) \prec
w^{-1}vC_{\preceq} (G)$, that is, $h^{-1}\bar{h}h\notin C_{\preceq}
(G)$.

\vsp\vsp

\noindent{\underbar{Case 2.}} One has $a\prec w^{-1}v$, but $w^{-1}
g^m w \preceq a$ for all $m > 0$.

We claim that, in this case, the proposition holds for $h = a$ and
$\bar{h}=w^{-1}g^{M+1}f^N w$. This may be checked in the very same
way as in Case 1 by noticing that, if $a\prec w^{-1}v$ but
$w^{-1}g^m w \succeq a$ for all $m > 0$, then $(w^{-1}fw,
w^{-1}gw,w^{-1}u,a,id)$ is a crossing.

\vsp\vsp

\noindent{\underbar{Case 3.}} One has $a\prec w^{-1}v$ and $w^{-1}
g^m w\succ a$ for some $m > 0$. (Note that the first condition
follows from the second one.)

We claim that, in this case, the proposition holds for $h=a$ and
$\bar{h}=w \notin C_{\preceq} (G)$. Indeed, we have $g^{m}w\succ ha$
(and $w\prec ha$), and since $g^{m}w\prec v\prec a$, we have
$wa\prec a$, which means that $h^{-1}\bar{h}h \prec id$. Finally,
from Lemmas \ref{lema1} and \ref{lema4} we get \esp $wa C_{\preceq}
(G)\preceq g^m w C_{\preceq} (G) \prec v C_{\preceq} (G) \preceq a
C_{\preceq} (G)$. \esp This implies that \esp $a^{-1} w a
C_{\preceq} (G) \prec C_{\preceq} (G)$, \esp which means that
$h^{-1} \bar{h} h \notin C_{\preceq} (G)$. $\hfill\square$


\subsection{Finitely many or uncountably many left-orderings}


The goal of this final short section is to use the previously
developed ideas to give an alternative proof of the following result
due to Linnell; see \cite{linnell}.

\vspace{0.2cm}

\vsp\vsp

\noindent {\bf Theorem E (Linnell).} {\em If the space of orderings
of an orderable group is infinite, then it is uncountable.}

\vsp\vsp

\noindent {\em Proof:} Let us fix an ordering $\preceq$ on an
orderable group $G$. We need to analyze two different cases.

\vspace{0.2cm}

\noindent{\underbar{Case 1.}} The Conradian soul of $C_{\preceq}(G)$
is nontrivial and has infinitely many left-orderings.

\vspace{0.1cm}

This case was settled in \cite{navas} (see Proposition 4.1 therein)
using ideas going back to Zenkov \cite{zenkov}. Alternatively we can
use Theorem C and D to conclude that $C_\preceq(G)$ has no isolated
left-orderings, so it is uncountable. By proposition \ref{nice
prop}, the same is true for the space of left-orderings of $G$.

\vspace{0.2cm}

\noindent{\underbar{Case 2.}} The Conradian soul of $C_{\preceq}(G)$
has only finitely many orderings.

\vspace{0.1cm}

If $\preceq$ is Conradian, then $G = C_{\preceq} (G)$ has finitely
many orderings. If not, then Theorems \ref{primero} and \ref{final}
imply that there exists an ordering $\preceq^*$ on $G$ which is an
accumulation point of its conjugates. The closure in
$\mathcal{LO}(G)$ of the set of conjugates of $\preceq^*$ is then a
compact set without isolated points. By a well-known fact in General
Topology, such a set must be uncountable. Therefore, $G$ admits
uncountably many orderings. \hfill$\square$


\chapter[{\small Left-orders on groups with finitely many $\ce$-orders}]{Left-orders on groups with finitely many $\ce$-orders}

The main result of this chapter is motivated by the following

\vsp

\begin{qs} Is it true that for left-orderable, solvable groups,
having an isolated left-ordering is equivalent to having only
finitely many left-orderings?
\end{qs}

Indeed, as far as the author knows, the only examples of groups
having infinitely many left-orderings together with isolated
left-orderings are braid groups, \cite{braids,dd}, and the groups
introduced in \cite{navas-hecke}. Both families of groups are not
solvable (actually they contain free subgroups). On the other hand,
the dichotomy holds for  nilpotent groups; see Theorem \ref{nilp}.

\vsp

Here we focus on a (very) restricted subfamily of solvable groups,
namely, groups having only finitely many Conradian orderings. These
groups are described in Theorem B. We show

\vsp\vsp

\noindent {\bf Theorem D.} \textit{If a $\ce$-orderable group has
only finitely many $\ce$-orderings, then its space of left-orderings
is either finite or homeomorphic to the Cantor set.}

\vsp\vsp

As shown in Theorem B, a group with finitely many Conradian
orderings admits a unique rational series. Therefore, it is
countable, so $\mathcal{LO}(G)$ is metrizable. Thus, in order to
prove Theorem D, we need to show that no left-ordering of $G$ is
isolated, unless there are only finitely many of them.

\vsp

We proceed by induction on the length of the rational series. In \S
3.1 we explore the case $n=2$. In this case we will give an explicit
description of $\mathcal{LO}(G)$. This extends \cite{rivas}, where
the space of orderings of the Baumslag-Solitar $B(1,\ell)$,
$\ell\geq 2$, is described. In \S 3.2.1, we obtain some technical
lemmata partially describing the inner automorphisms of a group with
a finite number of Conradian orderings. As a result we show that the
maximal convex subgroup of $G$ (with respect to a $\ce$-ordering) is
a group that fits into the classification made by Tararin, {\em
i.e.} a group with only finitely many left-orderings (a {\em Tararin
group}, for short); see Theorem \ref{teo T}. Finally, in \S 3.2.2,
we prove the inductive step. Section 3.2.3 is devoted to the
description of an illustrative example.


\section{The metabelian case}
\label{sec metabelian case}

Throughout this section, $G$ will denote a left-orderable, non
Abelian group with rational series of length $2$:
$$\{id\}=G_0\lhd G_1\lhd G_2=G.$$

\vsp

If $G$ is not bi-orderable, then for the rational series above the
quotient $G_2/G_0=G$ is non bi-orderable. Therefore, $G$ fits into
the classification made by Tararin, here Theorem \ref{teo T}, so it
has only finitely many left-orderings.

\vsp

For the rest of this section, we will assume that $G$ is not a
Tararin group, hence $G$ is bi-orderable. We have

\begin{lem} \label{lema 1}\textit{The group $G$ satisfies $G/G_1\simeq\Z$.}
\end{lem}

\noindent \textit{Proof:} Indeed, consider the action by conjugation
$\alpha: G/G_1\to Aut(G_1)$ given by $\alpha(gG_1)(h)=ghg^{-1}$.
Since $G$ is non Abelian, we have that this action is nontrivial,
{\em i.e.} $Ker(\alpha)\not=G/G_1$. Moreover, $Ker(\alpha)=\{id\}$,
since in the other case, as $G/G_1$ is rank-one Abelian, we would
have that $(G/G_1)/Ker(\alpha)$ is a torsion group. But the only
nontrivial, finite order automorphism of $G_1$ is the inversion,
which implies that $G$ is non bi-orderable, thus a Tararin group.

\vsp

The following claim is elementary and it we leave its proof to the
reader.

\vsp

\noindent \underbar{Claim.} If $G$ is a torsion-free, rank-one
Abelian group such that $G \not\simeq  \Z$, then for any $g\in G$,
there is an integer $n>1$ and $g_n\in G$ such that $g_n^n=g$.

\vsp

Now take any $b\in G\setminus G_1$ so that $\alpha (bG_1)$ is a
nontrivial automorphism of $G_1$. Since $G_1$ is rank-one Abelian,
for some positive $r=p/q\in \Q$, $r\not=1$, we must have that
$bab^{-1}=a^r$ for all $a\in G_1$. Suppose that $G/G_1\not\simeq
\Z$. By the previous claim, we have a sequence of increasing
integers $(n_1,n_2\ldots)$ and a sequence $(g_1,g_2,\ldots)$ of
elements in $G/G_1$ such that $g_i^{n_i}=bG_1$. But clearly this can
not happen since for $g_i$ we have that $g_iag_i^{-1}=a^{r_i}$,
where $r_i$ is a rational such that $r_i^{n_i}=r$, which is
impossible. (In other words, given $r$ we have found, among the
rational numbers, an infinite collection of $r_i$ solving the
equation $x^{n_i}-r=0$, but, by the Rational Roots Theorem or
Rational Roots Test this can not happen; see for instance
\cite[Proposition 5.1]{morandi}.) This finishes the proof of Lemma
1.1. $\hfill \square$

\vs

\begin{lem}\label{lema afin} {\em The group $G$ embeds in $Af_+(\R)$, the group of (orientation-preserving)
affine homeomorphism of the real line.}
\end{lem}

\noindent \textit{Proof:} We first embed $G_1$. Fix $a\in G_1$,
$a\not=id$. Define $\varphi_a:G_1\to Af_+(\R)$ by declaring
$\varphi_a(a)(x)=x+1$, and if $a^\prime \in G_1$ is such that
$(a^\prime)^q=a^p$, we declare $\varphi_a(a^\prime)(x)=x+p/q$.
Showing that $\varphi_a$ is an injective homomorphism is routine.

\vsp

Now let $b\in G$ be such that $\langle bG_1\rangle=G/G_1$. Let
$1\not=r\in \Q$ such that $ba^\prime b^{-1}=(a^\prime)^r$ for every
$a^\prime\in G_1$. Since $G$ is bi-orderable we have that $r>0$, and
changing $b$ by $b^{-1}$ if necessary, we may assume that $r>1$.
Then, given $w\in G$, there is a unique $n\in \Z$ and a unique
$\overline{w}\in G_1$ such that $w=b^n\overline{w}$.

\vsp

Define $\varphi_{b,a}:G\to Af_+(\R)$ by letting
$\varphi_{b,a}(b^n\overline{w}) =H^{(n)}_r\circ
\varphi_a(\overline{w}),$ where $H_r(x)=rx$ , and $H_r^{(n)}$ is the
$n$-th iterate of $H_r$ (by convention $H_r^{(0)}(x)=x$). We claim
that $\varphi_{b,a}$ is an injective homomorphism.

\vsp

Indeed, let $w_1,w_2$ in $G$, $w_1=b^{n_1}\overline{w}_1$,
$w_2=b^{n_2}\overline{w}_2$. Let $r_1\in \Q$ be such that
$\varphi_a(\overline{w}_1)(x)=x+r_1$. Then $ H_r^{(n)}\circ
\varphi_a(b^{-n}\overline{w}_1b^{n})(x)=H_r^{(n)} \circ
\varphi_a(\overline{w}_1^{(1/r)^{n}})(x)=r^n(x+
(1/r)^{n}r_1)=\varphi_a(\overline{w}_1)\circ H_r^{(n)}(x),$ for all
$n\in \Z$. Thus
\begin{eqnarray*}
\varphi_{b,a}(w_1w_2)
 &=& \varphi_{b,a}(b^{n_1}b^{n_2} \; b^{-n_2}\overline{w}_1 b^{n_2} \overline{w}_2)=H_r^{(n_1)} \circ H_r^{(n_2)}\circ \varphi_a(b^{-n_2}\overline{w}_1 b^{n_2})
 \circ\varphi_a( \overline{w}_2)\\
 &=& H_r^{(n_1)} \circ \varphi_a(\overline{w}_1)\circ H_r^{(n_2)}
 \circ\varphi_a( \overline{w}_2)=\varphi_{b,a}(w_1)\circ
 \varphi_{b,a}(w_2),
\end{eqnarray*}
which shows that $\varphi$ is a homomorphism. To see that it is
injective, suppose that
$\varphi(w_1)(x)=\varphi(b^{n_1}\overline{w}_1)(x) =
r^n\,x+r^nr_1=x$ for all $x\in \R$. Then $n=0$ and $r_1=0$, showing
that $w_1=id$. This finishes the proof of Lemma 1.2 $\hfill\square$

\vs

Once the embedding $\varphi=\varphi_{b,a}:G\to Af_+(\R)$ is fixed,
we can associate to each irrational number $\varepsilon$ an {\em
induced left-ordering} $\preceq_{\varepsilon}$ on $G$ whose set of
positive elements is $\{g \!\in\! G \mid \esp \varphi
(g)(\varepsilon)
> \varepsilon \}$. When $\varepsilon$ is rational, the
preceding set defines only a partial ordering. However, in this case
the stabilizer of the point $\varepsilon$ is isomorphic to
$\mathbb{Z}$, hence this partial ordering may be completed to two
total left-orderings $\preceq_{\varepsilon}^+$ and
$\preceq_{\varepsilon}^{-}$. These orderings were introduced by
Smirnov in \cite{smirnov}. Once the representation $\varphi$ is
fixed, we call these orderings, together with its corresponding
reverse orderings, Smirnov-type orderings.

\vsp

Besides the Smirnov-type orderings on $G$, there are four Conradian
(actually bi-invariant) orderings. Since $G_1$ is always convex in a
Conradian ordering, the sign of $b^na^s\in G$, $n\not=0$, depends
only on the sign of $b$ and $n$. Then it is not hard to check that
the four Conradian orderings are the following:

\vsp\vsp

1) $\preceq_{C_1}$, defined by $id \prec_{C_1} b^n a^s$ ($n\in \Z$,
$s\in \Q$) if and only if either $n\geq 1$, or $n=0$ and $s>0$.

\vsp

2) $\preceq_{C_2}$, defined by $id \prec_{C_2} b^na^s$ if and only
if either $n\leq -1$, or $n=0$ and $s>0$.

\vsp

3) $\preceq_{C_3}=\overline{\preceq}_{C_1} $ (the reverse ordering
of $\preceq_{C_1}$).

\vsp

4) $\preceq_{C_4}=\overline{\preceq}_{C_2}$.

\begin{prop} \label{C y S} {\em Let $U\subseteq \mathcal{LO}(G)$ be the set consisting of the four Conradian
orderings together with the Smirnov-type orderings. Then any
ordering in $U$ is non isolated in $U$.}
\end{prop}

\noindent {\em Proof:} We first show that the Conradian orderings
are non isolated.

\vsp

Indeed, we claim that $\preceq_\varepsilon \to \preceq_{C_1}$ as
$\varepsilon\to \infty$. To show this, it suffices to show that any
positive element in the $\preceq_{C_1}$ ordering becomes
$\preceq_\varepsilon$-positive for any $\varepsilon$ large enough.

\vsp

By definition of $\preceq_\varepsilon$, we have that
$id\prec_\varepsilon b^na^s$ if and only if
$r^n(\varepsilon+s)=\varphi(b^na^s)(\varepsilon)>\varepsilon$, where
$r>1$. Now, assume that $id\prec_{C_1}b^na^s$. If $n=0$, then $s>0$
and $\varepsilon+s>\varepsilon$. If $n\geq 1$, then
$r^n(\varepsilon+s)>\varepsilon$ for $\varepsilon
>\frac{-r^ns}{r^n-1}$. Thus the claim follows.

\vsp

In order to approximate the other three Conradian orderings, we
first note that, arguing just as before, we have
$\preceq_\varepsilon \to \preceq_{C_2}$ as $\varepsilon\to -\infty$.
Finally, the other two Conradian orderings
$\overline{\preceq}_{C_1}$ and $\overline{\preceq}_{C_2}$ are
approximated by $\overline{\preceq}_\varepsilon$ as $\varepsilon \to
\infty$ and $\varepsilon\to -\infty $, respectively.

\vsp

Now let $\preceq_S$ be an Smirnov-type ordering and let
$\{g_1,\ldots, g_n\}$ be a set of $\preceq_S$-positive elements.

\vsp

Suppose first that $\preceq_S=\preceq_\varepsilon$, where
$\varepsilon$ is irrational. Then we have that
$\varphi(g_i)(\varepsilon)>\varepsilon$ for all $1\leq i\leq n$.
Thus, if $\varepsilon^\prime$ is such that
$\varepsilon<\varepsilon^\prime< \min\{\varphi(g_i)(\varepsilon)\}$,
$1\leq i\leq n$, then we still have that $\varphi(g_i)
(\varepsilon^\prime)
>\varepsilon^\prime$, hence $g_i\succ_{\varepsilon^\prime} id$ for
$1\leq i \leq n$. To see that
$\succ_{\varepsilon^\prime}\not=\succ_{\varepsilon}$, first notice
that $\varphi(G_1)(x)$ is dense in $\R$ for all $x\in \R$. In
particular, taking $g\in G_1$ such that $\varepsilon< \varphi(g)(0)
<\varepsilon^\prime$, we have that
$\varphi(gb^ng^{-1})(\varepsilon)=\varphi(g)
(r^n\varphi(g)^{-1}(\varepsilon))=r^n\varphi(g)^{-1}(\varepsilon)+\varphi(g)(0)$.
Since $\varphi(g)^{-1}(\varepsilon)<0$, we have that for $n$ large
enough, $gb^ng^{-1}\prec_\varepsilon id$. The same argument shows
that $gb^ng^{-1} \succ_{\varepsilon^\prime}id$. Therefore,
$\preceq_{\varepsilon^\prime}$ and $\preceq_{\varepsilon}$ are
distinct.

\vsp

The remaining case is $\preceq_S=\preceq_\varepsilon^\pm$, where
$\varepsilon$ is rational. In this case we can order the set
$\{g_1,\ldots, g_n\}$ in such a way that there is $i_0$ with
$\varphi(g_i)(\varepsilon)>\varepsilon$ for $1\leq i\leq i_0$, and
$\varphi(g_i)(\varepsilon)=\varepsilon$ for $i_0+1\leq i \leq n$.
That is, $g_i\in Stab(\varepsilon)\simeq \Z$ for $i_0+1\leq i \leq
n$. Let $\varepsilon^\prime>\varepsilon$.

\vsp

We claim that either $\varphi(g_i)(\varepsilon^\prime)>
\varepsilon^\prime$ for all $i_0+1\leq i \leq n$ or
$\varphi(g_i)(\varepsilon^\prime)< \varepsilon^\prime$ for all
$i_0+1\leq i \leq n$. Indeed, since $\varphi$ gives an affine
action, it can not be the case that a nontrivial element of $G$
fixes two points. Hence, we have that
$\varphi(g_i)(\varepsilon^\prime) \not=\varepsilon^\prime$ for each
$i_0+1\leq i \leq n$. Now, suppose for a contradiction that there
are $g_{i_0}, \; g_{i_1}$ in $Stab(\varepsilon)$ with
$g_{i_0}(\varepsilon^\prime) <\varepsilon^\prime$ and
$g_{i_1}(\varepsilon^\prime)
>\varepsilon^\prime$. Let $n,m$ in $\N$ be such that $g_{i_0}^n
=g_{i_1}^m$. Then $\varepsilon^\prime
<\varphi(g_{i_1})^m(\varepsilon^{\prime}) = \varphi(g_{i_0})^n
(\varepsilon^{\prime})<\varepsilon^\prime$, which is a
contradiction. Thus the claim follows.

\vsp

Now assume that $\varphi(g_i)(\varepsilon^\prime)>
\varepsilon^\prime,$ for all $i_0+1\leq i\leq n$. If
$\varepsilon<\varepsilon^\prime< \min\{\varphi(g_i)(\varepsilon)\}$,
with $1\leq i\leq i_0$, then $g_i\succ_{\varepsilon^\prime} id$ for
$1\leq i \leq n$, showing that $\preceq_S$ is non isolated. In the
case where $\varphi(g_i)(\varepsilon^\prime)< \varepsilon^\prime,$
for all $i_0+1\leq i\leq n$, we let $\tilde{\varepsilon}$ be such
that $\max \{\varphi(g_i)^{-1}(\varepsilon)\}< \tilde{\varepsilon}
<\varepsilon$ for $1\leq i \leq i_0$. Then we have that
$g_i\succ_{\tilde{\varepsilon}} id$ for $1\leq i \leq n$. This shows
that, in any case, $\preceq_S=\preceq_\varepsilon^\pm$ is non
isolated in $U$. $\hfill\square$

\vs

The following theorem shows that the space of left-orderings of $G$
is made up by the Smirnov-type orderings together with the Conradian
orderings. This generalizes \cite[Theorem 1.2]{rivas}.

\vsp

\begin{thm}\label{laprop}\textit{Suppose $G$ is a non Abelian group with rational series of length 2.
If $G$ is bi-orderable, then its space of left-orderings has no
isolated points. Moreover, every non-Conradian ordering is equal to
an induced, Smirnov-type ordering, arising from the affine action of
$G$ on the real line given by $\varphi$ above.}

\end{thm}

\vsp

To prove Theorem \ref{laprop}, we will use the ideas (and notation)
involved in the {\em dynamical realization of an ordering}, here
Proposition \ref{real din}.

\vs

\noindent \textit{Proof of Theorem \ref{laprop}: } First fix $a\in
G_1$ and $b\in G$ exactly as above, that is, such that
$bab^{-1}=a^r$, where $r\in \Q$, $r>1$, and $\varphi(a)(x)=x+1$,
$\varphi(b)(x)=rx$. Now let $\preceq$ be a left-ordering on $G$, and
consider its dynamical realization.  To prove Theorem \ref{laprop},
we will distinguish two cases:

\vs

\noindent \textbf{Case 1.} The element $a \in G$ is cofinal (that
is, for every $g \in G$, there are $n_1,\, n_2$ in $\mathbb{Z}$ such
that $a^{n_1}\prec g\prec a^{n_2}$).

\vs

Note that, in a Conradian ordering, $G_1$ is convex, hence $a$
cannot be cofinal. Thus, in this case we have to prove that
$\preceq$ is an Smirnov-type ordering.

\vsp

For the next two claims, recall that for any measure $\mu$ on a
measurable space $X$ and any measurable function $f: X\to X$, the
{\em push-forward measure} $f_*(\mu)$ is defined by $f_*(\mu)(A) =
\mu(f^{-1}(A))$, where $A\subseteq X$ is a measurable subset. Note
that $f_*(\mu)$ is trivial if and only if $\mu$ is trivial.
Moreover, one has $(fg)_*(\mu) = f_*(g_*(\mu))$ for all measurable
functions $f,g$.

\vsp

Similarly, the {\em push-backward measure} $f^*(\mu)$ is defined by
$f^*(\mu)(A)=\mu(f(A))$.

\vsp\vsp\vsp\vsp

\noindent{\underbar{Claim 1.}} The subgroup $G_1$ preserves a Radon
measure $\nu$ ({\em i.e.,} a measure that is finite on compact sets)
on the real line which is unique up to a scalar multiplication and
has no atoms.

\vsp\vsp\vsp

Since $a$ is cofinal and $G_1$ is rank-one Abelian, its action on
the real line is {\em free} (that is, no point is fixed by any
nontrivial element of $G_1$). By H\"older's theorem (see
\cite[Theorem 6.10]{ghys} or \cite[\S 2.2]{book}), the action of
$G_1$ is semi-conjugated to a group of translations. More precisely,
there exists a non-decreasing, continuous, surjective function $\rho
\!: \mathbb{R} \rightarrow \mathbb{R}$ such that, to each $g \in
G_1$, one may associate a translation parameter $c_g$ so that, for
all $x \in \mathbb{R}$,
$$ \rho(g(x))=\rho(x)+ c_g.$$
Now since the Lebesgue measure $Leb$ on the real line is invariant
under translations, the {\em push-backward measure} $ \nu =
\rho^*(Leb)$ is invariant by $G_1$. Since $Leb$ is a Radon measure
without atoms, this is also the case of $\nu$.

\vsp

To see the uniqueness of $\nu$ up to scalar multiple, we follow
\cite[\S 2.2.5]{book}. Given any $G_1$-invariant measure $\mu$, we
consider the associated {\em translation number homomorphism}
$\tau_\mu:G_1 \to \R$ defined by
$$ \tau_{\mu}(g)= \left\{ \begin{array}{c c} \mu([x,g(x)]) &
\text{ if $g(x) > x$, }  \\0 & \text{if } g(x)=x, \\    -\mu([g(x),
x]) & \text{ if $g(x) < x. $}
\end{array} \right.$$
One easily checks that this definition is independent of $x\in \R$,
and that the kernel of $\tau_\mu$ coincides with the set of elements
having fixed points, which in this case reduces to the identity
element of $G_1$. Now, from \cite[Proposition 2.2.38]{book}, to
prove the uniqueness of $\nu$, it is enough to show that, for any
nontrivial Radon measure $\mu$ invariant under the action of $G_1$,
$\tau_\mu(G_1)$ is dense in $\R$. But since $G_1$ is rank-one
Abelian and $G_1\not\simeq \Z$, any nontrivial homomorphism from
$G_1$ to $\R$ has a dense image. In particular $\tau_\mu(G_1)$ is
dense in $\R$. So Claim 1 follows.

\vs

\noindent{\underbar{Claim 2.}} For some $\lambda\not= 1$, we have
$b_*(\nu) = \lambda \nu$.

\vsp\vsp\vsp

Since $G_1 \lhd G$, for any $a' \in G_1$ and all measurable $A
\subset \mathbb{R}$, we must have
$$b_*(\nu)(a' (A)) = \nu(b^{-1}a'
(A))=\nu(\bar{a}(b^{-1}(A)))=\nu(b^{-1}(A))=b_*(\nu)((A))$$ for some
$\bar{a} \in G_1$. (Actually, $a' = \bar{a}^r$.) Thus $b_*(\nu)$ is
a measure that is invariant under $G_1$. The uniqueness of the
$G_1$-invariant measure up to a scalar factor yields
$b_*(\nu)=\lambda \nu$ for some $\lambda > 0$. Assume for a
contradiction that $\lambda$ equals 1. Then the whole group $G$
preserves $\nu$. In this case, the {\em translation number
homomorphism} is defined on $G$. The kernel of $\tau_{\nu}$ must
contain the commutator subgroup of $G$, and, since $a^{r-1}=[a,b]
\in [G,G]$, we have that $\tau_{\nu}(a^{r-1})=0$, hence
$\tau_{\nu}(a)=0$. Nevertheless, this is impossible, since the
kernel of $\tau_{\nu}$ coincides with the set of elements having
fixed points on the real line (see \cite[\S 2.2.5]{book}). So Claim
2 is proved.

\vsp\vsp\vsp

By Claims 1 and 2, for each $g \in G$ we have $g_*(\nu) = \lambda_g
(\nu)$ for some $\lambda_g > 0$. Moreover, $\lambda_a = 1$ and
$\lambda_b = \lambda\not=1$. Note that, since $(fg)_*(\nu)
=f_*(g_*(\nu))$, the correspondence $g\mapsto \lambda_g$ is a group
homomorphism from $G$ to $\R_+$, the group of positive real numbers
under multiplication. Since $G_1$ is in the kernel of this
homomorphism and any $g\in G$ is of the from $b^na^s$ for $n\in \Z,
\; s\in \Q$, we have that the kernel of this homomorphism is exactly
$G_1$.

\vsp

\begin{lem}\label{lema A} {\em Let $A:G\to Af_+(\R)$, $g\mapsto A_g$, be defined by
$$\label{afin}A_g(x)=\esp \frac{1}{\lambda_g} x +
\frac{sgn(g)}{\lambda_g}\, \nu([t(g^{-1}),t(id)]),$$ where
$sgn(g)=\pm 1$ is the sign of $g$ in $\preceq$. Then $A$ is an
injective homomorphism. }
\end{lem}

\noindent {\em Proof:} For $ g,h$ in $G$ both $\preceq$-positive, we
compute
\begin{eqnarray*}
A_{gh}(x)
 &=& \frac{1}{\lambda_{gh}} x +
\frac{1}{\lambda_{gh}} \nu([t((gh)^{-1}),t(id)]) \\
 &=&   \frac{1}{\lambda_{g}\lambda_h} x +\frac{1}{\lambda_{g}\lambda_h}
 \left[ (h_*\nu) ([t(g^{-1}),t(h)]) \right]   \\
 &=&   \frac{1}{\lambda_{g}\lambda_h} x +\frac{1}{\lambda_{g}\lambda_h}
 \left[ \lambda_h\nu ([t(g^{-1}),t(id)]) +\nu ([t(h^{-1}),t(id)])  \right]   \\
 &=&  \frac{1}{\lambda_{g}\lambda_h} x +\frac{1}{\lambda_{g}}\nu
 ([t(g^{-1}),t(id)])+
 \frac{1}{\lambda_g\lambda_h}\nu ([t(h^{-1}),t(id)]) \\
 &=& A_g(A_h(x)).
 \end{eqnarray*}

\vsp

The other cases can be treated analogously, thus showing that $A$ is
a group homomorphism.

\vsp

Now, assume that $A_g(x)=x$ for some nontrivial $g\in G$. Then
$\lambda_g=1$. In particular, $g\in G_1$, since the kernel of the
application $g\mapsto \lambda_g$ is $G_1$. But in this case we have
that $g$ has no fixed point, thus assuming that
$0=\lambda_g^{n-1}\nu([t(g^{-1}),t(id)]=\nu([t(g^{-n}),t(id)]$
implies that $\nu$ is the trivial (zero) measure. This contradiction
settles Lemma \ref{lema A} . $\hfill\square$

\vs

Now, for $x\in \R$, let $F(x) = sgn(x - t(id)) \cdot
\nu([t(id),x])$. (Note that $F(t(id)) = 0$.) By semi-conjugating the
dynamical realization by $F$ we (re)obtain the faithful
representation $A \!: G \to Af_+(\R)$. More precisely, for all $g
\in G$ and all $x \in \mathbb{R}$, we have \begin{equation}
\label{semiconj} F (g(x)) = A_g (F(x)).\end{equation} \noindent For
instance, if $x > t(id)$ and $g \succ id$, then
\begin{eqnarray*}
F (g(x))
 &=& \nu ([t(id),g(x)])\\
 &=& \frac{1}{\lambda_g} \nu ([t(g^{-1}),x])\\
 &=& \frac{1}{\lambda_g} \nu([t(g^{-1}),t(id)]) + \frac{1}{\lambda_g} \nu([t(id),x])\\
&=& \frac{1}{\lambda_g} F(x) + \frac{1}{\lambda_g}
\nu([t(g^{-1}),t(id)]).
\end{eqnarray*}

The action $A$ induces a (perhaps partial) left-ordering
$\preceq_A$, namely $g \succ_A id$ if and only if $A_g (0) > 0$.
Note that equation (\ref{semiconj}) implies that for every $g\in
G_1$ such that $g\succ id$, we have $A_g(0)>0 $, hence $g\succ_A
id$. Similarly, for every $f\in G$ such that $A_f(0)>0$, we have
$f\succ id$. In particular, if the orbit under $A$ of $0$ is free
(that is, for every nontrivial element $g\in G$, we have
$A_g(0)\not=0$), then (\ref{semiconj}) yields that $\preceq_A$ is
total and coincides with $\preceq$ (our original ordering).

\vsp

If the orbit of $0$ is not free (this may arise for example when $\,
t(id) \,$ does not belong to the support of $\nu$), then the
stabilizer of $0$ under the action of $A$ is isomorphic to $\Z$.
Therefore, $\preceq$ coincides with either $\preceq_A^+$ or
$\preceq_A^-$ (the definition of $\preceq_A^{\pm}$ is similar to
that of $\preceq_{\varepsilon}^{\pm}$ above).

\vsp

At this point we have that $\preceq$ can be realized as an induced
ordering from the action given by $A$. Therefore, arguing as in the
proof of Proposition \ref{C y S}, we have that $\preceq_A$, hence
$\preceq$, is non isolated.

\vsp

To show that $\preceq$ is an Smirnov-type ordering, we need to
determine all possible embeddings of $G$ into the affine group.
Recall that $bab^{-1}=a^r$, where $r=p/q>1$.

\vsp

\begin{lem} {\em Every faithful representation of $\,G$
in the affine group is given by}
$$ a\sim \left( \begin{array}{c c} 1 & \alpha \\ 0 &1 \end{array}
\right), \;\;\; b\sim \left( \begin{array}{c c} r&\beta  \\ 0&1
\end{array} \right)$$
for some $\alpha \not=0$ and $\beta \in \mathbb{R}$.
\end{lem}

\noindent \textit{Proof:} Arguing as in Lemma \ref{lema afin} one
may check that $\varphi^\prime_{a,b}:\{a,b\}\to Af_+(\R)$ defined by
$\varphi_{a,b}^\prime(a)(x)=x+\alpha$ and
$\varphi^\prime_{a,b}(b)(x)=rx+\beta$ may be (uniquely) extended to
an homomorphic embedding $\varphi_{a,b}^\prime:G\to Af_+(\R)$.
Conversely, let
$$a\sim \left( \begin{array}{c c} s & \alpha \\ 0 &1 \end{array}
\right),\;\;\; b\sim\left( \begin{array}{c c} t & \beta \\ 0 &1
\end{array} \right)$$
be a representation. Since we are dealing with
orientation-preserving affine maps, both $s$ and $t$ are positive
real numbers. Moreover, the following equality must hold:
$$a^p\sim \left( \begin{array}{c c} s^p & s^{p-1}\alpha +\ldots +s \alpha+\alpha \\ 0 &1 \end{array}
\right)= \left( \begin{array}{c c} s^q & s^{q-1}\alpha t+ s^{q-2}\alpha t +\ldots+ \alpha t - s^q\beta+\beta \\
0 & 1 \end{array} \right)\sim ba^qb^{-1}.$$ Thus $s = 1$ and $t =
p/q=r$. Finally, since the representation is faithful,
$\alpha\not=0$. $\hfill \square$

\vspace{0.35cm}

Let $\alpha,\beta$ be such that $A_a(x)=x+\alpha$ and
$A_b(x)=rx+\beta$. We claim that if the stabilizer of $0$ under $A$
is trivial --which implies in particular that $\beta \!\neq\! 0$-- ,
then $\preceq_A$ (hence $\preceq$) coincides with
$\preceq_{\varepsilon}$ if $\alpha > 0$ (resp.
$\overline{\preceq}_{\varepsilon}$ if $\alpha < 0$), where
$\varepsilon = \frac{\beta}{(r-1) \alpha} $. Indeed, if $\alpha >
0$, then for each $g = b^n a^s \in G$, $s\in \Q$, we have $A_g(0) =
r^n s\alpha +\beta \frac{r^n - 1}{r-1} $. Hence $A_g (0) > 0$ holds
if and only if
$$r^n \frac{\beta}{(r-1)\alpha } +r^n s >\frac{\beta}{(r-1)\alpha}.$$
Letting $\varepsilon = \frac{\beta}{(r-1)\alpha}$, one easily checks
that the preceding inequality is equivalent to $g\succ_\varepsilon
id$. The claim now follows.

\vsp

In the case where the stabilizer of $0$ under $A$ is isomorphic to
$\mathbb{Z}$, similar arguments to those given above show that
$\,\preceq \,$ coincides with either $\,\preceq_{\varepsilon}^+$, or
$\, \preceq_{\varepsilon}^{-}$, or $\,
\overline{\preceq}_{\varepsilon}^+ \,$, or $\,
\overline{\preceq}_{\varepsilon}^{-} \,$, where $\varepsilon$ again
equals $ \frac{\beta}{(r-1) \alpha}$.

\vspace{0.43cm}

\noindent\textbf{Case 2.} The element $a\in G$ is not cofinal.

\vsp

In this case, for the dynamical realization of $\,\preceq \,, \,$
the set of fixed points of $a$, denoted $Fix(a)$, is non-empty. We
claim that $b(Fix(a))=Fix(a)$. Indeed, let $r=p/q$, and let $x\in
Fix(a)$. We have
$$a^{p}(b(x))=a^{p}b(x)=ba^q(x)=b(x)\, .$$
Hence $a^p(b(x))=b(x)$, which implies that $a(b(x)) = b(x)$, as
asserted. Note that, since there is no global fixed point for the
dynamical realization, we must have $b(x)\not=x \,, \,$ for all $x
\in Fix(a) \,.$ Note also that, since $G_1$ is rank-one Abelian
group, $Fix(a)=Fix(G_1)$.

\vsp

Now let $x_{-1} = \inf \{t(g)  \mid g\in G_1\}$ and $x_1=\sup
\{t(g)\mid g\in G_1\}$. It is easy to see that both $x_{-1}$ and
$x_1$ are fixed points of $G_1$. Moreover, $x_{-1}$ (resp. $x_1$) is
the first fixed point of $a$ on the left (resp. right) of $t(id)$.
In particular, $b((x_{-1},x_1))\cap(x_{-1},x_1)=\emptyset$, since
otherwise there would be a fixed point inside $(x_{-1},x_1)$. Taking
the {\em reverse} ordering if necessary, we may assume $b\succ id$.
In particular, we have that $b(x_{-1})\geq x_1$.

\vsp

We now claim that $G_1$ is a convex subgroup. First note that, by
the definition of the dynamical realization, for every $g\in G$ we
have $\,t(g)=g(t(id))$. Then, it follows that for every $g\in G_1$,
$t(g) \!\in (x_{-1},x_1)$. Now let $m,s$ in $\mathbb{Z}$ and $g\in
G_1$ be such that $id \prec b^m g \prec a^s$. Then we have $t(id) <
b^m (t(g)) < t(a^s) < x_1$. Since $b(x_{-1}) \geq x_1$, this easily
yields $m = 0$, that is, $b^m g =g\in G_1$.

\vsp

We have thus proved that $G_1$ is a convex (normal) subgroup of $G$.
Since the quotient $G / G_1$ is isomorphic to $\mathbb{Z}$, an
almost direct application of Theorem \ref{teo C} shows that the
ordering $\, \preceq \,$ is Conradian. This concludes the proof of
Theorem \ref{laprop}. $\hfill\square$

\vsp

\begin{rem} It follows from Theorem \ref{laprop} and Proposition \ref{C y S}
that no left-ordering is isolated in $\mathcal{LO}(G)$. Therefore,
since any group with normal rational series is countable,
$\mathcal{LO}(G)$ is a totally disconnected Hausdorff and compact
metric space, thus homeomorphic to the Cantor set.

\end{rem}

\begin{rem} The preceding method of proof also gives a complete classification --up to
topological semiconjugacy-- of all actions of $G$ by
orientation-preserving homeomorphisms of the real line (compare
\cite{Na-sol}). In particular, all these actions come from
left-orderings on the group (compare Question 2.4 in \cite{navas}
and the comments before it). This has been recently used by Guelman
and Liousse to classify all $C^1$ actions of the solvable
Baumslag-Solitar groups on the circle \cite{GL}.

\end{rem}


\section{The general case}
\subsection{A technical proposition}


The main purpose of this section is to prove the following

\begin{prop} \label{max}\textit{ Let $G$ be a group with only finitely many $\ce$-orderings,
and let $H$ be its maximal proper convex subgroup (with respect to
any $\ce$-ordering). Then $H$ is a Tararin group, that is, a group
with only finitely many left-orderings.}
\end{prop}

Note that the existence of a maximal convex subgroup follows from
Theorem B. Note also that Proposition \ref{max} implies that no
group with only finitely many $\ce$-orderings, whose rational series
has length at least 3, is bi-orderable (see also \cite[Proposition
3.2]{rivas}).

\vsp

The proof of Proposition \ref{max} is a direct consequence of the
following

\begin{lem} \label{useful}\textit{Let $G$ be a group with only finitely many
$\ce$-orderings whose rational series has length at least three:
\begin{equation}\label{serie convexa} \{id\}=G_0\lhd G_1\lhd G_2\lhd \ldots \lhd G_n=G\, ,\;\; n\geq
3. \end{equation} Then, given $a\in G_1$ and $b\in G_i$, $ i\leq
n-1$, we have that $bab^{-1}= a^{\varepsilon}$, where
$\varepsilon=\pm 1$.}
\end{lem}

\noindent \textit{Proof:} We shall proceed by induction on $i$. For
$i=0,1$, the conclusion is obvious. Let us deal with the case $i=2$.
Let $b\in G_2$, and suppose that $bab^{-1}=a^r$, where $r\not=\pm 1$
is rational. Clearly, this implies that $b^nab^{-n}=a^{r^n}$ for all
$n\in \Z$.

\vsp

Since $G_3/G_1$ is non Abelian, there exists $c\in G_3$ such that
$cb^pc^{-1}=b^q w$, with $p\not=q$ integers and $w\in G_1$. Note
that $wa=aw$. We let $t\in \Q$ be such that $cac^{-1}=a^t$. Then we
have
$$a^{r^q}=b^qab^{-q}=b^q\,waw^{-1} b^{-q}=cb^pc^{-1} a\,
cb^{-p}c^{-1}= cb^p a^{1/t} b^{-p}c^{-1}= c
a^{\frac{r^p}{t}}c^{-1}=a^{r^p},$$ which is impossible since $r\not=
\pm 1$ and $p\not= q$. Thus the case $i=2$ is settled.

Now assume, as the induction hypothesis, that for any $w\in G_{i-1}$
we have that $waw^{-1}=a^\varepsilon$, $\varepsilon=\pm 1$. Suppose
also that there exists $b\in G_i$ such that $bab^{-1}=a^r$,
$r\not=\pm1$. As before, we have that $b^nab^{-n}=a^{r^n}$ for all
$n\in \Z$.

\vsp

Let $c\in G_{i+1}$ such that $cb^pc^{-1}=b^q w$, with $p\not=q$
integers and $w\in G_{i-1}$. Let $t\in \Q$ be such that
$cac^{-1}=a^t$. Then we have
$$a^{r^q}= b^qab^{-q}=b^q\,w\, w^{-1}aw\, w^{-1} b^{-q}=cb^pc^{-1} a^{\varepsilon}\,
cb^{-p}c^{-1}= cb^p a^{\varepsilon/t} b^{-p}c^{-1}= c
a^{\frac{\varepsilon r^p}{t}}c^{-1}=a^{\varepsilon r^p},$$ which is
impossible since $r\not= \pm 1$ and $p\not= q$ imply
$|r^p|\not=|r^q|$. This finishes the proof of Lemma \ref{useful}.
$\hfill\square$

\vs

\noindent {\em Proof of Proposition \ref{max}:} Since in any
Conradian ordering of $G$, the series of convex subgroups is
precisely the (unique) rational series associated to $G$, we have
that $H=G_{n-1}$ in (\ref{serie convexa}). So $H$ has a rational
normal series. Therefore, to prove that $H$ is a Tararin group, we
only need to check that no quotient $G_{i}/G_{i-2}$, $2\leq i\leq
n-1$, is bi-orderable.

\vsp

Now, if in (\ref{serie convexa}) we take the quotient by the normal
and convex subgroup $G_{i-2}$, Lemma \ref{useful} implies that
certain element in $G_{i-1}/G_{i-2}$ is sent into its inverse by the
action of some element in $G_{i}/G_{i-2}$. Thus $G_{i}/G_{i-2}$ is
non bi-orderable.$\hfill\square$

\begin{cor} \label{buen remark} \textit{A group $G$ having only finitely many
$\ce$-orderings, with rational series
$$\{id\}\lhd G_1\lhd\ldots\lhd G_{n-1}\lhd G_n=G,$$
is a Tararin group if and only if $G/G_{n-2}$ is a Tararin group.}
\end{cor}

\vsp


\subsection{The inductive step}


Let $G$ be a group with rational series
$$\{ id \} = G_0 \lhd G_{1} \lhd \ldots \lhd G_{n-1} \lhd G_n = G,\;\; n\geq3,$$
such that no quotient $G_i/G_{i-2}$ is Abelian. Moreover, assume $G$
is not a Tararin group. Let $\preceq$ be a left-ordering on $G$. To
show that $\preceq$ is non-isolated we will proceed by induction.
Therefore, we assume as induction hypothesis that no group with only
finitely many $\ce$-orderings, but infinitely many left-orderings,
whose rational series has length less than $n$, has isolated
left-orderings.

\vsp

The main idea of the proof is to find a convex subgroup $H$ such
that either $H$ has no isolated left-orderings or such that $H$ is
normal and $G/H$ has no isolated left-orderings. We will see that
the appropriate convex subgroup to look at is the {\em convex
closure of $G_1$} (with respect to $\preceq$), that is, the smallest
convex subgroup that contains $G_1$.

\vsp

For $x,y$ in $G$, consider the relation in $G$ given by $x\sim y $
if and only if there are $g_1,g_2$ in $G_1$ such that $g_1x\preceq y
\preceq g_2x$. We check that $\sim$ is an equivalence relation.
Clearly $x\sim x$ for all $x\in G$. If $x\sim y$ and $y\sim z$, then
there are $g_1,\, g_2,\, g_1^\prime, \,g_2^\prime$ in $G_1$ such
that $g_1x\preceq y \preceq g_2x$ and $g_1^\prime y \preceq z\preceq
g_2^\prime y$. Then $g_1^\prime g_1x\preceq z\preceq g_2^\prime g_2
x$, hence $x\sim z$. Finally, $g_1x\preceq y \preceq g_2x$ implies
$g_2^{-1} y \preceq x\preceq g_1^{-1} y$, thus $x\sim y \Rightarrow
y\sim x$.

\vsp

Now let $g,x,y$ in $G$ be such that $x\sim y$. Then $g_1x\preceq y
\preceq g_2x$, for some $g_1,g_2$ in $G_1$, hence $gg_1x\preceq
gy\preceq gg_2x$. Since $G_1$ is normal, we have that
$gg_1x=g_1^\prime gx$ and $ gg_2x=g_2^\prime g x$, for some
$g_1^\prime, g_2^\prime$ in $G_1$. Therefore, $g_1^\prime gx\preceq
gy \preceq g_2^\prime gx$, so $ gx\sim gy$. Thus, $G$ preserves the
equivalence relation $\sim$. Let $H=\{x\in G\mid x\sim id\}$.

\vs

\noindent \underbar{Claim 1.} For every $g\in G$, we have
$$gH\cap H=\left\{ \begin{array}{c l} \emptyset & \text{if $g\notin H$,}
\\ H& \text{if $g\in H$}.  \end{array}
\right. $$

\vsp

Indeed, if $g\in H$, then $g\in (gH\cap H)$. Now, since $x\sim
y\Leftrightarrow gx\sim gy$, we have that $gH=H$. Now suppose $g$ is
such that there is some $z\in gH\cap H$. Then $id \sim z\sim g$,
which implies $g\in H$. Therefore, Claim 1 follows.

\vsp

Claim 1 implies that $H$ is a convex subgroup of $G$ that contains
$G_1$. Moreover, we have

\vs

\noindent \underbar{Claim 2.} The subgroup $H$ is the convex closure
of the subgroup $G_1$.

\vsp

Indeed, let $C$ denote the convex closure of $G_1$ in $\preceq$.
Then $H$ is a convex subgroup that contains $G_1$. Thus $C\subseteq
H$.

\vsp

To show that $H\subseteq C$ we just note that, by definition, for
every $h\in H$, there are $g_1,g_2$ in $G_1$ such that $g_1\preceq h
\preceq g_2$. So $H\subseteq C$, and Claim 2 follows.

\vsp

Proceeding as in Lemma \ref{lema 1}, we conclude that there is $c\in
G$ such that $c\,G_{n-1}$ generates the quotient $G/G_{n-1}$. We
have

\vs

\noindent \underbar{Claim 3.} $H/G_1$ is either trivial or
isomorphic to $\Z$.

\vsp

By proposition \ref{max} $G_{n-1}$ is a Tararin group. Therefore, in
the restriction of $\preceq$ to $G_{n-1}$, $G_1$ is convex. So we
have that $H\cap G_{n-1}=G_1$. This means that for every $g\in
G_{n-1}\setminus G_1$, one has $gH\cap H=\emptyset$.

\vsp

Now, assume $H/G_1$ is nontrivial and let $g\in H\setminus G_1$. By
the preceding paragraph, we have that $g\notin G_{n-1}$. Therefore,
$g=c^{m_1} w_{m_1}$, for $m_1\in \Z$, $m_1\not=0$ and $w_{m_1}\in
G_{n-1}$.

\vsp

Let $m_0$ be the least positive $m\in \Z$ such that $c^mw_m\in H$,
for $w_m\in G_{n-1}$. Then, by the minimality of $m_0$, we have that
$m_1$ is a multiple of $m_0$, say $km_0=m_1$. Letting
$(c^{m_0}w_{m_0})^k=c^{m_0k}\overline{w_{m_0}}$, we have that
$(c^{m_0}w_{m_0})^{-k}c^mw_m=\overline{w_{m_0}}^{-1}w_m\in H $.
Since $\overline{w_{m_0}}^{-1}w_m\in G_{n-1}$, we have that
$\overline{w_{m_0}}^{-1}w_m\in G_1$. Therefore we conclude that
$(c^m_0 w_{m_0})^k \,G_1=c^mw_m\, G_1$, which proves our Claim 3.

\vs

We are now in position to finish the proof of the Theorem D.
According to Claim 3 above, we need to consider two cases.

\vs

\noindent \textbf{Case 1.} $H=G_1$.

\vsp

In this case, $G_1$ is a convex normal subgroup of $\preceq$ and,
since by induction hypothesis $G/G_1$ has no isolated
left-orderings, $\preceq$ is non-isolated.

\vs

\noindent \textbf{Case 2.} $H/G_1\simeq \Z$.

\vsp

In this case, $H$ has a rational series of length 2:
$$\{id\}=G_0\lhd G_1 \lhd H.$$
We let $a\in G_1$, $a\not=id$, and $h\in H$ be such that $hG_{1}$
generates $H/G_{1}$. Let $r\in \Q$ be such that $hah^{-1}=a^r$. We
have three subcases:

\vsp

\noindent \textit{Subcase 1.} $r<0$.

\vsp

Clearly, in this subcase, $H$ is non bi-orderable. Thus $H$ is a
Tararin group and $G_1$ is convex in $H$. However, as proved in
Claim 2, $H$ is the convex closure of $G_1$. Therefore, this subcase
does not arise.

\vsp

\noindent \textit{Subcase 2.} $r>0$.

\vsp

Since $r>0$, we have that $H$ is not a Tararin group, thus $H$ has
no isolated left-orderings. Therefore, $\preceq$ is non-isolated.

\vsp

\noindent \textit{Subcase 3.} $r=0$.

\vsp

In this case, $H$ is a rank-two Abelian group, so it has no isolated
orderings. Hence $\preceq$ is non-isolated.

\vsp

This finishes the proof of Theorem D.


\subsection{An illustrative example}


This subsection is aimed to illustrate the different kind of
left-orderings that may appear in a group as above. To do this, we
will consider a family of groups with eight $\ce$-orderings. We let
$G(n)=\langle a, b , c\mid bab^{-1}=a^{-1}, cbc^{-1}=b^3,
cac^{-1}=a^n\rangle$, where $n\in \Z$. It is easy to see that $G(n)$
has a rational series of length three:
$$\{id\}\lhd G_1=\langle a\rangle \lhd G_2=\langle a, b\rangle\lhd
G(n) .$$ In particular, in a Conradian ordering, $G_1$ is convex and
normal.

\vsp

Now we note that $G(n)/G_1\simeq B(1,3)$, where $B(1,3)=\langle
\beta, \gamma \mid \gamma \beta\gamma^{-1}=\beta^3\rangle$ is a
Baumslag-Solitar group, where the isomorphism is given by
$c\mapsto\gamma\,$, $\;\;b\mapsto \beta \,$, $\;\; a \mapsto id$.
Now consider the (faithful) representation $\varphi: B(1,3)\to
Homeo_+(\R)$ of $B(1,3) \simeq G(n)/G_1$ into $Homeo_+(\R)$ given by
$\varphi(\beta)(x)=x+1$ and $\varphi(\gamma)(x)=3x$. It is easy to
see that if $x\in \R$, then $Stab_{\varphi(B(1,3))}(x)$ is either
trivial or isomorphic to $\Z$.

\vsp

In particular, $Stab_{\varphi(B(1,3))}(\frac{-3k}{2})=\langle
\gamma\beta^k\rangle$, where $k\in \Z$. Thus $\langle \gamma\beta^k
\rangle$ is convex in the induced ordering from the point
$\frac{-3k}{2}$ (in the representation given by $\varphi$). Now,
using the isomorphism $G(n)/G_1 \simeq B(1,3)$, we have induced an
ordering on $G(n)/G_1$ with the property that $\langle
cb^k\,G_1\rangle$ is convex. We denote this left-ordering by
$\preceq_2$. Now, extending $\preceq_2$ by the initial Conradian
ordering on $G_1$, we have created an ordering $\preceq$ on $G(n)$
with the property that $H(n)=\langle a,cb^k\rangle$ is convex.
Moreover, we have:

\vsp

\noindent - If $n=1$ and $k=0$, then $H(n)=\langle a,c \rangle \leq
G(n)$ is convex in $\preceq$ and $ca=ac$, as in Subcase 3 above.

\vsp

\noindent - If $n\geq 2$, and $k=0$, then $H(n)=\langle a,c\rangle
\leq G(n) $ is convex in $\preceq$ and $cac^{-1}=a^2$, as in Subcase
2 above.

\vsp

\noindent - If $n\leq -1$ and $k$ is odd, then $H(n)=\langle a,
cb^k\rangle \leq G(n)$ is convex and $cb^k \, a\, b^{-k}c^{-1}=
a^{-n}$, (again) as in Subcase 2 above.


\chapter[{\small On the space of left-orderings of the free group}]
{On the space of left-orderings of the free group}

As announced in the Introduction, in this chapter we give an
explicit construction leading to a proof of the following theorem
obtained by Clay in \cite{clay 2}:

\vsp\vsp

\noindent{\bf Theorem F (Clay).}  \textit{The space of
left-orderings of the free group on two or more generators $F_n$ has
a dense orbit under the natural conjugacy action of $F_n$.}

\vsp\vsp

We must note that, though stated for $F_n$, $n\geq2$, we will only
deal with the case $n=2$. Nevertheless, our method extends in a
rather obvious way for the general case (see Remark \ref{F_2 implica
F_n}).

\vsp

As a Corollary, and following \cite{clay 2}, we next re-prove
McCleary's theorem from \cite{mccleary} asserting that the space of
left-orderings of the free group on two or more generators has no
isolated points, hence it is homeomorphic to a Cantor set. Indeed,
we have the following general

\begin{prop} \label{clay -> mccleary} {\em Suppose $G$ is an infinite, left-orderable group such that
$\mathcal{LO}(G)$ contains a dense orbit for the conjugacy action of
$G$ on $\mathcal{LO}(G)$. Then $\mathcal{LO}(G)$ contains no
isolated points.}
\end{prop}

\noindent {\em Proof:} Let $\preceq_D$ be an ordering with dense
orbit in $\mathcal{LO}(G)$. We distinguish two cases.

\vsp

\noindent \textbf{Case 1.} $\preceq_D$ is non isolated.

\vsp

In this case, since the action of $G$ on $\mathcal{LO}(G)$ is by
homeomorphism, we have that no point in the orbit of $\preceq_D$ is
isolated. In particular, no point in $\mathcal{LO}(G)$ is isolated.

\vsp

\noindent \textbf{Case 2.} $\preceq_D$ is isolated.

\vsp

If $\preceq_D$ is isolated, then its reverse ordering
$\overline{\preceq}_D$ is also isolated (recall that, for any $f\in
G$, $id\prec_D f$ if and only if $id\, \overline{\succ}_D \,f$).
This implies that $\overline{\preceq}_D \in Orb_G(\preceq_D)$.
Hence, there exists $g\in G$ such that
$g(\preceq_D)=\overline{\preceq}_D$. By the definition of the
action, this means that $gfg^{-1}\prec_D id$, for any $f\succ_D id$.
But this is impossible, since the $\preceq_D$-signs of $g$ and
$g^{-1}$ are preserved under conjugation by $g$. Thus Case 2 never
arises. $\hfill\square$


\section{Constructing a dense orbit}


Let $\preceq$ be a left-ordering on $F_2=\langle a,b\rangle$. Let
$D:F_2 \to Homeo_+(\R)$ be an homomorphic embedding with the
property that {\em there exists $x \in \R$ such that $g\succ id$ if
and only if $D(g)(x)>x$.} We call $D$ a {\em dynamical
realization-like homomorphism} for $\preceq$. The point $x$ is
called {\em reference point} for $D$.

\vsp

\begin{ex} The embedding given by any dynamical realization of any
countable left-ordered group $(G,\preceq)$ is a dynamical
realization-like homomorphism for $\preceq$ with reference point
$t_\preceq(id)$.
\end{ex}

\begin{defn} Let $B_n=\{w\in F_2=\langle a,b\rangle\mid |w|\leq n\}$,
where $|w|$ represents the length of the element $w$, be the {\em
ball} of radius $n$ in $F_2$. Given $B_n\subseteq F_2$ and a
left-ordering $\preceq$ of $F_2$, let
$$g_{(B_n,\preceq)}^-= \min_\preceq \{ w\in B_n\} ,\;\;\;\; g_{(B_n,\preceq)}^+= \max_\preceq \{ w\in B_n \}.$$

Now let $D$ be a dynamical realization-like homomorphism for
$\preceq$, with reference point $x$. Then, we will refer to the
square
$[D(g_{(B_n,\preceq)}^-)(x),D(g_{(B_n,\preceq)}^+)(x)]^2\subset
\R^2$ as the $(B_n,\preceq)$-box.
\end{defn}

\vsp

We now proceed to the construction of a nice action of $F_2$ on
$\R$. Let $\mathcal{D}=\{\preceq_1,\preceq_2,\ldots\}$ be a
countable dense subset of $\mathcal{LO}(F_2)$. Let $\mathcal{B}=
\{B_n\}_{n=1}^\infty$ be the (countable) set of all balls in $F_2$.
Let $\eta:\Z\to \mathcal{B}\times \mathcal{D}$ be a surjection, with
$\eta(k)=(B_{n_k},\preceq_{m_k})$.

\vsp

Note that, if $D$ is any dynamical realization-like homomorphism for
$\preceq$, with reference point $x$, and if $\varphi:\R\to \R$ is
any increasing continuous function, then the conjugated homomorphism
$D_\varphi$ defined by $D_\varphi(g)= \varphi D(g) \varphi^{-1}$ is
again a dynamical realization-like homomorphism for $\preceq$ but
with reference point $\varphi(x)$. Therefore, for
$\eta(k)=(B_{n_k},\preceq_{m_k})$, we may let $D_{\eta(k)}:F_2\to
Homeo_+(\R)$ be a dynamical realization-like homomorphism for
$\preceq_{m_k}$ such that:

\vsp

\noindent $(i)$ The reference point for $D_{\eta(k)}$ is $k$.

\vsp

\noindent $(ii)$ The $\eta(k)$-box coincides with the square
$[k-1/3, k+1/3]^2$.

\vs

The next lemma shows that, in the action given by $D_{\eta(k)}$, the
$\preceq_{m_k}$-signs of elements in $B_{n_k}$ are contained as part
of the information of the {\em graphs}\footnote{As usual, for $f\in
Homeo_+(\R)$, the set $\{(x,f(x))\mid x\in \R\} \subset \R^2$ is
called the graph of $f$. } of $D_{\eta(k)}(a)$ and $D_{\eta(k)}(b)$
inside the $\eta(k)$-box.

\begin{lem} \label{signo}{\em Let $\eta(k)=(B_{n_k},\preceq_{m_k})$, and let $D_{\eta(k)}$
be dynamical realization-like homomorphisms satisfying properties
$(i)$ and $(ii)$ above. Then, for every $w\in B_{n_k}$, we have that
$D_{\eta(k)}(w)(k)$ belongs to $[k-1/3, k+1/3]$, and
$D_{\eta(k)}(w)>k$ if and only $w \succ_{m_k} id$. Moreover, we have
$D_{\eta(k)}(g_{\eta(k)}^{+ })(k) = k+ 1/3$ and
$D_{\eta(k)}(g_{\eta(k)}^{- })(k) = k- 1/3$.}

\end{lem}

\noindent {\em Proof:} Since $D_{\eta(k)}$ is a dynamical
realization-like homomorphism, property ($i)$ above implies that,
for any $w\in F_2$, we have that $D_{\eta(k)}(w)(k)>k$ if and only
if $w\succ_{m_k} id$.

\vsp

The fact that $D_{\eta(k)}(g_{\eta(k)}^{+ })(k) = k+ 1/3$ and
$D_{\eta(k)}(g_{\eta(k)}^{- })(k) = k- 1/3$ is a direct consequence
of property $(ii)$ above.

\vsp

Finally, note that for every $w\in B_{n_k}$ we have $g_{\eta(k)}^{-}
\preceq_{m_k} w\preceq_{m_k} g_{\eta(k)}^{+ }$. In particular,
$D_{\eta(k)}(w)(k)\in [k-1/3, k+1/3]$. $\hfill\square$

\vsp

\begin{rem} \label{remark signo} Note that every initial segment $w_1$ of any (reduced)
word $w\in B_{n_k}$ lies again in $B_{n_k}$. Hence, the iterates of
$k$ along the word $w$ are independent of the graphs of
$D_{\eta(k)}(a)$ and $D_{\eta(k)}(b)$ outside the $\eta(k)$-box
$=[k-1/3,k+1/3]^2$. In particular, for any representation $D:F_2\to
Homeo_+(\R)$ such that the graphs of $D(a)$ and $D(b)$ coincide with
the graphs of $D_{\eta(k)}(a)$ and $D_{\eta(k)}(b)$ inside
$[k-1/3,k+1/3]^2$, respectively, the conclusion of Lemma \ref{signo}
holds for $D$ instead of $D_{\eta(k)}$.
\end{rem}

\vsp


Theorem F is a direct consequence of the following

\begin{prop} \label{the action} {\em Let $F_2=\langle a,b\rangle$.
Then there is an homomorphic embedding $D:F_2\to Homeo_+(\R)$ such
that, for each $k\in \Z$, the graphs of $D(a)$ and $D(b)$ inside
$[k-1/3,k+1/3]^2$ coincide with the graphs of $D_{\eta(k)}(a)$ and
$D_{\eta(k)}(b)$, respectively. In this action, all the integers lie
in the same orbit.}
\end{prop}

\noindent {\em Proof of Theorem F from Proposition \ref{the
action}:} Let $(x_0,x_1,\ldots )$ be a dense sequence in $\R$ such
that $x_0=0$ (note that $0$ may not have a free orbit), and let $D$
be the homomorphic embedding given by Proposition \ref{the action}.
Let $\preceq$ be the induced ordering on $F_2=\langle a,b\rangle$
from the action $D$ and the reference points $( x_0,x_1,x_2,\ldots)$
(see the comments after Theorem \ref{Cohn}). In particular, for
$g\in F_2$, we have that $D(g)(0)>0\Rightarrow g\succ id$. We claim
that $\preceq$ has a dense orbit under the natural action of $F_2$
on $\mathcal{LO}(F_2)$.

\vsp

Clearly, to prove our claim it is enough to prove that the orbit of
$\preceq$ accumulates at every $\preceq_m \in \mathcal{D}$. That is,
given $\preceq_m$ and any finite set $\{ h_1,h_2,... ,h_N \}$ such
that $id \prec_m h_j$, for $1\leq j \leq N$, we need to find  $w\in
F_2$ such that $h_j\succ_w id$ for every $1\leq j\leq N$, where, as
defined in \S \ref{acting on the space}, $g\succ_w id $ if and only
if $wgw^{-1}\succ id$.

\vsp

Let $n\in \N$ be such that $h_1,\ldots , h_N$ belongs to $B_n$. Let
$k$ be such that $\eta(k)=(B_n,\preceq_m)$. By Proposition \ref{the
action}, there is $w_k\in F_2$ such that $D(w_k)(0)=k$. Also by
Proposition \ref{the action}, the graphs of $D(a)$ and $D(b)$,
inside $[k-1/3,k+1/3]^2=\eta(k)$-box, are the same as those of
$D_{\eta(k)}(a)$ and $D_{\eta(k)}(b)$, respectively. Then, Lemma
\ref{signo} implies that for each $h_j$, $1\leq j\leq N$, we have
that $h_i\succ_m id $ if and only if $D(h_j)(k)>k$. But this is the
same as saying that $D(h_j) (D(w_k)(0))> D(w_k)(0)$, which implies
that $D(w_k^{-1})\circ D(h_j)\circ D(w_k) (0)>0$, where $\circ$ is
the composition operation. Therefore, by definition of $\preceq$, we
have that $w_k^{-1} h_j w_k \succ id$ for every $1\leq j \leq N$.
Now, by definition of the action of $F_2$ on $\mathcal{LO}(F_2)$,
this implies that $\preceq_{w^{-1}_k}$ is a left-ordering such that
$h_j\succ_{w^{-1}_k} id$. This finishes the proof of Theorem F.
$\hfill\square$

\vsp\vsp

Before proving Proposition \ref{the action} we need one more lemma.
Let $\hat{a}$ and $\hat{b}$ be two increasing continuous functions
of the real line such that, for each $k\in \Z$, the graphs of
$\hat{a}$ and $\hat{b}$ inside $[k-1/3,k+1/3]^2$ coincide with the
graphs of $D_{\eta(k)}(a)$ and $D_{\eta(k)}(b)$, respectively. We
have

\begin{lem} \label{inductivo} {\em For each $k\in \Z$, we can modify the
homeomorphisms $\hat{a}$ and $\hat{b}$ inside $[k-1/3, k+1+1/3]^2$
but outside $[k-1/3,k+1/3]^2 \cup [k+1-1/3,k+1+1/3]^2$ (see Figure
4.1) in such a way that the modified homeomorphisms, which we still
denote $\hat{a}$ and $\hat{b}$, have the following property $P$:

\vsp

There is a reduced word $w$ in the free group generated by $\{
\hat{a}, \hat{b}\}$ such that $w(k)=k+1$. Moreover, the iterates of
$k$ along $w$ remain inside $[k-1/3,k+1+1/3]$.}

\vsp


\end{lem}

\vspace{0.8cm}


\beginpicture

\setcoordinatesystem units <1cm,1cm>


\putrule from 2.5 0 to 4 0 \putrule from 4 0 to 4 -1.5 \putrule from
2.5 -1.5 to 4 -1.5 \putrule from 2.5 -1.5 to 2.5 0

\putrule from 4.5 0.5 to 6 0.5 \putrule from 4.5 0.5 to 4.5 2
\putrule from 6 0.5 to 6 2 \putrule from 4.5 2 to 6 2

\putrule from 0.5 -2.1 to 8 -2.1

\setdots

\plot 2.5 -1.5 2.5 2 / \plot 2.5 2  6 2 / \plot 6 2  6 -1.5 / \plot
6 -1.5 2.5 -1.5 /

\plot 1.8 -2.2 6.4 2.4 /

\plot 4 0  4 -2.8 / \plot 4 0  1.4  0 /

\plot 6 0.5 7 0.5 / \plot 4.5 0.5 4.5 2.5 /

\put{Figure 4.1} at 4 -3.5 \put{} at -4.2 0

\small

\put{$\eta(k)$} at 3.1 -0.7 \put{-box} at 3.4 -1.1

\put{$\eta(k+1)$} at 5.3 1.3 \put{-box} at 5.3 1

\put{$\bullet$} at  5.3 -2.1 \put{$\bullet$} at  3.3 -2.1
\put{$k+1$} at  5.3 -2.4 \put{$k$} at  3.3 -2.4

\put{m} at 2.9 1.2 \put{o} at 3.2 1.05 \put{d} at 3.4 0.9 \put{i} at
3.7 0.7 \put{f} at 4 0.5 \put{i} at 4.25 0.3 \put{c} at 4.5 0.1
\put{a} at 4.75 -0.1 \put{t} at 5 -0.2 \put{i} at 5.2 -0.3 \put{o}
at 5.4 -0.4 \put{n} at 5.6 -0.55 \put{s} at 5.8 -0.7

\endpicture


\vspace{0.8cm}

\noindent {\em Proof:} Note that, since for each $k\in \Z$ the
graphs of $\hat{a}$ and $\hat{b}$ coincide with the graphs of
$D_{\eta(k)}(a)$ and $D_{\eta(k)}(b)$, respectively, by Remark
\ref{remark signo} we have that $k, k+1/3$ and $k-1/3$ are in the
same orbit (for the action of the free group generated by $\hat{a}$
and $\hat{b}$). Therefore, to show this lemma, it is enough to show
that we can modify $\hat{a}$ and $\hat{b}$ in such a way that the
following property $P^\prime$ holds:

\vs

{\em There is a reduced word $w$ in the free group generated by $\{
\hat{a}, \hat{b}\}$ such that $w(k+1/3)=k+1-1/3$. Moreover, the
iterates of $k+1/3$ along $w$ remain inside $[k-1/3,k+1+1/3]$.}

\vs

For $h\in \{\hat{a}^{\pm 1},\, \hat{b}^{\pm 1}\}$ define $l_h=\sup\{
x\in [k-1/3,k+1/3] \mid h(x)\leq k+1/3 \}$ and $r_h= \inf \{ x \in
[k+1-1/3,k+1+1/3] \mid h(x)\geq k+1-1/3\}$. Let $x_0\in\, ]k+1/3,
k+1-1/3[$. To modify $\hat{a}$ and $\hat{b}$, we proceed as follows:

\vsp

\noindent \textbf{Case 1:} There is $h\in \{\hat{a}^{\pm1},\,
\hat{b}^{\pm1}\}$ such that $l_h< k+1/3$ and $r_h=k+1-1/3$.

\vsp

In this case, we (re)define $h$ linearly from $(l_h, h(l_h))=(l_h,
k+1/3)$ to $( k+1/3 ,x_0)$, then linearly from $( k+1/3 ,x_0)$ to
$(x_0, k+1-1/3)$, and then linearly from $(x_0,k+1-1/3)$ to
$(k+1-1/3, h(k+1-1/3))= (r_h, h(r_h))$; see Figure 4.2 (a). The
other generator, say $f$, may be extended linearly from
$(l_f,f(l_f))$ to $(r_f, f(r_f))$.

\vsp

Note that in this case we have $h(k+1/3)=x_0$ and $h(x_0)=k+1-1/3$.
This shows that $P^\prime$ holds for $w=h^2$.

\vs

We note that, for $h\in \{\hat{a}^{\pm1},\, \hat{b}^{\pm1}\}$, we
have that $l_h=k+1/3 \Leftrightarrow l_{h^{-1}}<k+1/3\;$ and $\;
r_h=k+1-1/3\Leftrightarrow r_{h^{-1}}>k+1-1/3$. Therefore, if there
is no $h$ as in Case 1, then we are in


\vsp

\noindent \textbf{Case 2:} There are  $f,h \in \{\hat{a}^{\pm1},
\hat{b}^{\pm 1} \}$ such that $l_h< k+1/3$, $\;r_h>k+1-1/3$,
$\;l_f<k+1/3$ and $r_f>k+1-1/3$.

\vsp

In this case we define $h$ linearly from $(l_h,h(l_h))$ to $(k+1/3,
x_0)$, and then linearly from $(k+1/3, x_0)$ to $(r_h, h(r_h))$. For
$f$, we define it linearly from $(l_f, f(l_f))$ to $(k+1-1/3, x_0)$,
and then linearly from $(k+1-1/3, x_0)$ to $(r_f,f(r_f))$; see
Figure 4.2 (b).

\vsp

Note that $h(k+1/3)=x_0=f(k+1-1/3)$. This shows that $P^\prime$
holds for $w=f^{-1}h$. $\hfill\square$

\vspace{0.8cm}


\beginpicture

\setcoordinatesystem units <1cm,1cm>

\plot 1.42 1.7 1.8 2.1 / \plot 1.8 2.1 2.4 2.2 / \plot -1.7 -2.4
-1.3 -1.4 / \plot -1.3 -1.4 -1 0 / \plot -1 0  0.2 1.35 / \plot 0.2
1.35 1.42 1.7 /

\plot 7.5 -2 8 -1.24 / \plot 8 -1.24 8.55 0 / \plot 8.55 0 11.5 1.27
/ \plot 11.5 1.27 12.3 2 /

\plot 7.9 -2.2 8.3 -1.24 / \plot 8.3 -1.24 10.95 0 / \plot 10.95 0
11.2 2.3 /
\putrule from -2.5 0 to 2.8 0 \putrule from 0.2 2.5 to 0.2 0
\putrule from 0.2 -2.7 to 0.2 -0.5

\putrule from -2.5 -1.35 to -1 -1.35 \putrule from -1 -1.35 to -1
-2.7

\putrule from 1.4 2.7 to 1.4 1.35 \putrule from 1.4 1.35 to 2.9 1.35

\putrule from 7 0 to 12.5 0 \putrule from 9.74 2.5 to 9.74 0
\putrule from 9.74 -2.7 to 9.74 -0.5

\putrule from 7 -1.24 to 8.55 -1.24 \putrule from 8.55 -1.24 to 8.55
-2.7

\putrule from 10.95 2.7 to 10.95 1.27 \putrule from 10.95 1.27 to
12.5 1.27

\setdots \plot -2.3 -2.7 2.7 2.7 /    \plot 7 -2.8 12.5 2.8 /

\put{Figure 4.2 (a)} at 0.3 -3.5 \put{Figure 4.2 (b)} at 9.8 -3.5
\small

\put{$x_0$} at 0.2 -0.3 \put{$\bullet$} at 0.2  0

\put{$l_h$} at -1.45 -0.3 \put{$\bullet$} at -1.3 0
\put{$k+\frac{1}{3}$} at -1.2 0.5 \put{$\bullet$} at -1 0
\put{$k+1-\frac{1}{3}$} at 1.6 -0.3 \put{$\bullet$} at 1.4 0

\put{$x_0$} at 9.78 -0.3 \put{$\bullet$} at 9.75  0 \put{$\bullet$}
at 8.05 -0.01 \put{$l_h$} at 7.9 -0.3 \put{$\bullet$} at 11.45 0
\put{$r_h$} at 11.45 0.3

\put{$k+\frac{1}{3}$} at 8.6 0.5 \put{$\bullet$} at 8.55 0
\put{$k+1-\frac{1}{3}$} at 11.3 -0.5 \put{$\bullet$} at 10.95 0

\put{$h$} at -1.5 -2.6 \put{$h$} at 2.6 2.2

\put{$h$} at 7.3 -2.1  \put{$f$} at 7.9 -2.5 \put{$h$} at 12.5 2.1
\put{$f$} at 11.25 2.6

\put{} at -3 0
\endpicture


\vspace{0.5cm}

\noindent {\em Proof of Proposition \ref{the action}:} We let
$\hat{a}$ and $\hat{b}$ be as in Lemma \ref{inductivo}. For each
$k\in \Z$ we apply inductively Lemma \ref{inductivo} to modify
$\hat{a}$ and $\hat{b}$ inside $[k-1/3, k+1+1/3]^2$. Note that
property $P$ implies that these modifications are made in such a way
that they do not overlap one with each other. In particular, for
each $k\in \Z$, the graphs of $\hat{a}$ and $\hat{b}$ coincides with
the graphs of $D_{\eta(k)}(a)$ and $D_{\eta(k)}(b)$ inside
$[k-1/3,k+1/3]^2$. Moreover, since $k$ and $k+1$ lie in the same
orbit for all $k$, we have that all the integers are in the same
orbit.

\vsp

We thus let $D:F_2\to Homeo_+(\R)$ be the homomorphism defined by
$D(a)=\hat{a}$ and $D(b)=\hat{b}$. To see that $D$ is an embedding
we just have to note that any $w\in B_{n_k}$, where
$\eta(k)=(B_{n_k},\preceq_{m_k})$, acts nontrivially at the point
$k\in \R$. Indeed, since $D(w)(k)=D_{\eta(k)}(w)(k)$, Lemma
\ref{signo} applies. This finishes the proof of Proposition \ref{the
action} $\hfill\square$

\begin{rem}\label{F_2 implica F_n} We point out that only two technical
facts were needed in the proof above. The first one is that the
partial dynamics in the $(B_n,\preceq_m)$-box contains the
information of the $\preceq_m$-sings of the elements of $B_n$. The
second is that we can glue all these boxes together in a sole action
of $F_2$ so that there is a single orbit containing all the centers
of the boxes.

\vsp

In the case of a general free group $F_n = \langle a_1, \ldots, a_n
\rangle$, $n\geq2$,the first of these facts clearly holds. The
second fact can be ensured by performing the same construction
taking $a = a_1$, and $b = a_2$, whereas the remaining generators
are extended linearly from the edge of one box to the edge of the
following box. This gives an action of $F_n$ for which there is a
single orbit containing the centers of the boxes, which allows
conclude as in the case of $F_2$.
\end{rem}



\chapter[\small Describing all bi-orderings on Thompson's group $\efe$] {Describing all bi-orderings on Thompson's group $\efe$}


In this chapter, we focus on a remarkable bi-orderable group, namely
Thompson's group $\mathrm{F}$, and we provide a complete description
of all its possible bi-orderings. This is essentially taken from
\cite{thompson}.

\vsp

Recall that $\mathrm{F}$ is the group of orientation-preserving
piecewise-linear homeomorphisms $f$ of the interval $[0,1]$ such
that:

\vspace{0.08cm}

\noindent -- the derivative of $f$ on each linearity interval is an
integer power of $2$,

\vspace{0.08cm}

\noindent -- $f$ induces a bijection of the set of dyadic rational
numbers in $[0,1]$.

\vspace{0.08cm}

\noindent For each nontrivial $f \!\in\! \mathrm{F}$ we will denote
by $x^{-}_{f}$ (resp. $x^{+}_f$) the leftmost point $x^-$ (resp. the
rightmost point $x^+$) for which $f'_+ (x^{-}) \neq 1$ (resp.
$f'_{-} (x^{+}) \neq 1$), where $f'_+$ and $f'_{-}$ stand for the
corresponding lateral derivatives. One can then immediately
visualize four different bi-orderings on (each subgroup of)
$\mathrm{F}$, namely:\\

\vspace{0.1cm}

\noindent -- the bi-ordering $\preceq_{x^{-}}^{+}$ for which \esp $f
\succ id $
\esp \esp if and only if \esp $f'_+ (x^{-}_f) > 1$,\\

\vspace{0.1cm}

\noindent -- the bi-ordering $\preceq_{x^{-}}^{-}$ for which \esp $f
\succ id $
\esp \esp if and only if \esp $f'_+ (x^{-}_f) < 1$,\\

\vspace{0.1cm}

\noindent -- the bi-ordering $\preceq_{x^{+}}^{+}$ for which \esp $f
\succ id $
\esp \esp if and only if \esp $f'_{-} (x^{+}_f) < 1$,\\

\vspace{0.1cm}

\noindent -- the bi-ordering $\preceq_{x^{+}}^{-}$ for which \esp $f
\succ id $
\esp \esp if and only if \esp $f'_{-} (x^{+}_f) > 1$.\\

\vspace{0.1cm}

\noindent Although $\mathrm{F}$ admits many more bi-orderings than
these, the case of its derived subgroup $\mathrm{F}'$ is quite
different. As discussed in the Introduction, this particular case is
related to Dlab's work \cite{dlab}. In \S \ref{biord F prime} we
show

\begin{thm} \label{teo dlab}{\em The only bi-orderings on $\mathrm{F}'$ are
$\preceq_{x^{-}}^{+}$, $\preceq_{x^{-}}^{-}$, $\preceq_{x^{+}}^{+}$
and $\preceq_{x^{+}}^{-}$.}
\end{thm}


\vspace{0.1cm}

Remark that there are also four other ``exotic" bi-orderings on $\efe$, namely:\\

\vspace{0.1cm}

\noindent -- the bi-ordering $\preceq_{0,x^{-}}^{+,-}$ for which $f
\succ id$ if and only
if either $x^{-}_f = 0$ and $f'_+(0) > 1$, or $x^{-}_f \neq 0$ and $f'_+ (x^{-}_f) < 1$,\\

\vspace{0.1cm}

\noindent -- the bi-ordering $\preceq_{0,x^{-}}^{-,+}$ for which $f
\succ id$ if and only
if either $x^{-}_f = 0$ and $f'_+(0) < 1$, or $x^{-}_f \neq 0$ and $f'_+ (x^{-}_f) > 1$,\\

\vspace{0.1cm}

\noindent -- the bi-ordering $\preceq_{1,x^{+}}^{+,-}$ for which $f
\succ id $ if and only
if either $x^{+}_f = 1$ and $f'_{-}(1) < 1$, or $x^{+}_f \neq 1$ and $f'_{-} (x^{+}_f) > 1$,\\

\vspace{0.1cm}

\noindent -- the bi-ordering $\preceq_{1,x^{+}}^{-,+}$ for which $f
\succ id $ if and only
if either $x^{+}_f = 1$ and $f'_{-}(1) > 1$, or $x^{+}_f \neq 1$ and $f'_{-} (x^{+}_f) < 1$.\\

\vspace{0.1cm}

\noindent Notice that, when restricted to $\mathrm{F}'$, the
bi-ordering $\preceq_{0,x^{-}}^{+,-}$ (resp.
$\preceq_{0,x^{-}}^{-,+}$, $\preceq_{1,x^{+}}^{+,-}$, and
$\preceq_{1,x^{+}}^{-,+}$) coincides with $\preceq_{x^{-}}^{-}$
(resp. $\preceq_{x^{-}}^{+}$, $\preceq_{x^{+}}^{-}$, and
$\preceq_{x^{+}}^{+}$). Let us denote the set of the previous eight
bi-orderings on $\efe$ by $\mathcal{BO}_{Isol} (\mathrm{F})$.

\begin{rem}
 As the reader can easily check, the bi-ordering $\preceq_{0,x^{-}}^{+,-}$
appears as the extension by $\preceq_{x^{-}}^{+}$ of the restriction
of its reverse ordering $\bar{\preceq}_{x^{-}}^{+}$ (which coincides
with $\preceq_{x^{-}}^{-}$) to the maximal proper
${\preceq_{x^{-}}^{+}}$-convex subgroup $\mathrm{F}^{max} = \{f
\!\in\! \efe \!:\ f'_+(0) = 1\}$. The bi-orderings
$\preceq_{0,x^{-}}^{-,+}$, $\preceq_{1,x^{+}}^{+,-}$, and
$\preceq_{1,x^{+}}^{-,+}$ may be obtained in the same way starting
from $\preceq_{x^{-}}^{-}$, $\preceq_{x^{+}}^{+}$, and
$\preceq_{x^{+}}^{-}$, respectively.
\end{rem}

\vsp

\vsp

There is another natural procedure for creating bi-orderings on
$\efe$. For this, recall the well-known fact that $\mathrm{F}'$
coincides with the subgroup of $\mathrm{F}$ formed by the elements
$f$ satisfying $f'_+(0) = f'_{-}(1) = 1$. Now let
$\preceq_{\mathbb{Z}^2}$ be any bi-ordering on $\mathbb{Z}^2$, and
let $\preceq_{\mathrm{F}'}$ be any bi-ordering on $\mathrm{F}'$. It
readily follows from Theorem  \ref{teo dlab}, that $\preceq_{\efe'}$
is invariant under conjugacy by elements in $\efe$. Hence, from
Corollary \ref{nice coro}, we may define a bi-ordering $\preceq$ on
$\efe$ by declaring that $f \succ id$ if and only if either $f
\notin \mathrm{F}'$ and $\big( \log_2 (f'_{+}(0)), \log_2
(f'_{-}(1))\big) \succ_{\mathbb{Z}^2} \big( 0,0 \big)$, or $f \in
\mathrm{F}'$ and $f \succ_{\mathrm{F}'} id$.

\vsp

All possible ways of ordering finite-rank Abelian groups have been
described in \cite{sikora,teh} (see Example \ref{LOZ^2} for the
description of the space of orderings of $\Z^2$). In particular,
when the rank is greater than one, the corresponding spaces of
bi-orderings are homeomorphic to the Cantor set. Since there are
only four possibilities for the bi-ordering $\preceq_{\mathrm{F}'}$,
the preceding procedure gives four natural copies (which we will
coherently denote by $\Lambda_{x^{-}}^{+}$, $\Lambda_{x^{-}}^{-}$,
$\Lambda_{x^{+}}^{+}$, and $\Lambda_{x^{+}}^{-}$) of the Cantor set
in the space of bi-orderings of $\efe$. The main result of this
chapter establishes that these bi-orderings, together with the
special eight bi-orderings previously introduced, fill out the list
of all possible bi-orderings on $\efe$.

\vsp\vsp

\noindent {\bf Theorem G .} {\em The space of bi-orderings of \esp
$\mathrm{F}$ is the disjoint union of the finite set
$\mathcal{BO}_{Isol} (\mathrm{F})$ (whose elements are isolated
bi-orderings) and the copies of the Cantor set
$\Lambda_{x^{-}}^{+}$, $\Lambda_{x^{-}}^{-}$, $\Lambda_{x^{+}}^{+}$,
and $\Lambda_{x^{+}}^{-}$.}

\vsp\vsp

The first ingredient of the proof of this result comes from the
theory of Conradian orderings. Indeed, since $\mathrm{F}$ is
finitely generated, see \cite{CFP}, every bi-ordering $\preceq$ on
it admits a maximal proper convex subgroup
$\mathrm{F}^{max}_{\preceq}$. More importantly, this subgroup may be
detected as the kernel of a nontrivial, non-decreasing group
homomorphism into $(\mathbb{R},+)$; see Theorem \ref{teo C}. Since
$\mathrm{F}'$ is simple (see for instance \cite{CFP}) and non
Abelian, it must be contained in $\mathrm{F}^{max}_{\preceq}$. The
case of coincidence is more or less transparent: the bi-ordering on
$\mathrm{F}$ is contained in one of the four canonical copies of the
Cantor set, and the corresponding bi-ordering on $\mathbb{Z}^2$ is
of {\em irrational type} (see Example \ref{LOZ^2}). The case where
$\mathrm{F}'$ is strictly contained in $\mathrm{F}^{max}_{\preceq}$
is more complicated. The bi-ordering may still be contained in one
of the four canonical copies of the Cantor set, but the
corresponding bi-ordering on $\mathbb{Z}^2$ must be of {\em rational
type} ({\em e.g.}, a lexicographic ordering). However, it may also
coincide with one of the eight special bi-orderings listed above.
Distinguishing these two possibilities is the hardest part of the
proof. For this, we strongly use the internal structure of
$\mathrm{F}$, in particular the fact that the subgroup consisting of
elements whose support is contained in a prescribed closed dyadic
interval is isomorphic to $\mathrm{F}$ itself.

\begin{rem}
In general, if $\Gamma$ is a finitely generated (nontrivial) group
endowed with a bi-ordering $\preceq$, one can easily check that the
ordering $\preceq^{*}$ obtained as the extension by $\preceq$ of
$\bar{\preceq}$ restricted to $\Gamma^{max}_{\preceq}$ is
bi-invariant. This bi-ordering (resp. its conjugate
$\bar{\preceq}_*$) is always different from $\bar{\preceq}$ (resp.
from $\preceq$), and it coincides with $\preceq$ (resp. with
$\bar{\preceq}$) if and only if the only proper $\preceq$-convex
subgroup is the trivial one; by Conrad's theorem, $\Gamma$ is
necessarily Abelian in this case. We thus conclude that every non
Abelian finitely generated bi-orderable group admits at least four
different bi-orderings. Moreover, (nontrivial) torsion-free Abelian
groups having only two bi-orderings are those of rank-one (in higher
rank one may consider lexicographic type orderings). \label{cuatro}
\end{rem}

From Section \ref{acting on the space} we have that $Out (\efe)$
could be useful for understanding $\mathcal{BO}(\efe)$.
Nevertheless, in the case of Thompson's group $\mathrm{F}$, the
action of $Out (\mathrm{F})$ on $\mathcal{BO} (\efe)$ is almost
trivial. Indeed, according to \cite{ihes}, the group $Out
(\mathrm{F})$ contains an index-two subgroup $Out_+ (\mathrm{F})$
whose elements are (equivalence classes of) conjugacies by certain
orientation-preserving homeomorphisms of the interval $[0,1]$.
Although these homeomorphisms are dyadically piecewise-affine on
$]0,1[$, the points of discontinuity of their derivatives may
accumulate at $0$ and/or $1$, but in some ``periodically coherent"
way. It turns out that the conjugacies by these homeomorphisms
preserve the derivatives of nontrivial elements $f \!\in\!
\mathrm{F}$ at the points $x_f^-$ and $x_f^+$: this is obvious when
these points are different from $0$ and $1$, and in the other case
this follows from the explicit description of $Out (\mathrm{F})$
given in \cite{ihes}. So, according to Theorem G , this implies that
the action of $Out_+ (\mathrm{F})$ on $\mathcal{BO}(\mathrm{F})$ is
trivial.

\vsp

The set \esp $Out (\mathrm{F}) \setminus Out_+ (\mathrm{F})$ \esp
corresponds to the class of the order-two automorphism $\sigma$
induced by the conjugacy by the map $x \mapsto 1-x$. One can easily
check that \begin{equation} \label{invol} (\preceq_{x^-}^+)_{\sigma}
= \preceq_{x^+}^{-}, \quad (\preceq_{x^-}^-)_{\sigma} =
\preceq_{x^+}^{+}, \quad (\preceq_{0,x^-}^{+,-})_{\sigma} =
\preceq_{1,x^+}^{-,+}, \quad \mbox{and} \quad
(\preceq_{0,x^-}^{-,+})_{\sigma} =
\preceq_{1,x^+}^{+,-}.\end{equation} Moreover, $\sigma
(\Lambda_{x^-}^{+}) = \Lambda_{x^+}^{-}$ and $\sigma
(\Lambda_{x^-}^{-}) = \Lambda_{x^+}^{+}$, and the action on the
bi-orderings of the $\mathbb{Z}^2$-fiber can be easily described. We
leave the details to the reader.


\section{Bi-orderings on $\mathrm{F}'$}
\label{biord F prime}


For each dyadic (open, half-open, or closed) interval $I$, we will
denote by $\mathrm{F}_I$ the subgroup of $\mathrm{F}$ formed by the
elements whose {\em support}\footnote{The support of an element
$f\in \efe$ is the smallest closed set containing all the points
which are not fixed by $f$.} is contained in $I$. Notice that if $I$
is closed, then $\mathrm{F}_I$ is isomorphic to $\mathrm{F}$.
Therefore, for every closed dyadic interval $I \! \subset \esp
]0,1[$, every bi-ordering $\preceq^*$ on $\mathrm{F}'$ gives rise to
a bi-ordering on $\mathrm{F} \sim \mathrm{F}_I$. Moreover, if we fix
such an $I$, then the induced bi-ordering on $\mathrm{F}_I$
completely determines $\preceq^*$ (this is due to the invariance by
conjugacy). The content of Theorem \ref{teo dlab} consists of the
assertion that only a few (namely four) bi-orderings on
$\mathrm{F}_I$ may be extended to bi-orderings on $\mathrm{F}'$. To
prove this result, we will first focus on a general property of
bi-orderings on $\mathrm{F}$.

\vsp

Let $\preceq$ be a bi-ordering on $\efe$. Since bi-invariant
orderings are Conradian and $\mathrm{F}$ is finitely generated,
Theorem \ref{teo C} provides us with a (unique up to a positive
scalar factor) non-decreasing group homomorphism $\tau_{\preceq}
\!\!: \mathrm{F} \rightarrow (\mathbb{R},+)$, called the {\em Conrad
homomorphism}, whose kernel coincides with the maximal proper
$\preceq$-convex subgroup of $\efe$. Since $\mathrm{F}'$ is a non
Abelian simple group \cite{CFP}, this homomorphism factors through
$\mathrm{F} / \mathrm{F}' \sim \mathbb{Z}^2$, where the last
isomorphism is given by \esp $f \hspace{0.05cm} \mathrm{F}' \mapsto
\big( \log_2 (f'_+(0)), \log_2 (f'_{-}(1)) \big)$. \esp Hence, we
may write (each representative of the class of) $\tau$ in the form
$$\tau_{\preceq} (f)= a \log_2 (f'_+ (0)) + b \log_2 (f'_{-}(1)).$$
A canonical representative is obtained by taking $a,b$ so that \esp
$a^2 + b^2 = 1$. \esp We will call this the {\em normalized Conrad
homomorphism} associated to $\preceq$. In many cases, we will
consider this homomorphism as defined on $\mathbb{Z}^2 \sim
\mathrm{F} / \mathrm{F}'$, so that \esp $\tau_{\preceq} \big( (m,n)
\big) = am + bn$, and we will identify $\tau_{\preceq}$ to the
ordered pair $(a,b)$.

\vsp

Now let $\preceq^*$ be a bi-ordering on $\mathrm{F}'$. Let
$I_0\subset]0,1[$ be a closed diadic interval, and consider the
induced bi-ordering on $\mathrm{F}\sim \mathrm{F}_{I_0}$, which we
will just denote by $\preceq$. Let $I\subset]0,1[$, be any other
closed diadic interval, and consider $\tau_{\preceq,I}$ the
corresponding normalized Conrad homomorphism defined on $F_I$. Since
each $\mathrm{F}_I$ is conjugate to $\mathrm{F}_{I_0}$ by an element
of $\mathrm{F}^\prime$, we have that
$\tau_{\preceq,I}=\tau_{\preceq,I_0}$ (as ordered pairs). Also, by
definition, $\tau_\preceq=\tau_{\preceq,I_0}$.

\begin{lem} {\em If \esp $\tau_{\preceq}$ corresponds to the pair
$(a,b)$, then either \esp $a \!=\! 0$ \esp or \esp $b \! = \! 0$.}
\end{lem}

\noindent{\em Proof:} Assume by contradiction that $a>0$ and $b>0$
(all the other cases are analogous). Fix $f \!\in\!
\mathrm{F}_{[1/2,3/4]}$ such that $f'_+ (1/2) > 1$ and $f'_{-} (3/4)
< 1$, and denote $I_1 = [1/4,3/4]$ and $I_2 = [1/2,7/8]$. Viewing
$f$ as an element in $\mathrm{F}_{I_1} \sim \mathrm{F}$, we have
$$\tau_{\preceq,I_1} (f) = b \log_2 \big( f'_{-} (3/4)) < 0.$$
Since Conrad's homomorphism is non-decreasing, this implies that $f$
is negative with respect to the restriction of $\preceq^*$ to
$\mathrm{F}_{I_1}$, and therefore $f \prec^* id$. Now viewing $f$ as
an element in $\mathrm{F}_{I_2} \sim \mathrm{F}$, we have
$$\tau_{\preceq,I_2} (f) = a \log_2 \big( f'_{+} (1/2)) > 0,$$
which implies that $f \succ^* id$, thus giving a contradiction.
$\hfill\square$

\vspace{0.45cm}

We may now pass to the proof of Theorem \ref{teo dlab}. Indeed,
assume that for the Conrad's homomorphism above one has $a \!>\! 0$
and $b \!=\! 0$. We claim that $\preceq^*$ then coincides with
$\preceq_{x^{-}}^+$. To show this, we need to show that a nontrivial
element $f \!\in\! \mathrm{F}'$ is positive with respect to
$\preceq^{*}$ if and only if $f'_{+} (x_f^{-}) \!>\! 1$. Now such an
$f$ may be seen as an element in $\mathrm{F}_{[x_f^{-},x_f^{+}]}$,
and viewed in this way Conrad's homomorphism yields
$$\tau_{\preceq,[x_f^{-},x_f^{+}]} (f) = a \log_2 (f'_+ (x_f^{-})).$$
Since $a \!>\! 0$, if $f'_{+} (x_f^{-}) \!>\! 1$ then the right-hand
member in this equality is positive. Since Conrad's homomorphism is
non-decreasing,  we have that $f$ is positive with respect to
$\preceq^*$. Analogously, if $f'_{+} (x_f^{-}) \!<\! 1$ then $f$ is
negative with respect to $\preceq^*$.

\vsp

Similar arguments show that the case $a \!<\! 0, \esp b \!=\! 0$
(resp. $a \!=\! 0, \esp b \!>\! 0$, and $a \!=\! 0, \esp b \!<\! 0$)
necessarily corresponds to the bi-ordering $\preceq_{x^{-}}^{-}$
(resp. $\preceq_{x^+}^{-}$, and $\preceq_{x^+}^{+}$), which
concludes the proof of Theorem \ref{teo dlab}.

\vspace{0.1cm}

\begin{qs} A bi-ordering whose positive cone is finitely
generated as a normal semigroup is completely determined by finitely
many inequalities ({\em i.e} it is isolated in the space of
bi-orderings). This makes it natural to ask whether this is the case
for the restrictions to $\efe'$ of $\preceq_{x^{-}}^{+}$,
$\preceq_{x^{-}}^{-}$, $\preceq_{x^+}^{+}$, and $\preceq_{x^+}^{-}$.
A more sophisticated question is the existence of generators \esp
$f,g$ \esp of $\efe'$ such that:

\vspace{0.1cm}

\noindent -- $f'_+ (x_f^{-}) > 1$, \esp $g'_+ (x_g^{-}) > 1$, \esp
$f'_{-} (x_f^{+}) < 1$, \esp and \esp $g'_{-} (x_{g}^{+}) > 1$,

\vspace{0.1cm}

\noindent -- $\efe' \! \setminus \! \{id\}$ is the disjoint union of
$\langle \{f,g\} \rangle^+_N$ and $\langle \{f^{-1}, g^{-1}\}
\rangle^+_{N}$,

\vspace{0.1cm}

\noindent -- $\efe' \! \setminus \! \{id\}$ is also the disjoint
union of $\langle \{f^{-1},g\} \rangle^+_N$ and $\langle \{f,
g^{-1}\} \rangle^+_{N}$.

\vspace{0.1cm}

\noindent A positive answer for the this question would immediately
imply Theorem \ref{teo dlab}. Indeed, any bi-ordering $\preceq$ on
$\efe'$ would be completely determined by the signs of $f$ and $g$.
For instance, if \esp $f \!\succ\! id$ \esp and \esp $g \!\succ\!
id$ \esp then $P_{\preceq}^+$ would necessarily contain $\langle
\{f,g\} \rangle_N^+$, and by the second property above this would
imply that $\preceq$ coincides with $\preceq_{x^{-}}^{+}$.
\label{prima}
\end{qs}


\section{Bi-orderings on $\mathrm{F}$}

\subsection{Isolated bi-orderings on $\mathrm{F}$}

Before classifying all bi-orderings on $\efe$, we will first give a
proof of the fact that the eight elements in $\mathcal{BO}_{Isol}
(\mathrm{F})$ are isolated in $\mathcal{BO} (\mathrm{F})$. As in the
case of $\mathrm{F}'$, this proof strongly uses Conrad's
homomorphism.

\vsp

We just need to consider the cases of $\preceq_{x^{-}}^{+}$ and
$\preceq_{0,x^{-}}^{+,-}$. Indeed, all the other elements in
$\mathcal{BO}_{Isol} (\mathrm{F})$ are obtained from these by the
action of the (finite Klein's) group generated by the involutions
$\preceq \hspace{0.05cm} \mapsto \bar{\preceq}$ and $\preceq
\hspace{0.05cm} \mapsto \preceq_{\sigma}$; see equation
(\ref{invol}).

\vsp

Let us first deal with $\preceq_{x^{-}}^{+}$, denoted $\preceq$ for
simplicity. Let $(\preceq_k)$ be a sequence in
$\mathcal{BO}(\mathrm{F})$ converging to $\preceq$, and let $\tau_k
\!\sim\! (a_k,b_k)$ be the normalized Conrad's homomorphism for
$\preceq_k$ (so that \esp $\tau_{k} (m,n) = a_k m + b_k n$ \esp and
\esp $a_k^2 + b_k^2 = 1$). \esp

\vspace{0.35cm}

\noindent{\underbar{Claim 1.}} For $k$ large enough, one has
$b_k\!=\!0$.

\vspace{0.1cm}

Indeed, let $f,g$ be two elements in $\mathrm{F}_{]1/2,1]}$ which
are positive with respect to $\preceq$ and such that $f'_{-}(1)=1/2$
and $g'_{-}(1)=2$. For $k$ large enough, these elements must be
positive also with respect to $\preceq_k$. Now notice that
$$\tau_k(f) = -b_k \quad \mbox{ and } \quad \tau_k (g) = b_k.$$
Thus, if $b_k \! \neq \! 0$, then either $f \prec_k id$ or $g
\prec_k id$, which is a contradiction. Therefore, $b_k \! = \! 0$
for $k$ large enough.

\vspace{0.2cm}

Let us now consider the bi-ordering $\preceq^*$ on $\mathrm{F} \sim
\mathrm{F}_{[1/2,1]}$ obtained as the restriction of $\preceq$. Let
$\tau^* \!\sim\! (a^*,b^*)$ be the corresponding normalized Conrad's
homomorphism.

\vspace{0.35cm}

\noindent{\underbar{Claim 2.}} One has $b^* \!=\! 0$.

\vspace{0.1cm}

Indeed, for the elements $f,g$ in $\mathrm{F}_{]1/2,1]}$ above, we
have
$$\tau^* (f) = -b^* \quad \mbox{ and } \quad \tau^* (g) = b^*.$$
If $b^* \!\neq\! 0$, this would imply that one of these elements is
negative with respect to $\preceq^*$, and hence with respect to
$\preceq$, which is a contradiction. Thus, $b^* \!=\! 0$.

\vspace{0.2cm}

Denote now by $\preceq_k^*$ the restriction of $\preceq_k$ to
$\mathrm{F}_{[1/2,1]}$, and let $\tau_k^* \! \sim \! (a_k^*,b_k^*)$
be the corresponding normalized Conrad's homomorphism.

\vspace{0.35cm}

\noindent{\underbar{Claim 3.}} For $k$ large enough, one has $b_k^*
\!=\! 0$.

\vspace{0.1cm}

Indeed, the sequence $(\preceq_k^*)$ clearly converges to
$\preceq^*$. Knowing also that $b^* \!=\! 0$, the proof of this
claim is similar to that of Claim 1.

\vspace{0.35cm}

\noindent{\underbar{Claim 4.}} For $k$ large enough, one has $a_k
\!>\! 0$ and $a_k^* \!>\! 0$.

\vspace{0.2cm}

Since Conrad's homomorphism is nontrivial, both $a_k$ and $a_k^*$
are nonzero. Take any $f \!\in\! \mathrm{F}$ such that $f'_+(0)
\!=\! 2$. We have $\tau_k (f) \!=\! a_k$. Hence, if $a_k < 0$, then
\esp $f \prec_k id$, \esp while \esp $f \succ id$... Analogously, if
$a_k^* < 0$, then one would have \esp $g \prec_k id$ \esp and \esp
$g \succ id$ \esp for any $g \in \mathrm{F}_{[1/2,1]}$ satisfying
\esp $g'(1/2) = 2$.

\vspace{0.35cm}

\noindent{\underbar{Claim 5.}} If $a_k$ and $a_k^*$ are positive and
$b_k$ and $b_k^*$ are zero, then $\preceq_k$ coincides with
$\preceq$.

\vspace{0.2cm}

Given $f \!\in\! \mathrm{F}$ such that $f \succ id$, we need to show
that $f$ is positive also with respect to $\preceq_k$. If $x_f^{-}
\!=\! 0$, then $f'_{+}(0) > 1$, and since $a_k \!>\! 0$, this gives
$\tau_k (f) = a_k \log_2 (f'_{+}(0)) > 0$, and thus $f \succ_k id$.

\vsp

If $x_f^{-} \neq 0$, then $f'_{+}(x_f^{-}) > 1$. In the case
$x_f=1/2$, since $a_k^* \!>\! 0$, we have that $\tau_k^{*} (f) =
a_k^* \log_2 (f'_{+}(x_f^{-}))
> 0$, and therefore one has $f \succ_k id$. In the case
$0<x_f\not= 1/2$, we can conjugate $f$ by $h\in \efe$ such that
$x_{hfh^{-1}}=1/2$. As before, we get $\tau_k^{*} (hfh^{-1})> 0$,
and therefore one still has $f\succ_k id$.

\vspace{0.35cm}

The proof for $\preceq_{0,x^{-}}^{+,-}$ is similar to the above one.
Indeed, Claims 1, 2, and 3 still hold. Concerning Claim 4, one now
has that $a_k \!>\! 0$ and $a_k^* \!<\! 0$ for $k$ large enough.
Having this in mind, one easily concludes that $\preceq_k$ coincides
with $\preceq_{0,x^{-}}^{+,-}$ for $k$ very large.


\subsection{Classifying all bi-orderings on $\mathrm{F}$}

To simplify, we will denote by $\Lambda$ the union of
$\Lambda_{x^{-}}^{+}$, $\Lambda_{x^{-}}^{-}$, $\Lambda_{x^{+}}^{+}$,
and $\Lambda_{x^{+}}^{-}$. To prove our main result, fix a
bi-ordering $\preceq$ on $\efe$, and let $\tau_{\preceq} \!\!:
\mathrm{F} \rightarrow (\mathbb{R},+)$ be the corresponding
normalized Conrad's homomorphism. Since $\tau_{\preceq} \!\sim\!
(a,b)$ is nontrivial and factors through $\mathbb{Z}^2 \! \sim \!
\mathrm{F} / \mathrm{F}'$, there are two different cases to be
considered.

\vspace{0.45cm}

\noindent{\bf Case 1.} The image $\tau_{\preceq} (\mathbb{Z}^2)$ has
rank two.

\vspace{0.15cm}

This case appears when the quotient $a/b$ is irrational. In this
case, $\preceq$ induces the bi-ordering of irrational type
$\preceq_{a/b}$ on $\mathbb{Z}^2$ viewed as $\mathrm{F} /
\mathrm{F}'$ . Indeed, for each $f \in \mathrm{F} \setminus
\mathrm{F}'$ the value of $\tau_{\preceq} (f)$ is nonzero, and hence
it is positive if and only if $f \succ id$.

The kernel of $\tau_{\preceq}$ coincides with $\mathrm{F}'$. By
Theorem \ref{teo dlab}, the restriction of $\preceq$ to
$\mathrm{F}'$ must coincide with one of the bi-orderings
$\preceq_{x^{-}}^{+}$, $\preceq_{x^{-}}^{-}$, $\preceq_{x^{+}}^{+}$,
or $\preceq_{x^{+}}^{-}$. Therefore, $\preceq$ is contained in
$\Lambda$, and the bi-ordering induced on the $\mathbb{Z}^2$-fiber
is of irrational type.

\vspace{0.45cm}

\noindent{\bf Case 2.} The image $\tau_{\preceq} (\mathbb{Z}^2)$ has
rank one.

\vspace{0.15cm}

This is the difficult case: it appears when either $a/b$ is rational
or $b \!=\! 0$. There are two sub-cases.

\vspace{0.3cm}

\noindent{\em Subcase 1.} Either $a \!=\! 0$ or $b \!=\! 0$.

\vspace{0.15cm}

Assume first that $b \!=\! 0$. Denote by $\preceq^*$ the bi-ordering
induced on $\mathrm{F}_{[1/2,1]}$, and let $\tau_{\preceq^*}
\!\sim\! (a^*,b^*)$ be its normalized Conrad's homomorphism. We
claim that either $a^{*}$ or $b^*$ is equal to zero. Indeed, suppose
for instance that $a^* \!>\! 0$ and $b^* \!>\! 0$ (all the other
cases are analogous). Let $m,n$ be integers such that \esp $n > 0$
\esp and \esp $a^* m - b^* n > 0,$ \esp and let $f$ be an element in
$\mathrm{F}_{[3/4,1]}$ such that $f'_+(3/4) = 2^m$ and $f'_{-}(1) =
2^{-n}$. Then \esp $\tau_{\preceq^*} (f) = -b^* n < 0$, \esp and
hence $f \prec id$. On the other hand, taking $h \!\in\! \mathrm{F}$
such that $h(3/4) = 1/2$, we get that $h^{-1} f h \!\in\!
\mathrm{F}_{[1/2,1]}$, and
$$\tau_{\preceq^*} (h^{-1} f h) = a^* \log_2 ((h^{-1}fh)'_+ (1/2)) +
b^* \log_2 ((h^{-1}fh)'_{-} (1)) = am - bn > 0.$$ But this implies
that $h^{-1} f h$, and hence $f$, is positive with respect to
$\preceq$, which is a contradiction.

\vspace{0.15cm}

\noindent{\em (i)}  If $a \!>\! 0$ and $a^* \! > \! 0$: We claim
that $\preceq$ coincides with $\preceq_{x^{-}}^+$ in this case.
Indeed, let $f \!\in\! \mathrm{F}$ be an element which is positive
with respect to $\preceq_{x^{-}}^+$. We need to show that $f \succ
id$. Now, since $a > 0$, if $x_f^{-} \!=\! 0$ then
$$\tau_{\preceq} (f) = a \log_2 (f'_+ (0)) > 0,$$
and hence $f \succ id$. If $x_f^{-} \! \neq \! 0$ then taking $h
\!\in\! \mathrm{F}$ such that $h (x_f^{-}) \!=\! 1/2$ we obtain that
$h^{-1} f h \in \mathrm{F}_{[1/2,1]}$, and
$$\tau_{\preceq^*} (h^{-1} f h) = a^* \log_2 ((h^{-1} f h)'(1/2)) =
a^* \log_2 (f'(x_f^{-})).$$ Since $a^* \!>\! 0$, the value of the
last expression is positive, which implies that $h^{-1} f h$, and
hence $f$, is positive with respect to $\preceq$.

\vspace{0.15cm}

\noindent{\em (ii)} If $a \!>\! 0$ and $a^{*} \! < \! 0$: Similar
arguments to those of (i) above show that $\preceq$ coincides with
$\preceq_{0,x^{-}}^{+,-}$ in this case.

\vspace{0.15cm}

\noindent{\em (iii)} If $a \!>\! 0$ and $b^{*} \!>\! 0$: We claim
that $\preceq$ belongs to $\Lambda$, and that the induced
bi-ordering on the $\mathbb{Z}^2$-fiber is the lexicographic one. To
show this, we first remark that if $f \in \efe \setminus \efe'$ is
positive, then either $f'_+ (0) \!>\! 1$, or $f'_+ (0) \!=\! 1$ and
$f'_{-}(1) \!>\! 1$. Indeed, if $f'_{+}(0) \!\neq\! 1$, then the
value of \esp $\tau_{\preceq} (f) = a \log_2 (f'_{+}(0)) \neq 0$
\esp must be positive, since Conrad's homomorphism is
non-decreasing. If $f'_{+} (0) = 1$, we take $h \! \in \! \efe$ such
that $h (1/2) \!=\! x_f^{-}$. Then $h^{-1}fh$ belongs to
$\efe_{[1/2,1]}$, and the value of
$$\tau_{\preceq^*} (h^{-1}fh) = b^{*} \log_2 ((h^{-1}fh)'_{-}(1)) = b^{*} \log_2 (f'_{-} (1)) \neq 0$$
must be positive, since $f$ (and hence $h^{-1} f h$) is a positive
element of $\efe$.

To show that $\preceq$ induces a bi-ordering on $\mathbb{Z}^2$, we
need to check that $\efe'$ is $\preceq$-convex. Let $g \!\in\!
\efe'$ and $h \!\in\! \efe$ be such that \esp $id \preceq h \preceq
g$. If $h$ was not contained in $\efe'$, then $hg^{-1}$ would be a
negative element in $\efe \setminus \efe'$. But since
$$(hg^{-1})'_{+} (0) = h'_{+}(0) \quad \mbox{ and } \quad
(hg^{-1})'_{-}(1) = h'_{-}(1),$$ this would contradict the remark
above. Therefore, $h$ belongs to $\efe'$, which shows the
$\preceq$-convexity of $\efe'$. Again, the remark above shows that
the induced bi-ordering on $\mathbb{Z}^2$ is the lexicographic one.

\vspace{0.15cm}

\noindent{\em (iv)} If $a \!>\! 0$ and $b^{*} \!<\! 0$: As in (iii)
above, $\preceq$ belongs to $\Lambda$, and the induced bi-ordering
$\preceq_{\mathbb{Z}^2}$ on the $\mathbb{Z}^2$-fiber is the one for
which $(m,n) \succ_{\mathbb{Z}^2} (0,0)$ if and only if either $m
\!>\! 0$, or $m \!=\! 0$ and $n \!<\! 0$.

\vspace{0.15cm}

\noindent{\em (v)} If $a \!<\! 0$ and $a^* \! > \! 0$: As in (i)
above, $\preceq$ coincides with $\preceq_{0,x^{-}}^{-,+}$ in this
case.

\vspace{0.15cm}

\noindent{\em (vi)} If $a \!<\! 0$ and $a^{*} \! < \! 0$: As in (i)
above, $\preceq$ coincides with $\preceq_{x^{-}}^{-}$ in this case.

\vspace{0.15cm}

\noindent{\em (vii)} If $a \!<\! 0$ and $b^{*} \!>\! 0$: As in (iii)
above, $\preceq$ belongs to $\Lambda$, and the induced bi-ordering
$\preceq_{\mathbb{Z}^2}$ on the $\mathbb{Z}^2$-fiber is the one for
which $(m,n) \succ_{\mathbb{Z}^2} (0,0)$ if and only if either $m
\!<\! 0$, or $m \!=\! 0$ and $n \!>\! 0$.

\vspace{0.15cm}

\noindent{\em (viii)} If $a \!<\! 0$ and $b^{*} \!<\! 0$: As in
(iii) above, $\preceq$ belongs to $\Lambda$, and the induced
bi-ordering $\preceq_{\mathbb{Z}^2}$ on the $\mathbb{Z}^2$-fiber is
the one for which $(m,n) \succ_{\mathbb{Z}^2} (0,0)$ if and only if
either $m \!<\! 0$, or $m \!=\! 0$ and $n \!<\! 0$.

\vspace{0.15cm}

The case $a \!=\! 0$ is analogous to the preceding one. Letting now
$\preceq^*$ be the restriction of $\preceq$ to
$\mathrm{F}_{[0,1/2]}$, for the normalized Conrad's homomorphism
$\tau_{\preceq^*} \! \sim \! (a^*,b^*)$ one may check that either
$a^* \!=\! 0$ or $b^* \! = \! 0$.

Assume that $b \!>\! 0$. In the case $b^* \!>\! 0$ (resp. $b^* \!<\!
0$), the bi-ordering $\preceq$ coincides with $\preceq_{x^+}^{-}$
(resp. $\preceq_{1,x^+}^{-,+}$). If $a^* \!>\! 0$ (resp. $a^* \!<\!
0$), then $\preceq$ corresponds to a point in $\Lambda$ whose
induced bi-ordering $\preceq_{\mathbb{Z}^2}$ on the
$\mathbb{Z}^2$-fiber is the one for which \esp $(m,n)
\succ_{\mathbb{Z}^2} (0,0)$ \esp if and only if either $n \!>\! 0$,
or $n \!=\! 0$ and $m \!>\! 0$ (resp. either $n \!>\! 0$, or $n
\!=\! 0$ and $m \!<\! 0$).

Assume now that $ b\!<\! 0$. In the case $b^* \!>\! 0$ (resp. $b^*
\!<\! 0$), the bi-ordering $\preceq$ coincides with
$\preceq_{1,x^+}^{+,-}$ (resp. $\preceq_{x^+}^{+}$). If $a^* \!>\!
0$ (resp. $a^* \!<\! 0$), then $\preceq$ corresponds to a point in
$\Lambda$ whose induced bi-ordering $\preceq_{\mathbb{Z}^2}$ on the
$\mathbb{Z}^2$-fiber is the one for which \esp $(m,n)
\succ_{\mathbb{Z}^2} (0,0)$ \esp if and only if either $n \!<\! 0$,
or $n \!=\! 0$ and $m \!>\! 0$ (resp. either $n \!<\! 0$, or $n
\!=\! 0$ and $m \!<\! 0$).

\vspace{0.3cm}

\noindent{\em Subcase 2.} Both $a$ and $b$ are nonzero.

\vspace{0.15cm}

The main issue here is to show that $\mathrm{F}'$ is necessarily
$\preceq$-convex in $\efe$. Now since $ker(\tau_{\preceq})$ is
already $\preceq$-convex in $\efe$, to prove this it suffices to
show that $\efe'$ is $\preceq$-convex in $ker (\tau_{\preceq})$.
Assume by contradiction that $f$ is a positive element in $ker
(\tau_{\preceq}) \setminus \efe'$ that is smaller than some $h \in
\efe'$. Since both $a,b$ are non zero, we have that
$f^\prime_+(0)\not=1$ and $f^\prime _-(1) \not= 1$. Suppose first
that $\preceq$ restricted to $\efe'$ coincides with either
$\preceq_{x^{-}}^{+}$ or $\preceq_{x^{-}}^{-}$, and denote by $p$
the leftmost fixed point of $f$ in $]0,1]$. We claim that $f$ is
smaller than any positive element $g \in \efe_{]0,p[}$. Indeed,
since $\preceq$ coincides with either $\preceq_{x^{-}}^{+}$ or
$\preceq_{x^{-}}^{-}$ on $\efe'$, the element $f$ is smaller than
any positive $\bar{h} \in \efe_{]0,p[}$ such that $x_{\bar{h}}^+$ is
to the left of $x_h^{-}$; taking $n \in \mathbb{Z}$ such that
$f^{-n} (x_{h}^{-})$ is to the right of $x_g^{-}$, this gives \esp
$f = f^{-n} f f^n \prec f^{-n} \bar{h} f^n \prec g.$

Now take a positive element $h_0 \in \efe_{]0,p[}$ such that for
\esp $\bar{f}=h_0 f$ \esp there is no fixed point in $]0,p[$ (it
suffices to consider a positive $h_0 \in \efe_{[
\frac{p}{4},\frac{3p}{4}]}$ whose graph is very close to the
diagonal). Then \esp $id \prec \bar{f} \prec h_0 g$ \esp for every
positive $g \!\in\! \efe_{]0,p[}$. The argument above then shows
that $\bar{f}$ is smaller than every positive element in
$\efe_{]0,p[}$. In particular, since $h_0 = \bar{f} f^{-1}$ is in
$\efe_{]0,p[}$ and is positive, this implies that \esp $\bar{f}
\prec \bar{f} f^{-1}$, \esp and hence $f \prec id$, which is a
contradiction.

If the restriction of $\preceq$ to $\efe'$ coincides with either
$\preceq_{x^{+}}^+$ or $\preceq_{x^{+}}^{-}$, one proceeds similarly
but working on the interval $[q,1]$ instead of $[0,p]$, where $q$
denotes the rightmost fixed point of $f$ in $[0,1[$. This concludes
the proof of the $\preceq$-convexity of $\mathrm{F}'$, and hence the
proof of Theorem G.

\vspace{0.1cm}

\begin{rem} Our arguments may be easily modified to show that the subgroup
$\efe_{-} \!=\! \{f \in \efe \! : f'_{+}(0) = 1 \}$ has six
different bi-orderings, namely (the restrictions of)
$\preceq_{x^{-}}^{+}$, $\preceq_{x^{-}}^{-}$, $\preceq_{x^{+}}^{+}$,
$\preceq_{x^{+}}^{-}$, $\preceq_{1,x^{+}}^{+,-}$, and
$\preceq_{1,x^{+}}^{-,+}$. An analogous statement holds for
$\efe_{+} \!=\! \{f \in \efe \! : f'_{-}(1) = 1 \}$. Finally, the
group of piecewise-affine orientation-preserving dyadic
homeomorphisms of the real line whose support is bounded from the
right (resp. from the left) admits only two bi-orderings, namely
(the natural analogues of) $\preceq_{x^{+}}^+$ and
$\preceq_{x^{+}}^{-}$ (resp. $\preceq_{x^{-}}^+$ and
$\preceq_{x^{-}}^{-}$); compare \cite{dlab}.
\end{rem}

\begin{small}

\vspace{0.1cm}


\vspace{0.37cm}

\noindent Crist\'obal Rivas\\

\noindent Dep. de Matem\'aticas, Fac. de Ciencias, Univ. de Chile\\

\noindent Las Palmeras 3425, \~Nu\~noa, Santiago, Chile\\

\noindent Email: cristobalrivas@u.uchile.cl

\end{small}

\end{document}